\documentclass[11pt]{article}
\usepackage[normalem]{ulem}
\usepackage[english]{babel}
\usepackage{latexsym,amsfonts,amssymb,amsmath,lscape,setspace,amsthm,footmisc}
\usepackage{longtable}
\usepackage{graphicx,tikz}
\usepackage{multirow}
\usepackage{color}
\usepackage{slashed,epsfig}
\usepackage{hyperref}
\newcommand{\mathsym}[1]{{}}
\newcommand{\unicode}[1]{{}}
\newtheorem{congettura}{Conjecture}[section]
\newtheorem{definizione}{Definition}[section]
\newtheorem{teorema}{Theorem}[section]

\newtheorem{defin}{Definition}[section]

\newtheorem{proofteo}{Proof}[teorema]

\newcommand{\bd}{\begin{definizione}}
\newcommand{\ed}{\end{definizione}}

\def\IC{\relax\,\hbox{$\inbar\kern-.3em{\rm C}$}}
\def\IG{\relax\,\hbox{$\inbar\kern-.3em{\rm G}$}}
\def\IB{\relax{\rm I\kern-.18em B}}
\def\ID{\relax{\rm I\kern-.18em D}}
\def\IL{\relax{\rm I\kern-.18em L}}
\def\IF{\relax{\rm I\kern-.18em F}}
\def\IH{\relax{\rm I\kern-.18em H}}
\def\II{\relax{\rm I\kern-.17em I}}
\def\IN{\relax{\rm I\kern-.18em N}}
\def\IP{\relax{\rm I\kern-.18em P}}
\def\IQ{\relax\,\hbox{$\inbar\kern-.3em{\rm Q}$}}
\def\bfzero{\relax\,\hbox{$\inbar\kern-.3em{\rm 0}$}}
\def\IK{\relax{\rm I\kern-.18em K}}
\def\IG{\relax\,\hbox{$\inbar\kern-.3em{\rm G}$}}
 \font\cmss=cmss10 \font\cmsss=cmss10 at 7pt
\def\IR{\relax{\rm I\kern-.18em R}}
\def\ZZ{\relax\ifmmode\mathchoice
{\hbox{\cmss Z\kern-.4em Z}}{\hbox{\cmss Z\kern-.4em Z}}
{\lower.9pt\hbox{\cmsss Z\kern-.4em Z}} {\lower1.2pt\hbox{\cmsss
Z\kern-.4em Z}}\else{\cmss Z\kern-.4em Z}\fi}
\def\bfone{\relax{\rm 1\kern-.35em 1}}

\def\n010{\mathrm{N^{0,1,0}}}
\def\inbar{\vrule height1.5ex width.4pt depth0pt}
\def\bfzero{\relax{\rm I\kern-.18em 0}}
\def\bfone{\relax{\rm 1\kern-.35em 1}}

\DeclareFontFamily{U}{rsf}{} \DeclareFontShape{U}{rsf}{m}{n}{
  <5> <6> rsfs5 <7> <8> <9> rsfs7 <10-> rsfs10}{}
\DeclareMathAlphabet\Scr{U}{rsf}{m}{n}

\newcommand{\C}{{\mathbb C}}
\newcommand{\R}{{\mathbb R}}

\newcommand{\Z}{{\mathbb Z}}

\newcommand{\ft}[2]{{\textstyle\frac{#1}{#2}}}
\def\tilde{\widetilde}

\def\1bar{1\hskip -.275cm -}
\def\2bar{2\hskip -.275cm -}
\def\3bar{3\hskip -.275cm -}

\newsavebox{\uuunit}
\sbox{\uuunit}
                 {\setlength{\unitlength}{0.825em}
                      \begin{picture}(0.6,0.7)
                                      \thinlines
                                      \put(0,0){\line(1,0){0.5}}
                                      \put(0.15,0){\line(0,1){0.7}}
                                      \put(0.35,0){\line(0,1){0.8}}
                                     \multiput(0.3,0.8)(-0.04,-0.02){10}{\rule{0.5pt}{0.5pt}}
                      \end {picture}}

\makeatletter \@addtoreset{equation}{section} \makeatother

\def\bfone{\relax{\rm 1\kern-.35em 1}}

\def\bfone{\relax{\rm 1\kern-.35em 1}}
\font\cmss=cmss10 \font\cmsss=cmss10 at 7pt

\def\bfone{\relax{\rm 1\kern-.35em 1}}
\def\inbar{\vrule height1.5ex width.4pt depth0pt}
\def\IC{\relax\,\hbox{$\inbar\kern-.3em{\rm C}$}}
\def\ID{\relax{\rm I\kern-.18em D}}
\def\IF{\relax{\rm I\kern-.18em F}}
\def\IH{\relax{\rm I\kern-.18em H}}
\def\II{\relax{\rm I\kern-.17em I}}
\def\IN{\relax{\rm I\kern-.18em N}}
\def\IP{\relax{\rm I\kern-.18em P}}
\def\IQ{\relax\,\hbox{$\inbar\kern-.3em{\rm Q}$}}
\def\IR{\relax{\rm I\kern-.18em R}}
\font\cmss=cmss10 \font\cmsss=cmss10 at 7pt
\def\ZZ{\relax\ifmmode\mathchoice
{\hbox{\cmss Z\kern-.4em Z}}{\hbox{\cmss Z\kern-.4em Z}}
{\lower.9pt\hbox{\cmsss Z\kern-.4em Z}} {\lower1.2pt\hbox{\cmsss
Z\kern-.4em Z}}\else{\cmss Z\kern-.4em Z}\fi}

\def\tilde{\widetilde}
\def\bar{\overline}

\def\hat{\widehat}

\def\Coe#1.#2.{{#1\over #2}}

\def\coe#1.#2.{\relax{\textstyle {#1 \over #2}}\displaystyle}

\def\to{\rightarrow}
\def\notin{\hbox{{$\in$}\kern-.51em\hbox{/}}}

\def\IE{\relax{{\rm I\kern-.18em E}}}

\def\IGam{\relax{{\rm I}\kern-.18em \Gamma}}

\def\IA{\relax{\hbox{{\rm A}\kern-.82em {\rm A}}}}

\def\R{\hat{R}}

\newtheorem{remark}{Remark}[section]

\begin{document}
\begin{titlepage}
\begin{flushright}
\today
\end{flushright}
\begin{center}
\vskip 0.2cm
{\bf  RESOLUTION \`{A}  LA KRONHEIMER OF $\mathbb{C}^3/\Gamma$  SINGULARITIES \\[5pt]
AND   THE MONGE-AMP\`{E}RE EQUATION FOR RICCI-FLAT K\"AHLER  \\[5pt]  METRICS   
IN VIEW  OF D3-BRANE SOLUTIONS OF SUPERGRAVITY }
\par\bigskip
Massimo Bianchi ${}^{a}$, {Ugo~Bruzzo${}^{\; b,f}$,
Pietro~Fr\'e${}^{\; c,e,f}$, and Dario Martelli$^{\;d,e,f}$
\par\smallskip
\vskip 4mm
\small {\frenchspacing
{${}^a$\sl\small Dipartimento di Fisica, Universit\`{a} di Roma Tor
Vergata, Rome, Italy\\
INFN -- Sezione di Roma II}\\
\emph{e-mail:} \quad {\small {\tt massimo.bianchi@roma2.infn.it}}\\[2pt]
\vspace{-1pt}
{${}^b$\sl\small SISSA, Scuola Internazionale Superiore di Studi Avanzati,    \\
Via Bonomea 265, 34136 Trieste, Italy; Departamento de Mat\'ematica, Universidade \\ Federal da Para\'iba,
Jo\~ao Pessoa, PB, Brazil;  INFN -- Sezione di
Trieste; \\  IGAP (Institute for Geometry and Physics), Trieste, Italy}\\
\emph{e-mail:} \quad {\small {\tt ugo.bruzzo@sissa.it}}\\[2pt]
\vspace{-1pt}  {${}^c$\sl\small Dipartimento di Fisica,
Universit\`{a}
di Torino, Via P. Giuria, 1, 10125 Torino, Italy}\\
\emph{e-mail:} \quad {\small {\tt pfre@unito.it}}\\[2pt]
\vspace{-1pt}  {${}^d$\sl\small  Dipartimento di Matematica ``G.
Peano'', Universit\`{a} di Torino, \\  Via Carlo Alberto, 10, 10123 Torino, Italy}\\
\emph{e-mail:} \quad {\small {\tt dario.martelli@unito.it}}\\[2pt]
\vspace{-1pt} \vspace{3pt} {{\em $^{e}$\sl\small INFN --
 Sezione di Torino,
via P. Giuria 1, \ 10125 Torino,  Italy}}
 \centerline{$^f$ \!\!\it \sl\small Arnold-Regge Center,  via P. Giuria 1,  10125 Torino,
Italy}}}

\bigskip
\begin{center} \em
In affectionate and respectful memory of \\ Boris A. Dubrovin, a master of mathematics and life
\end{center}

\begin{abstract}  \footnotesize
In this paper we analyze the  relevance of the generalized
Kronheimer construction for the gauge-gravity correspondence.
We begin with the general structure of D3-brane solutions of type
IIB supergravity on smooth manifolds $Y^\Gamma$ that are supposed to be
the crepant resolution of quotient singularities
$\mathbb{C}^3/\Gamma$ with $\Gamma$ a finite subgroup
of $SU(3)$. We emphasize that non trivial 3-form
fluxes require the existence of imaginary self-dual harmonic forms $\omega^{2,1}$. Although
excluded in the classical Kronheimer construction they may
be reintroduced by means of mass deformations. Next we concentrate
on the other essential item for the D3-brane construction, namely, the existence of  a
Ricci-flat metric on $Y^\Gamma$. We study the issue of Ricci-flat 
K\"ahler metrics on such resolutions $Y^\Gamma$, with
particular attention to the case $\Gamma=\mathbb{Z}_4$. We advance
the conjecture that on the exceptional divisor of $Y^\Gamma$ the
Kronheimer K\"ahler metric and the Ricci-flat one, that is locally
flat at infinity, coincide. The conjecture is shown to be true in
the case of the Ricci-flat metric on ${\rm tot} K_{{\mathbb WP}[112]}$ that we construct, \textit{i.e.} the total space of  the canonical bundle
of the weighted projective space ${\mathbb WP}[112]$, which is a partial
resolution of $\mathbb{C}^3/\mathbb{Z}_4$. For the full resolution
we have $Y^{\mathbb{Z}_4}=\operatorname{tot} K_{\mathbb{F}_{2}}$, where
$\mathbb{F}_2$ is the second Hizebruch surface. We try to extend the
proof of the conjecture to this case using the one-parameter
K\"ahler metric on $\mathbb{F}_2$ produced by the Kronheimer
construction as initial datum in a Monge-Amp\`{e}re (MA) equation.
We exhibit three formulations of this MA equation, one in terms of the
K\"ahler potential, the other two in terms of the symplectic
potential but with two different choices of the variables. In both
cases one can establish a series solution in powers of the variable
along the fibers of the canonical bundle. The main property of the
MA equation is that it does not impose any condition on the initial
geometry of the exceptional divisor, rather it uniquely determines
all the subsequent terms as local functionals of this initial datum.
Although a formal proof is still missing, numerical and analytical
results support the conjecture. As a by-product of our investigation
we have identified some new properties of  this type of MA equations that we
believe to be so far unknown. 
\end{abstract}
\end{center}
\vfill
\parbox{6cm}{\hrulefill} \par
\noindent \small U.B.'s research is partly supported by PRIN ``Geometry of Algebraic Varieties'' and INdAM-GNSAGA. He is a member of the VBAC group.
\end{titlepage}
\tableofcontents \noindent {}
\newpage
\section{Introduction}
\label{Introibo} We report on the advances obtained
on the following special aspect of the
gauge/gravity correspondence, within  the context of quiver
gauge--theories
\cite{Morrison:1998cs,Bianchi:2014qma,Feng:2007ur,Feng:2000mw,Feng:2001xr}:
\textit{the relevance of the generalized Kronheimer
construction}\cite{Bruzzo:2017fwj,noietmarcovaldo} for the
resolution of $\mathbb{C}^3/\Gamma$ singularities. In particular,
after an introduction about D3-brane supergravity solutions,
we consider, within this framework, the issues of the  construction
of a Ricci-flat metric on the smooth resolution $Y^\Gamma$ of
$\mathbb{C}^3/\Gamma$. We begin with the general problem of
establishing holographic dual pairs whose members are
\begin{description}
  \item[A)] a gauge theory living on a D3-brane world volume,
  \item[B)] a classical D3-brane
solution of type IIB supergravity in D=10 supergravity.
\end{description}
Gauge theories based on quiver diagrams have been extensively
studied in the literature
\cite{Morrison:1998cs,Bianchi:2014qma,Feng:2007ur,Feng:2000mw,Feng:2001xr}
in connection with the problem of establishing holographic dual
pairs as described above. Indeed the quiver diagram is a powerful
tool which encodes the data of a K\"ahler  quotient describing the
geometry of the six directions transverse to the brane. The linear
data of such a K\"ahler (or HyperK\"ahler) quotient are the menu of
the dual supersymmetric gauge theory, as they specify:
\begin{enumerate}\itemsep=-2pt
  \item the gauge group factors,
  \item the matter multiplets,
  \item the representation assignments of the latter with respect to the gauge
  group factors.
 \end{enumerate}
The possibility of testing the holographic principle
\cite{Maldacena:1997re,Kallosh:1998ph,Ferrara:1998jm,Ferrara:1998ej,Ferrara:1997dh}
and resorting to the supergravity side of the correspondence in order
 to perform, \textit{classically and in the bulk}, quantum
calculations that pertain to the boundary gauge theory is tightly
connected with the quiver approach whenever the classical brane
solution has a conformal point corresponding to a limiting geometry
of the following type:
\begin{equation}\label{adsasak}
    M_{D} \, = \, \mathrm{AdS}_{p+2} \times \mathrm{SE}^{D-p-2}
\end{equation}
In the above equation by $\mathrm{AdS}_{p+2}$ we have denoted  anti
de Sitter space in $p+2$-dimensions while $\mathrm{SE}^{D-p-2}$
stands for a Sasaki-Einstein manifold in $D-p-2$-dimensions
\cite{Fabbri:1999hw}.
\par
Within the general scope of quivers a special subclass is that of
McKay quivers that are group theoretically defined by the embedding
of a finite discrete group $\Gamma$ in  an $n$-dimensional complex unitary group
\begin{equation}\label{hukko}
    \Gamma \hookrightarrow \mathrm{SU}(n)
\end{equation}
and are associated with the resolution of $\mathbb{C}^n/\Gamma$
quotient singularities by means of a Kronheimer-like construction
\cite{kro1,kro2,mango}.
\par
The case $n=2$ corresponds to the HyperK\"ahler quotient
construction of ALE-manifolds as the resolution of the
$\mathbb{C}^2/\Gamma$ singularities, the discrete group $\Gamma$
being a finite Kleinian subgroup of $\mathrm{SU(2)}$, as given by the
ADE classification\footnote{For a recent review of these matters see
chapter 8 of \cite{advancio}.}.
\par
The case $n=3$ was the target of many interesting and robust
mathematical developments starting from the middle of the nineties
up to the present day
\cite{itoriddo,62,Kin94,crawthesis,CrawIshii,SardoInfirri:1994is,SardoInfirri:1996ga,SardoInfirri:1996gb,degeratu}.
The main and most intriguing result in this context, which
corresponds to a generalization of the Kronheimer construction and
of the McKay correspondence, is the group theoretical prediction of
the cohomology groups
$\mathrm{H}^{(p,q)}\left(Y^\Gamma_{[3]}\right)$ of the crepant
smooth resolution $Y^\Gamma_{[3]}$ of the quotient singularity
$\mathbb{C}^3/\Gamma$.
Specifically, the main output of the generalized Kronheimer
construction for the crepant resolution of a singularity
$\mathbb{C}^3/\Gamma$ is a blowdown morphism:
\begin{equation}\label{caspiterina}
    \mathrm{BD}\, : \quad Y^\Gamma_{[3]} \, \longrightarrow \,
    \frac{\mathbb{C}^3}{\Gamma}
\end{equation}
where $Y^\Gamma_{[3]}$ is a noncompact smooth three-fold with
trivial canonical bundle. On such a complex three-fold 
 a Ricci-flat K\"ahler metric
\begin{equation}\label{carigno}
   \text{ds}^2_{\mathrm{RFK}}(Y^\Gamma_{[3]})\, = \,
    \mathrm{\mathbf{g}}^{\mathrm{RFK}}_{\alpha\beta^\star} \,
    dy^\alpha \otimes dy^{\beta^\star}
\end{equation}
with asymptotically conical boundary conditions (Quasi-ALE) is guaranteed to exist (see e.g.~ \cite{Joyce-QALE}, Thm.~3.3),
{\it although it is not necessarily the one obtained by means of the
K\"ahler quotient.} According to the theorem proved by Ito-Reid
\cite{itoriddo,crawthesis,CrawIshii} and based on the concept of age
grading\footnote{For a recent review of these matters within a
general framework of applications to brane gauge theories see
\cite{Bruzzo:2017fwj,noietmarcovaldo}.}, the homology cycles of
$Y_{[3]}^\Gamma$ are all algebraic and its non vanishing cohomology
groups are all even and of type $\mathrm{H}^{(q,q)}$. We actually
have a correspondence between the cohomology classes of type $(q,q)$
and the discrete group conjugacy classes with age-grading $q$,
encoded in the statement:
\begin{eqnarray}\label{vecchioni}
    \mbox{dim} \, \mathrm{H}^{1,1}\left(Y^\Gamma_{[3]}\right)& =
    & \# \, \mbox{ of junior conjugacy classes in $\Gamma$;}\nonumber\\
\mbox{dim} \, \mathrm{H}^{2,2}\left(Y^\Gamma_{[3]}\right)& =
    & \# \, \mbox{ of senior conjugacy classes in $\Gamma$;}
    \nonumber\\
    && \mbox{all other cohomology groups are trivial}
\end{eqnarray}
The absence of harmonic forms of type $(2,1)$ implies that the
three-folds $Y^\Gamma_{[3]}$ admit no infinitesimal deformations of
their complex structure and is also a serious obstacle, as we
discuss in section \ref{cimice} to the construction of supergravity
D3-brane solutions based on $Y^\Gamma_{[3]}$ that have transverse
three-form fluxes.
\par
 There is however a possible way out that is provided by the
existence of mass-deformations. This is the main point of another
line of investigation which we hope to report on soon. If the McKay quiver diagram has certain properties,
the superpotential $\mathcal{W}(\Phi)$ on the gauge-theory side of
the correspondence can be deformed by well defined mass-terms and,
after (gaussian) integration of the massive fields, the McKay quiver
is remodeled into a new non-McKay quiver associated with the
K\"ahler or HyperK\"ahler quotient description of smooth K\"ahler
manifolds, like the resolved conifold, that admit harmonic
$(2,1)$-forms and sustain adequate D3-brane solutions. 
\par
On the basis of the above remarks we can spell-out the scope of the
present  paper in the following way. The embedding (\ref{hukko})
determines in a unique way a McKay quiver diagram which determines:
\begin{enumerate} \itemsep=-2pt
  \item the gauge group $\mathcal{F}_\Gamma$,
  \item the matter field content $\Phi^I$ of the gauge theory,
  \item the representation assignments of all the matter fields $\Phi^I$,
  \item the possible (mass)-deformations of the superpotential
  $\mathcal{W}(\phi)$, 
  \item the Ricci-flat metric on $Y$ can be inferred, by means of the Monge-Amp\`{e}re equation, from the
  K\"ahler metric on the exceptional compact divisor (in those cases where it exists) in the
  resolution of $\mathbb{C}^3/\Gamma$, which, on its turn, is determined by the
  McKay quiver through the Kronheimer construction.
\end{enumerate}
In relation with point 4) of the above list, to be discussed in
a future paper, for the case $\mathbb{C}^3/\mathbb{Z}_4$ we
anticipate the following. By means of gaussian integration we get a
new quiver diagram that is not directly associated with a discrete
group, yet it follows from the McKay quiver of $\Gamma$ in a unique
way. The group theoretical approach allows us to identify
deformations of the superpotential and introduce new directions in
the moduli space of the crepant resolution. In this sense, we go beyond the
Ito-Reid theorem. Both physically and mathematically this is quite
interesting and provides a new viewpoint on several results, some of
them well known in the literature. Most of the latter are based on
cyclic groups $\Gamma$ and rely on the powerful weapons of toric
geometry. Yet the generalized Kronheimer construction applies also
to non abelian groups $\Gamma\subset \mathrm{SU(3)}$ and so do the
cohomological theorems proved by Ito-Reid, Ishii and Craw. Hence
available mass-deformations are encoded also in the McKay quivers of
non abelian groups $\Gamma$ and one might explore the geometry of
the transverse manifolds emerging in these cases.
\par
In relation with point 5) of the above list, fully treated in this
paper for the same case  $\mathbb{C}^3/\mathbb{Z}_4$, we stress
that, although the Kronheimer metric on $Y^\Gamma_{[3]}$  is not Ricci-flat, yet
its restriction to the exceptional divisor provides the appropriate
starting point for an iterative solution of the Monge Amp\`{e}re
equation which determines the Ricci-flat metric.
\par
In view of the above considerations we can conclude that the McKay
quiver diagram does indeed provide a determination of both sides of
a D3-brane dual pair, the gauge theory side and the supergravity
side.
\par
In this  paper we  focus on two paradigmatic examples, namely
$\mathbb{C}^3/\mathbb{Z}_3$ (with $\mathbb{Z}_3$ diagonally embedded
in $\mathrm{SU(3)}$) and  $\mathbb{C}^3/\mathbb{Z}_4$. The latter
case was studied in depth in a recent publication
\cite{noietmarcovaldo}. Relying on those results here we concentrate
on the issue of the Ricci-flat K\"ahler metric.
\par
While in the case of HyperK\"ahler quotients (yielding
$\mathcal{N}=2$ gauge theories and corresponding to the original
Kronheimer construction of $\mathbb{C}^2/\Gamma$ resolutions) the
\textit{Kronheimer metric} is automatically Ricci-flat, in the case
of K\"ahler quotients and of the generalized Kronheimer construction
of $\mathbb{C}^3/\Gamma$ resolutions, the \textit{Kronheimer metric}
is not Ricci-flat and one needs to resort to different techniques in
order to find a Ricci-flat metric on the same three-fold
$Y^\Gamma_{[3]}$ that is algebraically determined by the
construction.
\par
The fascinating scenario that emerges from our combined analytical and 
numerical results is summarized in the following discussion.
\par
From the point of view of complex algebraic geometry the resolved
variety $Y^\Gamma_{[3]}$ is in many cases and, in particular in
those here analyzed, the total space of a line-bundle over a compact
complex two-fold, which coincides with the exceptional divisor
$\mathcal{ED}$ of the resolution of singularities:
\begin{equation}\label{strattone1}
    Y^\Gamma_{[3]} \,\stackrel{\pi}{\longrightarrow} \,
    \mathcal{ED}_{[2]} \quad ; \quad \forall \, p \, \in
    \,\mathcal{ED}_{[2]} \quad : \quad \pi^{-1}(p) \, \sim \, \mathbb{C}
\end{equation}
In the paradigmatic example, recently studied in
\cite{noietmarcovaldo}, of the resolution \`{a} la Kronheimer of the
$\mathbb{C}^3/\mathbb{Z}_4$ singularity, $\mathcal{ED}$ is indeed
the compact component of the exceptional divisor emerging from the
blow-up of the singular point in the origin of $\mathbb{C}^3$ and it
happens to be the second Hirzebruch surface $\mathbb{F}_2$. Other
cases are possible.
\par
Hirzebruch surfaces are $\mathbb{P}^1$ bundles over $\mathbb{P}^1$,
so that
\begin{equation}\label{strattone2}
    \mathcal{ED}_{[2]} \,\stackrel{\tilde{\pi}}{\longrightarrow} \,
    \mathbb{P}^1 \quad ; \quad \forall \, p \, \in
    \,\mathbb{P}^1 \quad : \quad \tilde{\pi}^{-1}(p) \, \sim \, \mathbb{P}^1
\end{equation}
This double fibration is illustrated in a pictorial fashion in
fig.\ref{Yconcept}.
\par
Given this hierarchical  structure, the sought for Ricci-flat metric
is constrained to possess the following continuous isometry group:
\begin{equation}\label{gisometry}
\mathrm{G}_{iso} \, = \, \mathrm{SU(2) \times U(1)_{v} \times U(1)_{w}}
\end{equation}
whose holomorphic algebraic action on the three coordinates $u,v,w$
is described later in eq.~(\ref{ciabattabuona}). The chosen isometry
group implies that the sought for Ricci-flat metric is toric, as
each of the three complex coordinates is acted on by an independent
$\mathrm{U(1)}$-isometry. Furthermore the enhancement of one of the
$\mathrm{U(1)}$'s to $\mathrm{SU(2)}$ guarantees that either the
K\"ahler potential $\mathcal{K}$ in the standard complex formulation
of K\"ahler geometry, or the symplectic potential $\mathcal{G}$, the
Legendre transform of the former appearing in the available
symplectic formalism  \cite{abreu},
are functions only of two invariant real variables (see section
\ref{sympatia} and \ref{kelleropotto}). Assuming that we possess
either one of these two real functions for the Ricci-flat metric
\footnote{For conventions see once again sections \ref{kelleropotto}
and \ref{sympatia}.} :
\begin{equation}\label{pinacoteca}
\mathcal{K}_{\text{Ricci-flat}}(\varpi,\mathfrak{f} )\quad \text{or}
\quad \mathcal{G}_{\text{Ricci-flat}}(\mathfrak{v},\mathfrak{w})
\end{equation}
we can reduce the corresponding geometry to that of the exceptional
divisor by setting a section of the $Y^\Gamma_{[3]}$ bundle to zero
as:
\begin{equation}\label{definiodivisore}
    w \, = \, 0 \quad \Leftrightarrow \quad \mathfrak{f}=0 \,\, , \,
    \, \mathfrak{w} \, = \, \ft 32
\end{equation}
The fascinating scenario we have alluded to some lines above is
encoded in the following:
\begin{figure}
\centering
\includegraphics[width=12cm]{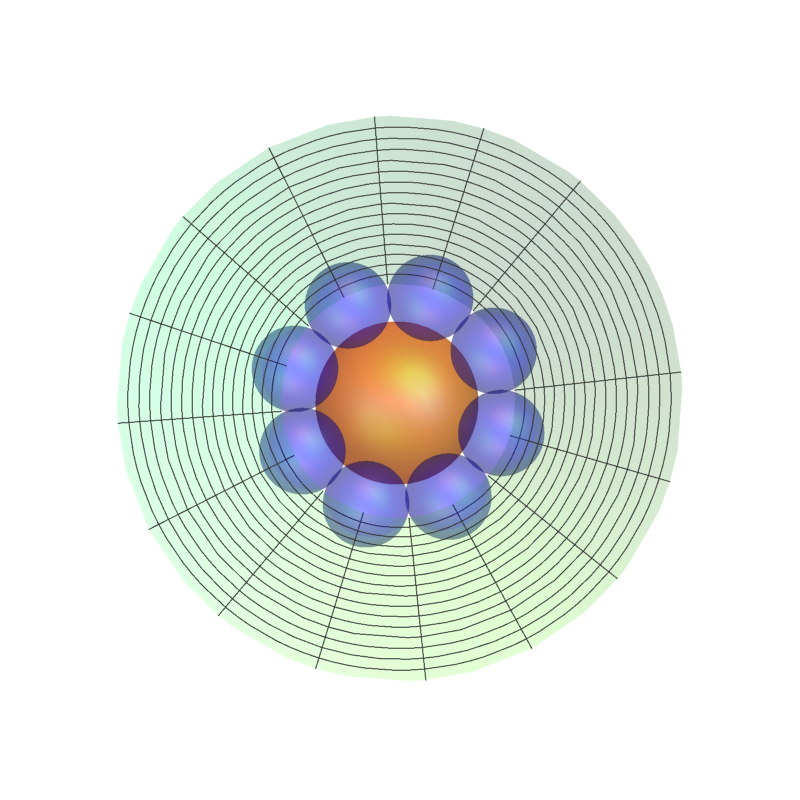}
\caption{\label{Yconcept} A conceptual picture of the resolved
three-fold $Y^\Gamma_{[3]}$ displaying its double fribration
structure. The orange sphere in the middle symbolizes the base
manifold of the bundle $\mathcal{ED}_{[2]}$. A dense complex
coordinate patch for this $\mathbb{P}^1$ is named $u$ in the main
body of the article. The blueish spheres around the orange one
symbolize the $\mathbb{P}^1$ fibers of $\mathcal{ED}_{[2]}$. A dense
complex coordinate patch for these fibers is named $v$ in the main
body of the article. Finally the greenish rays enveloping the
divisor $\mathcal{ED}_{[2]}$ symbolize the noncompact fibers of the
bundle $Y^\Gamma_{[3]}$. A dense coordinate patch for these fibers
is named $w$ in the main body of the article.}
\end{figure}
\begin{congettura}\label{cognisco} The Kronheimer K\"ahler metric
$\mathrm{ds}^2_{Kro}[Y^\Gamma_{[3]}]$ on the line bundle
(\ref{strattone1}) and the Ricci-flat one
$\mathrm{ds}^2_{Ricflat}[Y^\Gamma_{[3]}]$ on the same manifold, that
has the same isometries and is asymptotically locally
flat\footnote{More precisely, this metric is Quasi-ALE in the sense
of \cite{Joyce-QALE}.}
 are different, yet they 
 coincide on the exceptional divisor $\mathcal{ED}$.
\end{congettura}
The basic argument in favor of this conjecture is provided by an in
depth analysis of a particular orthotoric metric that we construct in this paper and
that is shown to describe the Ricci-flat
metric on a  degenerate limit of three-fold $Y^\Gamma_{[3]}$, as  described in
\cite{noietmarcovaldo}. This is a partial resolution of the
$\mathbb{C}^3/\mathbb{Z}_4$ singularity and it occurs when the stability
parameters (Fayet-Iliopolous parameters in the physics jargon) are restricted to be on the
unique type III wall\footnote{According to the terminology in
\cite{CrawIshii}, a wall in the space of stability 
parameters is of type III when it corresponds to a degeneration
which contracts divisors to curves. In this case the noncompact
component $\mathbb P^1\times \mathbb C$ of the exceptional divisor
shrinks to  $\mathbb C$.} appearing in the chamber structure
associated with the generalized Kronheimer construction for this
McKay quiver. From the algebraic geometry viewpoint, this variety $Y_{[3]}$ is the
total space of the canonical bundle over the weighted projective
space ${\mathbb WP}[112]$:
\begin{equation}\label{prunetto}
    Y_{[3]}\, = \, {\rm tot} K_{{\mathbb WP}[112]}
\end{equation}
and its exceptional divisor is ${\mathbb WP}[112]$. We show  that
the K\"ahler metric induced on ${\mathbb WP}[112]$ by our new 
Ricci-flat orthotoric metric is precisely identical with that
obtained from the Kronheimer construction once reduced to the
divisor.
\par
The various inspections of this known case within the framework of
different formalisms and using different coordinate patches provided
us with the means to make conjecture \ref{cognisco} more robust. The
main tool at our disposal is provided by the Monge-Amp\`{e}re (MA)
equation for Ricci-flatness of which we develop two versions, one in
terms of the K\"ahler potential $\mathcal{K}(\varpi,\mathfrak{f} )$
(see section \ref{kelleropotto}) and one in terms of the symplectic
potential\footnote{See section \ref{sympatia} for the definition of
the real variables $\mathfrak{v},\mathfrak{w}$.}
$\mathcal{G}(\mathfrak{v},\mathfrak{w})$(see section
\ref{sympotto}). In both cases we showed that the potential can be
developed in power series of the invariant associated with the non
compact fibers (either $\mathfrak{f}$ or $\mathfrak{w}-\ft 32$) and
that the MA equation imposes no restriction on the 0-th order
potentials $\mathcal{K}_0(\varpi)$ or $\mathcal{G}_0(\mathfrak{v})$,
namely on the geometry chosen for the exceptional divisor. Rather,
dealing carefully with the boundary conditions, we discovered that
in both cases the MA equation completely determines  all the other
terms once $\mathcal{K}_0(\varpi)$ or $\mathcal{G}_0(\mathfrak{v})$
are given. Hence we can start with $\mathcal{K}_0^{Kro}(\varpi)$ or
$\mathcal{G}_0^{Kro}(\mathfrak{v})$ as they are determined by the
Kronheimer construction and going through the power series treatment
of the MA equation we can construct  a corresponding Ricci-flat
metric.
\par
The only question which remains open is whether this Ricci-flat
metric is asymptotically locally flat. In the case of 
${\rm tot} K_{{\mathbb WP}[112]}$ it is. This supports the
conjecture. In order to transform the conjecture into a theorem one
should first resum the series and study the metric at large
distances.
\par
In this respect our study of the symplectic potential produced
encouraging results. First of all we were able to construct an
explicit form
$\mathcal{G}_{{\mathbb WP}[112]}(\mathfrak{v},\mathfrak{w})$ of such
potential for the orthotoric case.  The
function
$\mathcal{G}_{{\mathbb WP}[112]}(\mathfrak{v},\mathfrak{w})$, which
is relatively simply written in terms of elementary transcendental
functions, satisfies the MA equation and can be expanded in series
of $(\mathfrak{w}-\ft 32)$. The remarkably similar behavior of the
series truncations of the exact solution corresponding to ${\rm tot} K_{{\mathbb WP}[112]}$ with the same truncations
of the series determined by the MA equation for the smooth case
${\rm tot} K_{\mathbb{F}_{2}}$ suggests that also in
the latter case there exists a summation of the series to some
simple deformation of the function
$\mathcal{G}_{{\mathbb WP}[112]}(\mathfrak{v},\mathfrak{w})$.
\par
We postpone to future publications further attempts to sum the
series solution and prove, if possible, our conjecture.
\section{D3-brane supergravity solutions on resolved $\mathbb{C}^3/\Gamma$
singularities}\label{cimice}
An apparently  general property of the $Y_{[3]}^\Gamma$ manifolds that
emerge from the crepant resolution construction, at least when
$\Gamma$ is abelian and cyclic is the following. The non-compact
$Y_{[3]}^\Gamma$ corresponds to the total space of some line-bundle over a
complex two-dimensional compact base manifold $\mathcal{M}_2$:
\begin{equation}\label{sagnalotto}
    Y_{[3]}^\Gamma \, \stackrel{\pi}{\longrightarrow} \,\mathcal{M}_2
\end{equation}
According with this structure we name $u,v,w$ the three complex
coordinates of $Y_{[3]}^\Gamma$, $u,v$ being the coordinates of the base
manifold $\mathcal{M}_2$ and $w$ being the coordinate spanning the
fibers. We will use the same names also in more general cases even
 if the interpretation of $w$ as fiber coordinate will be lost. Hence
we have:
\begin{equation}\label{ycordi}
   \pmb{y} \equiv y^\alpha \, = \, \left\{u,v,w\right\} \quad;\quad \bar{\pmb{y}} \equiv y^{\bar{\alpha}}\, = \,
   \left\{ \bar{u},\bar{v},\bar{w}\right\}
\end{equation}
An important observation which ought to be done right at the
beginning is that other K\"ahler metrics
$\hat{\mathrm{\mathbf{g}}}_{\alpha\beta^\star}$ do exist on the
three-fold $Y_{[3]}$ that are not Ricci-flat, although the
cohomology class of the associated K\"ahler form $\hat{\mathbf{K}}$
can be the same as the cohomology class of
$\mathbf{K}_{\mathrm{RFK}}$. Within
the framework of the generalized Kronheimer construction, among such K\"ahler (non-Ricci flat) metrics we have the one determined by the K\"ahler quotient according to the formula of Hithchin, Karlhede,
Lindstr\"om and Ro\v{c}ek \cite{HKLR}. Indeed, as
we show later in explicit examples, the K\"ahler metric:
\begin{equation}\label{carigno}
     \text{ds} ^2_{\mathrm{HKLR}}(Y_{[3]})\, = \,\mathrm{\mathbf{g}}^{\mathrm{HKLR}}_{\alpha\beta^\star} \, dy^\alpha \otimes dy^{\beta^\star}
\end{equation}
which emerges from the mathematical K\"ahler quotient construction
and which is naturally associated with $Y_{[3]}$ when this latter is
interpreted as the \textit{space of classical vacua} of the D3-brane
gauge theory (set of extrema of the scalar potential), is
generically non Ricci-flat.
\par
On the other hand on \textit{the supergravity side} of the dual
D3-brane pair we need the Ricci-flat metric in order to construct a
bona-fide D3-brane solution of type IIB supergravity.
In particular, calling $Y^\Gamma_{[3]}$ the
crepant resolution of the $\mathbb{C}^3/\Gamma$ singularity, admitting a Ricci-flat metric, we can construct a
bona-fide D3 brane solution which is solely defined by a single
real function $H$ on $Y^\Gamma_{[3]}$, that should be harmonic with
respect to the Ricci-flat metric, namely:
\begin{equation}\label{merlatino}
    \Box_{\mathbf{g}^{\mathrm{RFK}}} \, H \, = \, 0
\end{equation}
Indeed the function $H(\mathbf{y})$ is necessary and sufficient to
introduce a flux of the Ramond $5$-form so as to produce
the splitting of the $10$-dimensional space into a $4$-dimensional
world volume plus a  transverse $6$-dimensional space that is
identified with the three-fold $Y^\Gamma_{[3]}$. This is the very
essence of the D3-picture.
\par
Yet there is another essential item that was pioneered in
\cite{Grana:2000jj,Gubser:2000vg,Grana:2001xn}
namely the consistent
addition of fluxes for the complex $3$-forms
$\mathcal{H}_\pm$ that appear in the field content of type IIB
supergravity. These provide relevant new charges on both sides of
the gauge/gravity correspondence. In
\cite{Bertolini:2001ma,Bertolini:2002pr} such fluxes were
constructed explicitely relying on a special kind of three-fold:
\begin{equation}\label{specialusto}
    Y_{[3]} \, = \, Y_{[1+2]} \, = \, \mathbb{C}\times \mathrm{ALE}_\Gamma
\end{equation}
where $\mathrm{ALE}_\Gamma$ denotes one of the
$\mathrm{ALE}$-manifolds constructed by Kronheimer \cite{kro1,kro2}
as HyperK\"ahler quotients resolving the singularity
$\mathbb{C}^2/\Gamma$ with $\Gamma \subset \mathrm{SU(2)}$ a finite
Kleinian subgroup.
\par
As we explain in detail below, the essential geometrical feature of
$Y_{[3]}$,  required to construct consistent fluxes of the
complex $3$-forms $\mathcal{H}_\pm$, is that $Y_{[3]}$ should admit
imaginary \textit{(anti)-self-dual, harmonic $3$-forms} $\Omega^{(2,1)}$,
which means:
\begin{equation}\label{alaguerre}
    \star_{\mathbf{g}^{\mathrm{RFK}}} \Omega^{(2,1)} \, = \, \pm \, {\rm i} \,  \Omega^{(2,1)}
\end{equation}
and simultaneously:
\begin{equation}\label{conspicua}
 d\Omega^{(2,1)} \, = \, 0 \quad \Rightarrow \quad
 d\star_{\mathbf{g}^{\mathrm{RFK}}}\Omega^{(2,1)}\, = \, 0
\end{equation}
Since the Hodge-duality operator involves the use of a metric, we
have been careful in specifying that (anti)-self-duality should
occur with respect to the Ricci-flat metric that is the one used 
in the rest of the supergravity solution construction.
\par
The reason why the choice (\ref{specialusto}) of the three-fold
allows the existence of harmonic anti-self dual $3$-forms is easily
understood recalling that the $\mathrm{ALE}_\Gamma$-manifold
obtained from the resolution of $\mathbb{C}^2/\Gamma$ has a compact
support cohomology group of type $(1,1)$ of the following dimension:
\begin{equation}\label{creativina}
    \mathrm{dim}\,\mathrm{H}^{(1,1)}_{comp}\left(\mathrm{ALE}_\Gamma \right) \, = \, r \quad
    \mbox{where}
    \quad r \, = \, \# \,\mbox{ of nontrivial conjugacy classes of
    $\Gamma$}
\end{equation}
Naming $z \in \mathbb{C}$ the coordinate on the factor $\mathbb{C}$
of the product (\ref{specialusto}) and $\omega_I^{(1,1)}$ a basis of
harmonic anti-self dual one-forms on $\mathrm{ALE}_\Gamma$, the
ansatz utilized in \cite{Bertolini:2001ma,Bertolini:2002pr} to
construct the required $\Omega^{(2,1)}$ was the following:
\begin{equation}\label{galeno}
     \Omega^{(2,1)} \, \equiv \, \partial_z
    \, \mathfrak{f}^I (z) \, dz \, \wedge \, \omega^{(1,1)}_I
\end{equation}
where  $\mathfrak{f}^I (z)$ is a set of holomorphic functions of
that variable.  As it is well known $r$ is also the  rank of the
corresponding Lie Algebra in the ADE-classification of the
corresponding Kleinian groups  and the $2$-forms $\omega^{(1,1)}_I$
can be chosen dual to a basis homology cycles $\mathcal{C}_I$
spanning $H_{2}\left(\mathrm{ALE}_\Gamma\right)$, namely we can set:
\begin{equation}\label{samotracia}
    \int_{\mathcal{C}_I} \,\omega^{(1,1)}_J \, = \, \delta_{IJ}
\end{equation}
The form $\Omega^{(2,1)}$ is closed by construction:
\begin{equation}\label{samellus}
   d \Omega^{(2,1)} \, = \,0
\end{equation}
and it is also anti-selfdual with respect to the Ricci-flat metric:
\begin{equation}\label{carondimonio}
     \text{ds}^2_{Y_{[1+2]}} \, = \, dz\otimes d\bar{z}\, + \,  \text{ds}^2_{\mathrm{ALE}_\Gamma}
\end{equation}
Hence the question whether we can construct sufficiently flexible
D3-solutions of supergravity with both $5$-form and $3$-form
fluxes depends on the nontriviality of the relevant cohomology
group:
\begin{equation}\label{ramingo21}
    \mbox{dim} \, \mathrm{H}^{(2,1)}\left(Y_{[3]}\right)\, > \, 0
\end{equation}
and on our ability to find harmonic (anti)-self dual representatives
of its classes (typically not with compact support and hence non
normalizable).
\par
At this level we find a serious difficulty. It seems therefore that
we are not able to find the required $\Omega^{(2,1)}$ forms on
$Y^\Gamma_{[3]}$ and that no D3-brane supergravity solution with
$3$-form fluxes can be constructed dual to the gauge theory obtained
from the Kronheimer construction dictated by $\Gamma \subset
\mathrm{SU(3)}$. Fortunately, the sharp conclusion encoded in
eq.~(\ref{vecchioni}) follows from a hidden mathematical assumption
that, in physical jargon, amounts to a rigid universal choice of the
holomorphic superpotential $\mathcal{W}(\Phi)$. Under appropriate
conditions that we plan to explain and which are detectable at the
level of the McKay quiver diagram,  the superpotential can be
deformed (mass deformation) yielding a family of three-folds
$Y^{\Gamma,\mu}_{[3]}$ which flow, for limiting values of the
parameter ($\mu \to \mu_0$) to a three-fold $Y^{\Gamma,\mu_0}_{[3]}$
admitting imaginary anti self-dual harmonic (2,1)-forms. Since the content and
the interactions of the gauge theory are dictated by the McKay
quiver of $\Gamma$  and by its associated Kronheimer construction,
we are entitled to see its mass deformed version and the exact
D3-brane supergravity solution built on $Y^{\Gamma,\mu_0}_{[3]}$ as
dual to each other.
\par
This will be the object of a future work.  Here we
begin with an accurate mathematical summary of the construction of
D3-brane solutions of type IIB supergravity using the geometric
formulation of the latter within the rheonomy framework
\cite{castdauriafre}.
\subsection{Geometric formulation of Type IIB supergravity }
\label{type2bsum} \setcounter{equation}{0} 
In order to discuss conveniently the D3 brane solutions of type IIB
that have as transverse space the crepant resolution of a
$\mathbb{C}^3/\Gamma$ singularity, we have to recall the geometric
Free Differential Algebra formulation of the chiral ten dimensional
theory fixing with care all our conventions, which is not
only a matter of notations but also of principles and geometrical
insight. Indeed the formulation of type IIB supergravity as it
appears in string theory textbooks
\cite{greenschwarz,Polchinskibook} is tailored for the comparison
with superstring amplitudes and is quite appropriate to this goal.
Yet, from the viewpoint of the general geometrical set up of
supergravity theories this formulation is somewhat unwieldy.
Specifically it neither makes the $\mathrm{SU(1,1)/U(1)}$ coset
structure of the theory manifest, nor does it  relate the
supersymmetry transformation rules  to the underlying algebraic
structure which, as in all other instances of supergravities, is a
simple and well defined {\sl Free Differential algebra}. The Free
Differential Algebra of type IIB supergravity was singled out many
years ago by Castellani in \cite{castella2b} and the geometric,
manifestly $\mathrm{SU(1,1)}$--covariant formulation of the theory
was constructed by Castellani and Pesando in \cite{igorleo}. Their
formulae and  their transcription from a complex $\mathrm{SU(1,1)}$
basis to a real $\mathrm{SL(2,\mathbb{R})}$ basis were summarized
and thoroughly  explained  in a dedicated chapter of a book authored
by one of us \cite{Fre:2013ika} which we refer the reader to.

\subsection{The D3-brane solution with a  $Y_{[3]}$ transverse
manifold} \label{3} 
In this section we discuss a D3-brane solution of type IIB
supergravity in which, transverse to the  brane world-manifold, we
place a smooth non compact three-fold $ Y_{[3]}$ endowed with a
Ricci-flat K\"ahler metric.
\par
The ansatz for the D3-brane solution is characterized by two kinds of
flux; in addition to the usual RR 5-form flux, there is a
non-trivial flux of the supergravity  complex 3-form field strengths
$ \mathcal{H}_{\pm}$.
\par
We separate the ten coordinates of space-time into the following
subsets:
\begin{equation}
 x^M = \left \{ \begin{array}{rcll}
x^\mu &:& \mu =0,1,2,3& \mbox{coordinates of the 3-brane world volume}   \\
y^\tau &:& \tau=4,5,6,7,8,9 & \mbox{real coordinates of the $Y$
variety}   \
\end{array} \right.
\label{coordisplit}
\end{equation}

\subsubsection{The D3 brane ansatz}
\label{sol3flux}  We make the following ansatz for the
metric\footnote{As explained in appendix A, the conventions for the
gamma matrices and the spinors are set with a mostly minus metric
$d\tau^2$. In the discussion of the solution, however, we use $ds^2
= -d\tau^2$ for convenience. We hope this does not cause any
confusion.}:
\begin{eqnarray}
\label{ansazzo}
 \text{ds}^2_{[10]}&=&H(\pmb{y},\bar{\pmb{y}})^{-\frac{1}{2}}\left
(-\eta_{\mu\nu}\,dx^\mu\otimes dx^\nu \right
)+H(\pmb{y},\bar{\pmb{y}})^{\frac{1}{2}} \, \left(
\mathrm{\mathbf{g}}^{\mathrm{RFK}}_{\alpha\beta^\star} \, dy^\alpha \otimes dy^{\beta^\star}\right) \, \nonumber\\
 \text{ds}^2_{Y}&=&\mathrm{\mathbf{g}}_{\alpha\beta^\star}^{\mathrm{RFK}} \, dy^\alpha \otimes dy^{\beta^\star}\nonumber\\
{\rm det}(g_{[10]})&=&H(\pmb{y},\bar{\pmb{y}}){\rm det}(\mathrm{\mathbf{g}^{\mathrm{RFK}}})\nonumber\\
\eta_{\mu\nu}&=&{\rm diag}(+,-,-,-)
\end{eqnarray}
where $\mathrm{\mathbf{g}}^{\mathrm{RFK}}$ is the K\"ahler metric of
the $Y_{[3]}$ manifold
\begin{equation}
  \mathrm{\mathbf{g}}_{\alpha \bar{\beta}}^{\mathrm{RFK}} \, = \, \partial_\alpha \,
  \partial_{\bar{\beta}} \,
  \mathcal{K}^{\mathrm{RFK}}\left(\pmb{y},\bar{\pmb{y}}\right)
\label{m6defi}
\end{equation}
the real function
$\mathcal{K}^{\mathrm{RFK}}\left(\pmb{y},\bar{\pmb{y}}\right)$ being
a suitable K\"ahler potential.
\subsubsection{Elaboration of the ansatz}
In terms of vielbein the ansatz (\ref{ansazzo}) corresponds to
\begin{equation}
V^{A}= \left \{\begin{array}{rcll}
  V^a  & = & H(\pmb{y},\bar{\pmb{y}})^{-1/4} \, dx^a & a=0,1,2,3\\
  V^\ell& = & H(\pmb{y},\bar{\pmb{y}})^{1/4} \,
  \mathbf{e}^\ell & \ell \, = \, 4,5,6,7,8,9
\end{array}\right.
\label{splittoviel}
\end{equation}
where   $\mathbf{e}^\ell$ are the vielbein $1$-forms of the manifold
$Y_{[3]}$. The structure equations of the latter are\footnote{The
hats over the spin connection and the Riemann tensor denote
quantities computed without the warp factor.}:
\begin{eqnarray}
0& = &  d \, \mathbf{e}^i  - \widehat{\omega}^{ij} \, \wedge \, \mathbf{e}^k \, \eta_{jk}\nonumber\\
\widehat{R}^{ij} & = & d \widehat{\omega}^{ij} -
\widehat{\omega}^{ik}\, \wedge \, \widehat{\omega}^{\ell j} \,
\eta_{k\ell} = \widehat{R}^{ij} _{\phantom{ij}\ell m } \,
\mathbf{e}^\ell \,\wedge \, \mathbf{e}^m \label{structeque}
\end{eqnarray}
The relevant property  of the $Y$ metric that we use in solving 
Einstein equations is that it is Ricci-flat:
\begin{equation}
 \widehat{R}^{im}_{\phantom{ij}\ell m } = 0 \label{ricciflatto}
\end{equation}
What we need in order to derive our solution and discuss its
supersymmetry properties is the explicit form of the spin connection
for the full $10$-dimensional metric (\ref{ansazzo}) and  the
corresponding Ricci tensor. From the torsion equation one can
uniquely determine the solution:
\begin{eqnarray}
\omega^{ab} & = & 0 \nonumber\\
\omega^{a\ell} & = &  \ft 1 4 \, H^{-3/2} \, dx^a  \eta^{\ell k} \,
\partial_k \, H
\nonumber\\
\omega^{\ell m} & = & \widehat{\omega}^{\ell m} + \Delta
\omega^{\ell m} \quad ; \quad \Delta \omega^{\ell m} =- \ft 1 2 \,
H^{-1} \, \mathbf{e}^{[\ell}  \, \eta^{m]k} \, \partial_k H
\label{spinconnect}
\end{eqnarray}
Inserting this result into the definition of the curvature $2$-form
we obtain\footnote{The reader should be careful with the indices.
Latin indices are always frame indices referring to the vielbein
formalism. Furthermore we distinguish the 4 directions of the brane
volume by using Latin letters from the beginning of the alphabet
while the 6 transversal directions are denoted by Latin letters from
the middle and the end of the alphabet. For the coordinate indices we
utilize Greek letters and we do exactly the reverse. Early Greek
letters $\alpha,\beta,\gamma,\delta,\dots$ refer to the 6 transverse
directions while Greek letters from the second half of the alphabet
$\mu,\nu,\rho,\sigma,\dots$ refer to the D3 brane world volume
directions as it is customary in $D=4$ field theories. }:
\begin{eqnarray}
R^{a}_{b} & = & -  \frac{1}{8}\, \left [ H^{-3/2} \Box_{\mathbf{g}}
\, H - H^{-5/2} \,
\partial_k H\partial^k H \right] \, \delta^a_b \nonumber\\
R^{a}_{\ell} & = & 0\nonumber\\
R_\ell^m   &=& \frac{1}{8}H^{-3/2} \Box_{\mathbf{g}} H\delta_\ell^m
             - \frac{1}{8} H^{-5/2}\partial_s H\partial^s
             H\delta_\ell^m
             +\frac{1}{4} H^{-5/2} \partial_\ell H\partial^m H
\label{riccius}
\end{eqnarray}
where for any function $f\left(\pmb{y},\bar{\pmb{y}}\right)$ with
support on $Y_{[3]}$:
\begin{equation}
\Box_{\mathbf{g}} \, f\left(\pmb{y},\bar{\pmb{y}}\right) \, =
\,\frac{1}{\sqrt{\mathrm{det}\mathbf{g}}}\, \left( \partial_\alpha
\left(\sqrt{\mathrm{det}\mathbf{g}}\,\,
\mathbf{g}^{\alpha\beta^\star} \,
\partial_{\beta^\star} \,f \right) \right) \label{laplacious}
\end{equation}
denotes the action on it of the Laplace--Beltrami operator  with
respect to the metric (\ref{m6defi}) which is the Ricci-flat one: we
have omitted the superscript $\mathrm{RFK}$ just for simplicity.
Indeed on the supergravity side of the correspondence we use
only the Ricci-flat metric and there is no ambiguity.
\subsubsection{Analysis of the field equations in geometrical terms}
The equations of motion for the scalar fields $\varphi$ and
$C_{[0]}$ and for the 3-form field strength $F^{NS}_{[3]}$ and
$F^{RR}_{[3]}$ can be better analyzed using the complex notation.
Defining, as we did above:
\begin{eqnarray}
{\mathcal{H}}_\pm & = & \pm 2 \,e^{-\varphi/2} F^{NS}_{[3]} + {\rm
i} 2 \,e^{\varphi/2} \,F^{RR}_{[3]}
\label{mcH} \\
P & =& \ft 1 2 \, d\varphi -{\rm i} \ft 12 \, e^\varphi \,
F_{[1]}^{RR} \label{Psc}
\end{eqnarray}
eq.s (\ref{dstarP})-(\ref{hodge2formeq}) can be  respectively
written as:
\begin{eqnarray}
d(\star P)- {\rm i} e^{\varphi} dC_{[0]}\wedge \star P + \ft 1 {16}
   {\mathcal{H}}_+ \, \wedge \, \star
  {\mathcal{H}}_+=0 \label{dstarP} \\
d \star {\mathcal{H}}_+ - \frac{{\rm i}}{2} e^{\varphi} dC_{[0]}
\wedge \, \star {\mathcal{H}}_+=
   {\rm i} \, {F}_{[5]}^{RR}\,
  \wedge \, {\mathcal{H}}_+  - P \wedge \star
  {\mathcal{H}}_- \label{hodge2formeq}
\end{eqnarray}
while  the equation for the 5-form becomes:
\begin{equation}
d\star F^{RR}_{[5]} = {\rm i} \, \ft 1 {8} \, {\mathcal{H}}_+ \wedge
{\mathcal{H}}_- \label{f5equazia}
\end{equation}
\par
Besides assuming the structure (\ref{ansazzo}) we also assume that
the two scalar fields, namely the dilaton $\varphi$ and the
Ramond-Ramond $0$-form $C_{[0]}$ are constant and vanishing:
\begin{equation}
  \varphi=0 \quad ; \quad C_{[0]}=0
\label{zerodilat}
\end{equation}
As we shall see, this assumption simplifies considerably the
equations of motion, although these two scalar fields can be easily restored \cite{Grana:2001xn}.
\subsubsection{The three-forms}
The basic ansatz  characterizing the  solution and  providing  its
interpretation as a D3-brane with three-form fluxes  is described
below.
\par
The ansatz for the complex three-forms of type IIB supergravity is
given below and is inspired by what was done in
\cite{Bertolini:2002pr,Bertolini:2001ma} in the case where $Y_{[3]}
= \mathbb{C}\times \mathrm{ALE}_\Gamma$:
\begin{equation}
{\mathcal{H}}_+ \, = \, \Omega^{(2,1)}  \label{hpmposiz}
\end{equation}
where $\Omega^{(2,1)}$ is localized on $Y_{[3]}$ and satisfies
eq.s~(\ref{alaguerre}-\ref{conspicua})
\par
If we insert the ans\"atze (\ref{zerodilat},\ref{hpmposiz}) into the
scalar field equation (\ref{dstarP}) we obtain:
\begin{eqnarray}
{\mathcal{H}}_+ \, \wedge \, \star_{10} {\mathcal{H}}_+=0
\label{holovincol}
\end{eqnarray}
This equation is automatically satisfied by our ansatz for a very
simple reason that we explain next. The form ${\mathcal{H}}_+$ is by
choice a three-form on $Y_{[3]}$ of type $(2,1)$. Let $\Theta^{[3]}$
be any three-form that is localized on the transverse
six-dimensional \footnote{For the sake of the present calculation
and the following ones where we have to calculate a Hodge dual, it
is more convenient to utilize a set of 6 real coordinates $t^I$
($I=1,\dots,6$) for the manifold $Y_{[3]}$. Let $\partial_I \equiv
\frac{\partial}{\partial t^I}$ denote the standard partial
derivatives with respect to such coordinates.} manifold $Y_{[3]}$:
\begin{equation}\label{terzaformina}
    \Theta^{[3]} \, = \, \Theta_{IJK} \, dt^I\wedge dt^J \wedge dt^K
\end{equation}
When we calculate the Hodge dual of $\Theta^{[3]}$ with respect to
the 10-dimensional metric (\ref{ansazzo}) we obtain a $7$-form with
the following structure:
\begin{equation}\label{salamepuzzoso}
    \star_{10}\,\Theta^{[3]} \, = \, H^{-1} \,
    \mbox{Vol}_{\mathbb{R}^{(1,3)}} \, \wedge \, \widetilde{\Theta}^{[3]}
\end{equation}
where:
\begin{equation}\label{cirimetto}
    \mbox{Vol}_{\mathbb{R}^{(1,3)}} \, = \, \ft {1}{4!} \, dx^{\mu} \wedge dx^{\nu}
    \wedge dx^\rho \wedge dx^\sigma \, \epsilon_{\mu\nu\rho\sigma}
\end{equation}
is the volume-form of the flat D3-brane and
\begin{equation}\label{zante}
    \widetilde{\Theta}^{[3]} \, \equiv \, \star_{\mathbf{g}} \, \Theta^{[3]}
\end{equation}
is the dual of the three-form $\Theta^{[3]}$ with respect to the
metric $\mathbf{g}$ defined on $Y_{[3]}$. Let us now specialize the
three-form $\Theta^{[3]}$ to be of type $(2,1)$:
\begin{equation}\label{specialindo}
    \Theta^{[3]} \, = \, \mathrm{Q}^{(2,1)}
\end{equation}
As shown in \cite{Grana:2000jj,Gubser:2000vg}, preservation of supersymmetry 
requires the complex three-form ${\mathcal{H}}_+$ to obey the condition\footnote{It also requires ${\mathcal{H}}_+$ to be primitive.}
\begin{equation}\label{gerundio}
    \star_{\mathbf{g}} \, \mathrm{Q}^{(2,1)} \, = - {\rm i} \,
   \mathrm{Q}^{(2,1)}
\end{equation}
Hence:
\begin{equation}\label{ciarlatano}
    {\mathcal{H}}_+ \wedge \star_{10}\, {\mathcal{H}}_+ \, = - {\rm i} \,
   \mathrm{Q}^{(2,1)} \, \wedge \,
    \mathrm{Q}^{(2,1)}
     \,\wedge \, H^{-1} \mbox{Vol}_{\mathbb{R}^{(1,3)}}\, = \, 0
\end{equation}
\subsubsection{The self-dual $5$-form} Next we consider the self-dual $5$-form $F_{[5]}^{RR}$
which by definition must satisfy the following
Bianchi identity:
\begin{equation}
d \, F_{[5]}^{RR} = {\rm i} \, \ft 1 8 \,{ \mathcal{H}}_+
\, \wedge \, { \mathcal{H}}_- \label{f5bianchi}
\end{equation}
Our ansatz for $F_{[5]}^{RR}$ is the following:
\begin{eqnarray}
F_{[5]}^{RR} & = & \alpha \left( U + \star_{10}\, U \right)  \label{genova1}\\
U & =  & d \left( H^{-1} \, \mbox{Vol}_{\mathbb{R}^{(1,3)}}
\right) \label{f5ansaz}
\end{eqnarray}
where $\alpha$ is a constant to be determined later. By construction
$F_{[5]}^{RR}$ is self-dual and its equation of motion is trivially
satisfied. What is not guaranteed is that also the Bianchi identity
(\ref{f5bianchi}) is fulfilled. Imposing it, results into a
differential equation for the function
$H\left(\pmb{y},\bar{\pmb{y}}\right)$. Let us see how this works.
\par
Starting from the ansatz (\ref{f5ansaz}) we obtain:
\begin{eqnarray}
  U &=& -\frac{1}{4!} \, \epsilon_{\mu\nu\rho\sigma} \,dx^\mu\wedge dx^\nu\wedge dx^\rho \wedge dx^\sigma \wedge \frac{dH}{H^2} \\
 U_{\mu\nu\rho\sigma  I}  &=& -\frac{1}{4!}
  \epsilon_{\mu\nu\rho\sigma} \, \frac{\partial_I H}{H^2} \quad ;
  \quad
  \mbox{all other components vanish}
\end{eqnarray}
 Calculating the components
of the dual form $\star_{10}\, U$ we find that they are non
vanishing uniquely in the six transverse directions:
\begin{eqnarray}
   \star_{10} U &=& \tilde{U}_{I_1\dots I_5} \,\,dt^{I_1} \wedge \dots \wedge dt^{I_5} \nonumber\\
   \tilde{U}_{I_1\dots I_5} &=& - \, \frac{\sqrt{\mbox{det} \,
   g_{(10)}}}{5!} \epsilon_{I_1\dots I_5 J}\,
   \epsilon_{\mu\nu\rho\sigma} \, g_{(10)}^{JK} \, g_{10}^{\mu\mu^\prime}\, g_{(10)}^{\nu\nu^\prime}\, g_{(10)}^{\rho\rho^\prime}\,
   g_{(10)}^{\sigma\sigma^\prime} \, U_{\mu^\prime\nu^\prime\rho^\prime\sigma^\prime  J}
   \nonumber\\
   &=& \frac{\sqrt{\mbox{det}\mathbf{g}}}{5!} \,\epsilon_{I_1\dots I_5 J} \,
   \mathbf{g}^{JK} \, \partial_K \,H
\end{eqnarray}
The essential point in the above calculation is that all powers of
the function $H$ exactly cancel so that $\star_{10} U$ is linear in
the $H$-derivatives \footnote{Note that we use $\mathbf{g}_{IJ}$ to
denote the components of the K\"ahler metric (\ref{m6defi}) in the
real coordinate basis $t^I$.}. Next using the same coordinate basis
we obtain:
\begin{eqnarray}\label{scolopio}
    d \, F_{[5]}^{RR}& =  & \alpha \, d\star U \, =
    \,\alpha \,
    \underbrace{\frac{1}{\sqrt{\mbox{det} \, \mathbf{g} }} \, \partial_I \,
    \left( \sqrt{\mbox{det} \, \mathbf{g} } \, \mathbf{g}^{IJ} \, \partial_J
    H\right)}_{\Box_{\mathbf{g}} \, H} \, \times \, \mbox{Vol}_{Y_{[3]}} \nonumber\\
    &=& \alpha \,\Box_{\mathbf{g}} \, H(\pmb{y},\bar{\pmb{y}}) \, \times
    \,\mbox{Vol}_{Y_{[3]}}
\end{eqnarray}
where:
\begin{eqnarray}\label{supercaffelatte}
    \mbox{Vol}_{Y_{[3]}} & \equiv &\sqrt{\mbox{det} \, \mathbf{g} } \,
    \frac{1}{6!} \epsilon_{I_1 \dots I_6} dt^{I_1}\wedge \dots
    \wedge dt^{I_6} \nonumber\\
    & = & \sqrt{\mbox{det} \, \mathbf{g}}  \,
    \frac{1}{(3!)^2}\,
    \epsilon_{\alpha\beta\gamma}\,
    dy^\alpha \wedge dy^\beta \wedge dy^\gamma \, \wedge \epsilon_{\bar{\alpha}\bar{\beta}\bar{\gamma}}\,
    d\bar{y}^{\bar{\alpha}} \wedge d\bar{y}^{\bar{\beta}} \wedge
    d\bar{y}^{\bar{\gamma}}
\end{eqnarray}
is the volume form of the transverse six-dimensional space. Once
derived with the use of real coordinates, the relation
(\ref{scolopio}) can be transcribed in terms of complex coordinates
and the Laplace-Beltrami operator $\Box_\mathbf{g}$  can be written
as in eq. (\ref{laplacious}). Let us now analyze the source terms
provided by the three-forms. With our ansatz we obtain:
\begin{eqnarray}
\ft 1 8 \,{ \mathcal{H}}_+  \, \wedge \, {
  \mathcal{H}}_- &= & \mathbb{J}\left(\pmb{y},\bar{\pmb{y}}\right)
  \, \times \,
  \mbox{Vol}_{Y_{[3]}} \nonumber\\
\mathbb{J}\left(\pmb{y},\bar{\pmb{y}}\right) & = & - \, \frac{1}{72
\, \sqrt{\mbox{det}\, \mathbf{g}}} \, \, \,
\Omega_{\alpha\beta\bar{\eta}}\,\, \bar{\Omega}_{\bar{\delta}
\bar{\theta}\gamma} \,\, \epsilon^{\alpha\beta\gamma} \,\,
\epsilon^{\bar{\eta}\bar{\delta}\bar{\theta}}\label{tempra}
\end{eqnarray}
we conclude that the Bianchi identity (\ref{f5bianchi}) is satisfied
by our ansatz if:
\begin{equation}
  \Box_\mathbf{g}\,  H = - \frac
  {1}{\alpha}\, \mathbb{J}\left(\pmb{y},\bar{\pmb{y}}\right)
\label{maindiffe}
\end{equation}
This is the main differential equation to which the entire
construction of the D3-brane solution can be reduced to. We are
going to show that the parameter $\alpha$ is determined by
Einstein's equations and fixed to $\alpha=1$.
\subsubsection{The equations for the three--forms}
Let us consider next the field equation for the complex three-form,
namely eq. (\ref{hodge2formeq}). Since the two scalar fields are
constant the $\mathrm{SU(1,1)/O(2)}$ connection vanishes and  we
have:
\begin{equation}
d \star {\mathcal{H}}_+ = {\rm i} \, F^{RR}_{[5]}\, \wedge
\, {\mathcal{H}}_+ \label{simpeq}
\end{equation}
Using our ansatz we immediately obtain:
\begin{eqnarray}
d \star {\mathcal{H}}_+ = & = &- 2\, {\rm i} H^{-2}
\mathrm{}dH \, \wedge \, \tilde{\Omega}^{(2,1)}, \wedge \,
\Omega_{\mathbb{R}^{1,3}}\, +2 {\rm i} \, H^{-1} \,
d\tilde{\Omega}^{(2,1)} \wedge \, \Omega_{\mathbb{R}^{1,3}}
 \nonumber\\
{\rm i} \, F^{RR}_{[5]}\, \wedge \, {\mathcal{H}}_+ & = & - 2\,
\alpha {\rm i} H^{-2} dH \, \wedge \, \Omega^{(2,1)} \, \wedge \,
\Omega_{\mathbb{R}^{1,3}}\, \label{pagnacco}
\end{eqnarray}
Hence if $\alpha=1$,  the field equations for the three-form reduces
to:
\begin{equation}\label{cannalotto}
    \tilde{\Omega}^{(2,1)} \, \equiv  \, \star_{\mathbf{g}} {\Omega}^{(2,1)} \, = \,
    -\, {\rm i}\,{\Omega}^{(2,1)} \quad ;
    \quad d\star_{\mathbf{g}} {\Omega}^{(2,1)} \, = \,0 \quad ;
    \quad d  {\Omega}^{(2,1)} \, = \,0
\end{equation}
which are nothing else but eq.s~(\ref{alaguerre}-\ref{conspicua}). In
other words the solution of type IIB supergravity with three-form
fluxes exists if and only if the transverse space admits
\textit{closed and imaginary anti-self-dual forms} $\Omega^{(2,1)}$ as we
already stated\footnote{By construction a closed anti-self-dual form
is also coclosed, namely it is harmonic.}.
\par
In order to show that also the Einstein's equation is satisfied by
our ansatz we have to calculate the (trace subtracted) stress energy
tensor of the five and three index field strengths. For this purpose
we need the components of $F_{[5]}^{RR}$. These are easily dealt
with. Relying on the ansatz (\ref{f5ansaz}) and on
eq.~(\ref{splittoviel}) for the vielbein we immediately get:
\begin{equation}
  F_{A_1 \dots A_5} =\left\{ \begin{array}{ccc}
    F_{i abcd} & = & \frac {1} {5!} \, f_i \, \epsilon_{abcd} \\
    F_{ j_1\dots i_5}  & = & \frac{1}{5!} \epsilon_{i j_1\dots j_5} \, f^i \\
    \mbox{otherwise} & = & 0 \
  \end{array} \right.
\label{F5intrinsic}
\end{equation}
where:
\begin{equation}
  f_i = - \alpha\, H^{-5/4} \, \partial_i H
\label{fidefi}
\end{equation}
Then by straightforward algebra we obtain:
\begin{eqnarray}
  T^{a}_{b}\left[ F_{[5]}^{RR}\right]  & \equiv & -75 \, F^{a \,\cdot\, \cdot\, \cdot \,\cdot} \,
  F_{b\,\cdot\, \cdot\, \cdot \,\cdot} = -  \frac 1 8 \, \delta^{a}_{b}
  \, f_p \, f^p \nonumber\\
  & = & - \alpha^2 \, \frac 1 8 \, \delta^{a}_{b}
  \, H^{-5/2} \partial_p H\, \partial^p H \nonumber\\
  T^{i}_{j}\left[  F_{[5]}^{RR}\right]  & \equiv & -75 \, F^{i \,\cdot\, \cdot\, \cdot \,\cdot} \,
  F_{j\,\cdot\, \cdot\, \cdot \,\cdot} = \frac 1 4\,  f^i \, f_j \, - \, \frac 1 8 \,
  \delta^i_j \, f_p \, f^p \nonumber\\
  & = & \alpha^2 \, \left[ \frac 1 4\, H^{-5/4}  \partial^i H \, \partial_j H  \, - \, \frac 1 8
  \, \delta^i_j \,  H^{-5/4}  \partial^p H \, \partial_p H \right]
\label{TofF5}
\end{eqnarray}
\par
Inserting eq.s (\ref{TofF5}) and (\ref{riccius}) into Einstein's
equations:
\begin{eqnarray}
  R^a_b & = &  T^{a}_{b}\left[F_{[5]}^{RR} \right]
 \nonumber\\
  R^i_j & = &T^{i}_{j}\left[F_{[5]}^{RR}\right]
\end{eqnarray}
we see that they are satisfied, provided
\begin{equation}
  \alpha= 1
\label{fixed}
\end{equation}
and the master equation (\ref{maindiffe}) is satisfied. This
concludes our proof that an exact D3-brane solution with a $Y$
transverse space does indeed exist.
\section{An example without mass deformations and no harmonic $\Omega^{(2,1)}$:
 $Y_{[3]}=\mathcal{O}_{\mathbb{P}^2}(-3)$} \label{parteOP23}
\begin{figure}
\begin{center}
\begin{tikzpicture}[scale=0.50]
\draw [thick] [fill=yellow] (-5,-5) circle (1.5cm); \node at (-5,-5)
{$\mathrm{U_2(N)}$}; \draw [thick] [fill=green] (5,-5) circle
(1.5cm); \node at (5,-5) {$\mathrm{U_3(N)}$}; \node at (0,0)
{$\mathrm{U_1(N)}$}; \draw [thick] [fill=red] (0,0) circle (1.5cm);
\node at (0,0) {$\mathrm{U_1(N)}$}; \draw [black, line width=0.07cm]
[->] (-6.5,-5) to [out=90,in=-180] (0,1.5); \draw [blue, line
width=0.07cm] [->] (-3.5,-5) to [out=90,in=-180] (0,-1.5); \draw
[red, line width=0.07cm] [->] (-5,-3.5) to [out=90,in=-180](-1.5,0);
\draw [black, line width=0.07cm] [->] (0,1.5) to [out=0,in=90]
(6.5,-5); \draw [blue, line width=0.07cm] [->] (0,-1.5) to
[out=0,in=90] (3.5,-5); \draw [red, line width=0.07cm] [->] (1.5,0)
to [out=0,in=90](5,-3.5);
\draw [blue, line width=0.07cm] [->] (3.5,-5) to (-3.5,-5); 
\draw [red, line width=0.07cm] [->] (5,-6.5) to (-5,-6.5);
\draw [black, line width=0.07cm] [->] (6.5,-5) to [out=-90,in=-90]
(-6.5,-5);
\end{tikzpicture}
\caption{ \label{quiveruMaster} The quiver diagram of the diagonal
embedding of the group $\mathbb{Z}_3 \to \mathrm{SU(3)}$ }
\end{center}
\end{figure}
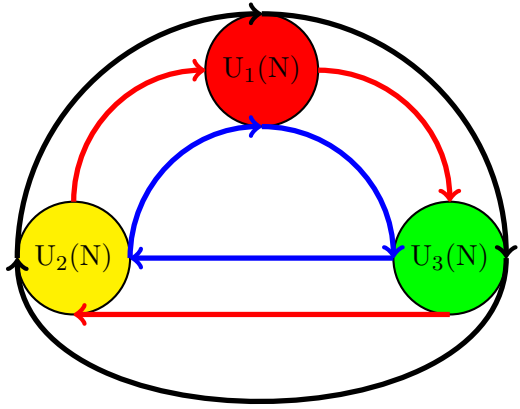

In \cite{Bruzzo:2017fwj} as a master example of the generalized
Kronheimer construction of crepant resolutions the following case
was considered:
\begin{equation}\label{rubelli}
    Y_{[3]} \, = \, \mathcal{O}_{\mathbb{P}^2}(-3)\, \longrightarrow
    \, \frac{\mathbb{C}^3}{\mathbb{Z}_3}
\end{equation}
the action of the group $\mathbb{Z}_3\subset \mathrm{SU(3)}$ on the
three-complex coordinates $\{x,y,z\}$  being generated by the
matrix:
\begin{equation}\label{calienteZ3}
  \mathfrak{g} \, = \, \left(
\begin{array}{ccc}
 e^{\frac{2 {\rm i} \pi }{3}} & 0 & 0 \\
 0 & e^{\frac{2 {\rm  i} \pi }{3}} & 0 \\
 0 & 0 & e^{\frac{2 {\rm i} \pi }{3}} \\
\end{array}
\right)
\end{equation}
Following the steps of the construction one arrives at the following
nine-dimensional flat K\"ahler manifold
\begin{equation}\label{sicuro}
    \mathcal{S}_{\mathbb{Z}_3} \, \equiv \,
    \mbox{Hom}\left(\mathcal{Q}\otimes R,R\right)^{\mathbb{Z}_3}
\end{equation}
where $\mathcal{Q}$ is the three dimensional representation of
$\mathbb{Z}_3$ generated by $\mathfrak{g}$, while $R$ denotes the
regular representation. The points of $\mathcal{S}_{\mathbb{Z}_3}$
are identified with the following triplet of matrices of $3\times 3
$ matrices:
\begin{equation}\label{romil}
 A\, = \,    \left(
\begin{array}{ccc}
 0 & 0 & \Phi^A_{1,3} \\
 \Phi^A_{2,1} & 0 & 0 \\
 0 & \Phi^A_{3,2} & 0 \\
\end{array}
\right) \quad ; \quad B\, = \,    \left(
\begin{array}{ccc}
 0 & 0 & \Phi^B_{1,3} \\
 \Phi^B_{2,1} & 0 & 0 \\
 0 & \Phi^B_{3,2} & 0 \\
\end{array}
\right) \quad ; \quad C\, = \,    \left(
\begin{array}{ccc}
 0 & 0 & \Phi^C_{1,3} \\
 \Phi^C_{2,1} & 0 & 0 \\
 0 & \Phi^C_{3,2} & 0 \\
\end{array}
\right)
\end{equation}
The nine complex coordinates of $\mathcal{S}_{\mathbb{Z}_3}$ are the
matrix entries $\Phi^{A,B,C}_{1,3}$, $\Phi^{A,B,C}_{2,1}$,
$\Phi^{A,B,C}_{3,2}$. With reference to the quiver diagram of
fig.\ref{quiveruMaster}
which is dictated by the McKay matrix $\mathcal{A}_{ij}$ appearing
in the decomposition \footnote{$\mathrm{D}_i$ denote the irreducible
representations of the group $\Gamma = \mathbb{Z}_3$ and each node
of the quiver diagram corresponds to one of them. The number of
lines going from node $i$ to node $j$ is equal to integer value of
$\mathcal{A}_{ij}$. In each node $i$ we have a component
$U_i(n_i\times N)$ of the gauge group $\mathcal{F}_\Gamma$ where
$n_i$ is the dimension of the irrep $D_i$ and $N$ is the number of
D3-branes. }:
\begin{equation}\label{quiverro}
    \mathcal{Q}\otimes \mathrm{D}_i \, = \, \bigoplus_{j=1}^3 \, \mathcal{A}_{ij}\,\mathrm{D}_j
\end{equation}
the entries $\Phi^{A,B,C}_{1,3}, \dots$ are interpreted as the
complex scalar fields of as many Wess-Zumino multiplets in the
bifundamental of the $\mathrm{U_i(N)}$ groups mentioned in the lower
suffix.
\par
In the case of a single brane (N=1) the quiver group
$\mathcal{G}_{\mathbb{Z}_3}$ has the following structure:
\begin{equation}\label{quiverogruppo}
\mathcal{G}_{\mathbb{Z}_3} \, = \, \mathbb{C}^\star \otimes
\mathbb{C}^\star \, \simeq\, \frac{\mathbb{C}^\star \otimes
\mathbb{C}^\star \otimes
\mathbb{C}^\star}{\mathbb{C}^\star_{central}}
\end{equation}
and its maximal compact subgroup $\mathcal{F}_{\mathbb{Z}_3}\subset
\, \mathcal{G}_{\mathbb{Z}_3}$ is the following:
\begin{equation}\label{gaugogruppo}
\mathcal{F}_{\mathbb{Z}_3} \, = \,\mathrm{ U(1)} \otimes
\mathrm{U(1)} \, \simeq \, \frac{\mathrm{ U(1)} \otimes
\mathrm{U(1)}\otimes \mathrm{U(1)}}{\mathrm{U(1)}_{central}}
\end{equation}
The gauge group $\mathcal{F}_{\mathbb{Z}_3}$ and its
complexification $\mathcal{F}_{\mathbb{Z}_3}$ are embedded into
$\mathrm{SL(3,\mathbb{C})}$ by defining the following two
generators:
\begin{equation}\label{carneadillo}
    \mathbf{t}_1 \, = \,\left(
\begin{array}{ccc}
 {\rm i} & 0 & 0 \\
 0 & -{\rm i} & 0 \\
 0 & 0 & 0 \\
\end{array}
\right) \quad ; \quad \mathbf{t}_2 \, = \,\left(
\begin{array}{ccc}
 0 & 0 & 0 \\
 0 & {\rm i} & 0 \\
 0 & 0 & -{\rm i} \\
\end{array}
\right)
\end{equation}
and setting:
\begin{equation}\label{massacarriola}
    \mathcal{F}_{\mathbb{Z}_3} \, = \, \exp\left[\theta_1 \, \mathbf{t}_1
    +\theta_2 \, \mathbf{t}_2\right ] \quad  \theta_{1,2} \in [0,2\pi]\quad ; \quad
    \mathcal{G}_{\mathbb{Z}_3} \, = \, \exp\left[w_1 \, \mathbf{t}_1
    +w_2 \, \mathbf{t}_2\right ] \quad  w_{1,2} \in \mathbb{C}
\end{equation}
\subsection{The HKLR K\"ahler potential}
The K\"ahler potential of the linear space
$\mathcal{S}_{\mathbb{Z}_3}$, which in the D3-brane gauge theory
provides the kinetic terms of the nine scalar fields
$\Phi^{A,B,C}_{1,2,3}$ is given by:
\begin{equation}\label{simelito}
    \mathcal{K}_0\left(\Phi\right)\, = \, \mbox{Tr}\left(A^\dagger\,A + B^\dagger \, B + C^\dagger C \right)
\end{equation}
where the three matrices $A,B,C$ are those of equation
(\ref{romil}). According with the principles of the Kronheimer
construction, the superpotential is given  by
$\mathcal{W}\left(\Phi\right)\, = \, \text{const}\times
\text{Tr}\left( \left[A,B\right] \, C\right)$. The final HKLR
K\"ahler metric, whose determination  requires two steps of physical
significance:
\begin{enumerate}
  \item Reduction to the critical surface of the superpotential
  \textit{i.e.} $\partial_\Phi \,\mathcal{W} \, = \, 0$
  \item Reduction to the level surfaces of the gauge group moment
  maps by solving the algebraic moment map equations,
\end{enumerate}
was calculated in \cite{Bruzzo:2017fwj} according with the general
theory there summarized, which is originally due to the authors of
\cite{HKLR}. The final form of HKLR K\"ahler potential is provided
by:
\begin{eqnarray}\label{HKLRformul}
    \mathcal{K}_{HKLR}(z,\bar{z}, \pmb{\zeta}) & = &\mathcal{K}_0 + \zeta_I \, \mathfrak{C}^{IJ}
    \, \log\left[ \Upsilon_J^{\alpha_{J,\pmb{\zeta}}}  \right]\nonumber\\
    & = & \alpha  \Big\{\left(2 \zeta _1-\zeta _2\right) \log \left[\Upsilon _1\right]-
    \left(\zeta _1-2 \zeta _2\right) \log \left[\Upsilon _2\right]\Big\}+\frac{\Sigma  \left(\Upsilon _1^3+\Upsilon
   _2^3+1\right)}{\Upsilon _1 \Upsilon _2}
\end{eqnarray}
where, as it was extensively discussed in \cite{noietmarcovaldo},
the coefficient $\alpha$ might be adjusted, chamber by chamber, in
chamber space, so as to make the periods of the tautological line
bundles integer on the homology basis.
\par
Setting:
\begin{equation}\label{camillus}
    \Sigma \, \equiv \, |z _1|^2+|z _2|^2+|z _3|^2 \quad ; \quad
    \Upsilon_{1,2} \, = \, \Upsilon _{1,2} \, = \,
    \Lambda_{1,2}\left(\Sigma,\pmb{\zeta}\right)
\end{equation}
where $z_{1,2,3}$ are the three complex coordinates and
$\pmb{\zeta}\, = \, \left\{\zeta_1,\zeta_2\right\}$ the two Fayet-Iliopoulos parameters.  Let us describe the explicit form of these
functions. To this effect let us name $\zeta_1=p$, $\zeta_2 = q$,
and let us introduce the following blocks:
\begin{eqnarray}\label{nastissa}
  \mathfrak{A} & = & \sqrt{p^6 \left(\left(2 p^3 q^3+9 p^2 q \Sigma ^3+9 p q^2
   \Sigma ^3+27 \Sigma ^6\right)^2-4 \left(p^2 q^2+3 p
   \Sigma ^3+3 q \Sigma ^3\right)^3\right)}\nonumber\\
\mathfrak{B} & = &   2 p^6 q^3+9 p^5 q \Sigma ^3+9 p^4 q^2 \Sigma
^3+27 p^3 \Sigma ^6+\mathfrak{A}
\end{eqnarray}
then we have:
\begin{eqnarray}
 \Lambda_{1}\left(\Sigma,p,q\right) &=&  \sqrt[3]{\frac{\sqrt[3]{2} p^4 q^2}{3 \sqrt[3]{\mathfrak{B}}
   \Sigma ^3}+\frac{\sqrt[3]{2}
   p^3}{\sqrt[3]{\mathfrak{B}}}+\frac{\sqrt[3]{2} p^2
   q}{\sqrt[3]{\mathfrak{B}}}+\frac{\sqrt[3]{\mathfrak{B}}}{3
   \sqrt[3]{2} \Sigma ^3}+\frac{p^2 q}{3 \Sigma ^3}+1} \label{Lammo1}\\
 \Lambda_{2}\left(\Sigma,p,q\right) & = & 
{1\over {18\ 6^{2/3} \mathfrak{B}^{2/3} p \Sigma
   ^5}}
\left[\frac{2^{2/3} \mathfrak{B}^{2/3}+2
   \sqrt[3]{\mathfrak{B}} \left(p^2 q+3 \Sigma ^3\right)+2
   \sqrt[3]{2} p^2 \left(p^2 q^2+3 p \Sigma ^3+3 q \Sigma
   ^3\right)}{\sqrt[3]{\mathfrak{B}} \Sigma
   ^3}\right]^{2/3}
   \times \nonumber\\
   &&\left[6 \mathfrak{B}^{2/3} p^2 \Sigma ^3 (p-q)-\sqrt[3]{2}
   \mathfrak{B}^{4/3}+2^{2/3} \mathfrak{B} \left(p^2 q+3 \Sigma
   ^3\right)+ \right. \nonumber \\ && \left.
   2 \sqrt[3]{2 \mathfrak{B}} p^2
   \left(p^4 q^3+3 p^3 q \Sigma ^3+6 p^2 q^2 \Sigma ^3+9
   p \Sigma ^6+9 q \Sigma ^6\right)
 \right.\nonumber\\&&\left.\qquad    -2\ 2^{2/3} p^4
   \left(p^2 q^2+3 p \Sigma ^3+3 q \Sigma ^3\right)^2
   \right]\label{Lammo2}
\end{eqnarray}
\subsection{The issue of the Ricci-flat metric}
One main question is whether the metric arising from the K\"ahler
quotient, which is encoded in eq.~(\ref{HKLRformul}) is Ricci-flat. A
Ricci-flat metric on the crepant resolution of the singularity
$\mathbb{C}^3/\mathbb{Z}_3$, namely on
$\mathcal{O}_{\mathbb{P}^2}(-3)$,  is known in explicit form from the work of Calabi\footnote{Such metrics were also re-discovered in the physics literature in \cite{Page:1985bq}.}
\cite{Calabi-Metriques},
 yet it is not a priori
obvious that the metric defined by the K\"ahler potential
(\ref{HKLRformul}) is that one. The true answer is that it is not,
as we show later on. Indeed we are able to construct directly the 
K\"ahler potential for the resolution of $\mathbb{C}^n/\mathbb{Z}_n$, for any $n\geq 2$, in particular determining the
unique Ricci-flat metric on $\mathcal{O}_{\mathbb{P}^2}(-3)$ with
the same isometries as the metric (\ref{HKLRformul}) and comparing
the two we see that they are different.  Here we stress that the
metric defined by (\ref{HKLRformul}) obviously depends on  the
level parameters $\zeta_1,\zeta_2$ while the Ricci-flat one is
unique up to an overall scale factor. This is an additional reason
to understand a priori that (\ref{HKLRformul}) cannot be the Ricci
flat metric.

Actually Calabi in \cite{Calabi-Metriques} found an easy form of the
K\"ahler potential of a Ricci-flat metric on the canonical bundle of
a K\"ahler-Einstein manifold, and that result applies to the cases
of the canonical bundle of $\mathbb P^2$.
However, in view of applications to cases where we shall
consider the canonical bundles of   manifolds which are not K\"ahler-Einstein,
in the section we stick with our strategy of using the metric coming from
the K\"ahler quotient as a starting point.

\subsubsection{The Ricci-flat metric on $Y_{[3]}
=\mathcal{O}_{\mathbb{P}^2}(-3)$} \label{noriccisiricci} As we have
noticed above the HKLR K\"ahler metric defined by the K\"ahler
potential (\ref{HKLRformul}) depends only on the variable $\Sigma$
defined in eq.~(\ref{camillus}). It follows that the HKLR K\"ahler
metric admits $\mathrm{U(3)}$ as an isometry group, which is the
hidden invariance of $\Sigma$. The already addressed question is
whether the HKLR metric can be Ricci-flat. An almost immediate
result is that a Ricci-flat K\"ahler metric depending only on the
sum of the squared moduli of the complex coordinates is unique (up to a
scale factor) and we can give a general formula for it.
\par
We can present the result in the form of a theorem.
\begin{teorema}
Let $\mathcal{M}_n$ be a non-compact $n$-dimensional K\"ahler
manifold admitting a dense open coordinate patch $z_i$,
$i=1,\dots,n$ which we can identify with the total space of the line
bundle $\mathcal{O}_{\mathbb{P}^{n-1}}(-n)$, the bundle structure
being exposed by the coordinate transformation:
\begin{equation}\label{bundellus}
    z_i\,=\, u_i \, w^{\ft 1 n} \quad , \quad  (i=1,\dots,n-1) \quad \quad ; \quad z_n
    \, = \, w^{\ft 1 n}
\end{equation}
where $u_i$ is a set of inhomogenous coordinates for
$\mathbb{P}^{n-1}$. The K\"ahler potential $\mathcal{K}_n$ of a
$\mathrm{U(n)}$ isometric K\"ahler metric on $\mathcal{M}_n$ must
necessarily be a real function of the unique real variable $\Sigma
\, = \, \sum_{i=1}^n |z_i|^2$. If we require that metric should be
Ricci-flat, the K\"ahler potential is uniquely defined and it is the
following one:
\begin{equation}\label{kaleropolo}
  \mathcal{K}_n(\Sigma)\, = \, k \, +  \frac{\left(\Sigma ^n+\ell ^n\right)^{-\frac{n-1}{n}}
   \left((n-1) \left(\Sigma ^n+\ell ^n\right)-\ell ^n
   \left(\Sigma ^{-n} \ell ^n+1\right)^{\frac{n-1}{n}} \,
   _2F_1\left(\frac{n-1}{n},\frac{n-1}{n};\frac{2
   n-1}{n};-\ell ^n \Sigma ^{-n}\right)\right)}{n-1}
\end{equation}
where $k$ is an irrelevant additive constant and $\ell>0$ is a
constant that can be reabsorbed by rescaling all the complex
coordinates by a factor $\ell$, namely $z_i\to \ell \tilde{z}_i$.
\end{teorema}
\begin{proofteo}
{\rm The proof of the above statement is rather elementary. It
suffices to recall that the Ricci tensor of any  K\"ahler metric
$\mathbf{g}_{ij^\star}\, = \, \partial_i\partial_{j^\star}
\mathcal{K}(z,\bar{z})$ can always be calculated as follows:
\begin{equation}\label{lombillo}
    \mathrm{Ric}_{ij^\star}[\mathbf{g}] \, = \,
    \partial_i\,\partial_{j^\star} \, \log\left[\mathrm{Det}\left[\mathbf{g}\right]\right]
\end{equation}
In order for the Ricci tensor to be zero it is necessary that
$\mathrm{Det}\left[\mathbf{g}\right]$ be the square modulus
of a holomorphic function $|F(z)|^2$, on the other hand under the
hypotheses of the theorem it is a real function of the real variable
$\Sigma$. Hence it must be a constant. It follows that we have to
impose the equation:
\begin{equation}\label{pomigliano}
    \mathrm{Det}\left[\mathbf{g}\right] \, = \, \ell^2 \, = \,
    \mbox{const}
\end{equation}
Let $\mathcal{K}(\Sigma)$ be the sought for K\"ahler potential,
calculating the K\"ahler metric and its determinant we find:
\begin{equation}\label{bagnasco}
  \mathrm{Det}\left[\mathbf{g}\right] \, = \, \Sigma ^{n-1}  \mathcal{K}(\Sigma)' \left(\Sigma ^2\mathcal{K}(\Sigma)''+\Sigma\,   \mathcal{K}(\Sigma)' \right)
\end{equation}
Inserting eq.~(\ref{bagnasco}) into eq.~(\ref{pomigliano}) we obtain a
non linear differential equation for $\mathcal{K}(\Sigma)$ of which
eq.~(\ref{kaleropolo}) is the general integral. This proves the
theorem. $\diamondsuit$}
\end{proofteo}
\subsubsection{Particular cases} It is interesting to analyze
particular cases of the general formula (\ref{kaleropolo}).
\paragraph{\sc The case $n=2$: Eguchi-Hanson.} The case $n=2$ yielding a Ricci
flat metric on $\mathcal{O}_{\mathbb{P}^1}(-2)$ is the Eguchi-Hanson
case namely the crepant resolution of the Kleinian singularity
$\mathbb{C}^2/\mathbb{Z}_2$. This is known to be a HyperK\"ahler
manifold and all HyperK\"ahler metrics are Ricci-flat. Hence also
the HKLR metric must be Ricci-flat and identical with the one
defined by eq.~(\ref{kaleropolo}). Actually we find:
\begin{eqnarray}\label{chicco}
\mathcal{K}_2(\Sigma)&=&\left(\Sigma ^2+\ell
^2\right)^{-\frac{1}{2}}
   \left( \left(\Sigma ^2+\ell ^2\right)-\ell ^2
   \left(\Sigma ^{-2} \ell ^2+1\right)^{\frac{1}{2}} \,
   _2F_1\left(\frac{1}{2},\frac{1}{2};\frac{3}{2};-\ell ^2 \Sigma
   ^{-2}\right)\right)\nonumber\\
   &=&\sqrt{\Sigma ^2+\ell ^2}-\ell  \log
   \left(\sqrt{\Sigma ^2+\ell ^2}+\ell
   \right)+\ell  \log (\Sigma ) \, + \, \mbox{const}
\end{eqnarray}
which follows from the identification of the hypergeometric function
with combinations of elementary transcendental functions occurring
for special values of its indices. The second transcription of the
function is precisely the K\"ahler potential of the Eguchi-Hanson
metric in its HKLR-form as it arises from the Kronheimer
construction (see for instance \cite{Bruzzo:2017fwj}).
\paragraph{\sc The case $n=3$: $\mathcal{O}_{\mathbb{P}^2}(-3)$.} The
next case is that of interest for the D3-brane solution. For $n=3$,
setting $\ell=1$, which we can always do by a rescaling of the
coordinates, we find:
\begin{eqnarray}\label{carriolarotta}
 \mathcal{K}_{Rflat}(\Sigma) &=&   \frac{2 \left(\Sigma ^3+1\right)-
   \left(\frac{1}{\Sigma ^3}+1\right)^{2/3} \,
   _2F_1\left(\frac{2}{3},\frac{2}{3};\frac{5}{3};-\frac{
   1}{\Sigma ^3}\right)}{2 \left(\Sigma ^3+1\right)^{2/3}}\nonumber\\
   &=&\frac{2 (\Sigma ^3+1)-\,
   _2F_1\left(\frac{2}{3},1;\frac{5}{3};\frac{1}{\Sigma
   ^3+1}\right)}{2 \left(\Sigma
   ^3+1\right)^{2/3}}
\end{eqnarray}
The second way of writing the K\"ahler potential follows from one of the standard Kummer relations among hypergeometric
functions. There is a third transcription that also in this case
allows to write it in terms of elementary transcendental functions.
Before considering it we use eq.~(\ref{carriolarotta}) to study the
asymptotic behavior of the K\"ahler potential for large values of
$\Sigma$. We obtain;
\begin{equation}\label{asintInf}
   \mathcal{K}_{Rflat}(\Sigma) \, \stackrel{\Sigma \to \infty}{\approx} \,  \Sigma-\frac{1 }{6
   \Sigma ^2} +\frac{1}{45 \Sigma ^5} +\mathcal{O}\left(\frac{1}{\Sigma^7}\right)
\end{equation}
Eq.(\ref{asintInf}) shows that the Ricci-flat metric is
asymptotically flat since the K\"ahler potential approaches that of
$\mathbb{C}^3$.
\par
\begin{figure}
\centering
\includegraphics[height=8cm]{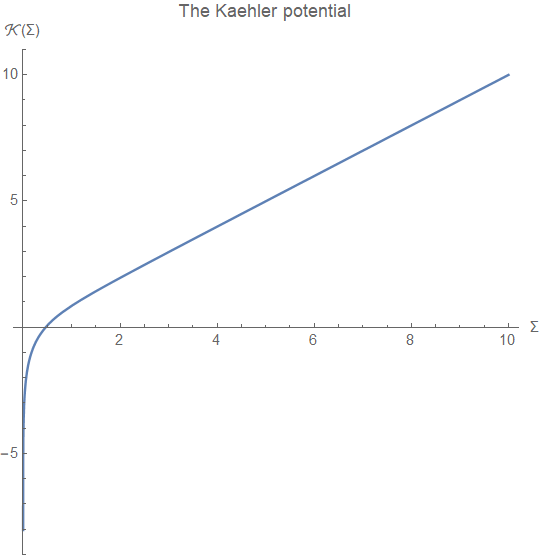}
\caption{ \label{kalleropotto} The plot of the K\"ahler potential
$\mathcal{K}_{Rflat}(\Sigma)$ for the Ricci-flat metric on
$\mathcal{O}_{\mathbb{P}^2}(-3)$. The asymptotic flatness of the
metric is evident from the plot. For large values of $\Sigma$ it
becomes a straight line with angular coefficient $1$. }
\end{figure}

As anticipated, there is an alternative way of writing the K\"ahler
potential (\ref{carriolarotta}) which is the following:
\begin{eqnarray}\label{pagliato}
  \mathcal{K}_{Rflat}(\Sigma)  &=&\frac{\pi }{2
   \sqrt{3}}+ \frac{1}{6} \left(6 \sqrt[3]{\Sigma ^3+1}+2 \log
   \left(\sqrt[3]{\Sigma ^3+1}-1\right)\right.\nonumber\\
   &&\left.-\log
   \left(\left(\Sigma ^3+1\right)^{2/3}+\sqrt[3]{\Sigma
   ^3+1}+1\right)-2 \sqrt{3} \tan ^{-1}\left(\frac{2
   \sqrt[3]{\Sigma
   ^3+1}+1}{\sqrt{3}}\right)\right)
\end{eqnarray}
The identity of eq.~(\ref{pagliato}) with eq.~(\ref{carriolarotta})
can be worked with analytic manipulations that we omit.  The
representation (\ref{pagliato}) is particularly useful to explore
the behavior of the K\"ahler potential at small values of $\Sigma$.
We immediately find that:
\begin{equation}\label{sigmuzero}
    \mathcal{K}_{Rflat}(\Sigma) \, \stackrel{\Sigma \to 0}{\approx} \,
    \log[\Sigma]+\frac{\pi }{2
   \sqrt{3}}+\mathcal{O}\left(\Sigma^6\right)
\end{equation}
The behavior of $\mathcal{K}_{Rflat}(\Sigma)$ is displayed in
fig.\ref{kalleropotto}.
\subsection{The harmonic function in the case $Y_{[3]}
=\mathcal{O}_{\mathbb{P}^2}(-3)$} Let us now consider the equation
for a harmonic function $H(z,\bar{z})$ on the background of the
Ricci-flat metric of $Y_{[3]}$ that we have derived in the previous
sections. Once again we suppose that $H=H(\Sigma)$ is a function
only of the real variable $\Sigma$, {\it viz.} $R=\sqrt{\Sigma}$. For the Ricci-flat metric the
Laplacian equation takes the simplified form: $\partial_i \left(
g^{ij^\star} \, \partial_{j^\star}
    H(\Sigma)\right) \, = \, 0$,
since the determinant of the metric is constant. Using the K\"ahler
metric that follows from the K\"ahler potential $K_{Rflat}(\Sigma)$
defined by eq.s~(\ref{carriolarotta}),(\ref{pagliato}), we obtain a
differential equation that upon the change of variable
$\Sigma=\sqrt[3]{r}$  takes the following form:
\begin{equation}\label{balzoblu}
     3 r (r+1) C''(r)+(5 r+3) C'(r)\, = \, 0
\end{equation}
The general integral  eq.~(\ref{balzoblu}) is displayed below:
\begin{equation}\label{integrai}
    C(r)\, = \, \kappa +\lambda  \left(\log
   \left(1-\sqrt[3]{r+1}\right)-\frac{1}{2} \log
   \left((r+1)^{2/3}+\sqrt[3]{r+1}+1\right)-\sqrt{3} \tan
   ^{-1}\left(\frac{2
   \sqrt[3]{r+1}+1}{\sqrt{3}}\right)\right)
\end{equation}
$\kappa,\lambda$ being the two integration constants. We fix these
latter with boundary conditions. We argue in the following way: if
the transverse space to the brane  were the
original $\mathbb{C}^3/\mathbb{Z}_3$ instead of the resolved variety
$\mathcal{O}_{\mathbb{P}^2}(-3)$, then the harmonic function
describing the D3-brane solution would be the following:
\begin{equation}\label{ciumiglio}
    H_{orb}(R) \, = \, 1 + \frac{1}{R^4} \quad ; \quad R \, \equiv \,
    \sqrt{\Sigma} \, = \, \sqrt{\sum_{i}^2 |z_i|^2} \, = \, \sqrt[6]{r}
\end{equation}
The asymptotic identification for $R\to \infty$ of the Minkowski
metric in ten dimension would be guaranteed, while at small values
of $R$ we would find (via dimensional transmutation) the standard
$\mathrm{AdS_5}$-metric times that of $\mathbb{S}^5$ (see the
following eq.s~(\ref{ciulandario}) and (\ref{adsroba})). In view of
this,  naming $R$ the square root of the variable $\Sigma$, we fix
the coefficients $\kappa,\lambda$ in the harmonic function
$H_{res}(R)$ in such a way that for large values of $R$ it
approaches the harmonic function pertaining to the orbifold case
(\ref{ciumiglio}). The asymptotic expansion of the function: $
H_{res}(R) \, \equiv \, C(r^6) $ is the following one:
\begin{equation}\label{orcoinfi}
   H_{res}(R)\, \stackrel{R\to\infty}{\approx} \,   \left(\lambda -\frac{\pi  \kappa }{2 \sqrt{3}}\right)-\frac{1}{2} \kappa
   \left(\frac{1}{R}\right)^4+O\left(\left(\frac{1}{R}\right)^5\right)
\end{equation}
Hence the function $H_{res}(R)$ approximates the function
$H_{orb}(R)$ if we set  $ \kappa \, = \, 2  \, , \, \lambda =
\frac{\pi}{\sqrt{3}} $. In this way we conclude that:
\begin{eqnarray}\label{resHar}
    H_{res}(R) & = & \frac{1}{3} \left(2 \log \left(\sqrt[3]{R^6+1}-1\right)-\log
   \left(\left(R^6+1\right)^{2/3}+\sqrt[3]{R^6+1}+1\right)\right.\nonumber\\
   &&\left.-2 \sqrt{3} \tan
   ^{-1}\left(\frac{2 \sqrt[3]{R^6+1}+1}{\sqrt{3}}\right)\right)+\frac{\pi}{\sqrt{3}}
\end{eqnarray}
\par
\begin{figure}
\centering
\includegraphics[height=6cm]{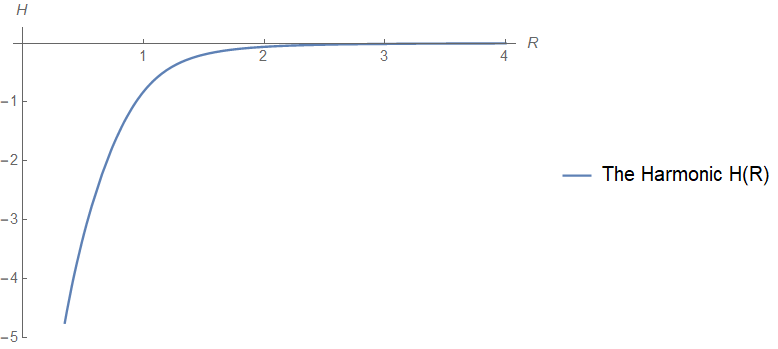}
\caption{ \label{harmoplotto} The plot of the harmonic function
$H_{res}(R)$ for the Ricci-flat metric on
$\mathcal{O}_{\mathbb{P}^2}(-3)$. }
\end{figure}

The overall behavior of the function $H_{res}(R)$ is displayed in
fig.\ref{harmoplotto}.
\subsection{The asymptotic limits of the Ricci-flat metric for the D3-brane solution on $\mathcal{O}_{\mathbb{P}^2}(-3)$
} In the case of a standard D3-brane on $Y_{[3]}\, = \,
\mathbb{C}^3\simeq \mathbb{R}^6$ one writes the same ansatz as in
eq.~(\ref{ansazzo}) and (\ref{genova1}-\ref{f5ansaz}) where now the
K\"ahler metric is $ \mathbf{g}_{\alpha\beta^\star} \, = \,
\delta_{\alpha\beta^\star} $ Rewriting the complex coordinates in
terms of polar coordinates $ z_1= e^{{\rm i} \varphi _1} R \cos \phi ,\,
z_2= e^{{\rm i} \varphi_2} R \cos \chi  \sin \phi ,\, z_3= e^{{\rm i}
\varphi_3} R \sin \chi    \sin \phi$ we obtain that:
\begin{eqnarray}\label{fiasaporita}
   \text{ds}^2_{\mathbb{C}^3} &\equiv& \sum_{i=1}^3 |dz_i|^2\,= \, dR^2+R^2 \, \text{ds}^2_{\mathbb{S}^5}
\end{eqnarray}
where:
\begin{equation}\label{cinquesfera}
  \text{ds}^2_{\mathbb{S}^5} \, = \,    d{\varphi_1}^2 \cos ^2 \phi +\sin ^2\phi 
   \left(d{\varphi_2}^2 \cos ^2\chi +d{\varphi_3}^2 \sin ^2\chi +d\chi^2\right)+d\phi^2
\end{equation}
is the $\mathrm{SO(6)}$-invariant metric of a $5$-sphere in polar
coordinates. In other words the Ricci-flat  K\"ahler metric
$ \text{ds}^2_{\mathbb{C}^3}$ (which is also Riemann-flat) is that of the
\textit{metric cone} on the Sasaki-Einstein metric of
$\mathbb{S}^5$.  At the same time the $\mathrm{SO(6)}$-invariant
harmonic function on $\mathbb{C}^3$ is given by the already quoted
$H_{orb}(R)$ in (\ref{ciumiglio}), and the complete $10$-dimensional
metric of the D3-brane solution takes the form:
\begin{equation}\label{ciulandario}
    \text{ds}^2_{10|orb}\,
    =\,\frac{1}{\sqrt{1+\frac{1}{R^4}}}\,  \text{ds}^2_{\mathrm{Mink}_{1,3}}
    \, + \,\sqrt{1+\frac{1}{R^4}}\, \left(dR^2+R^2 \,
      \text{ds}^2_{\mathbb{S}^5}\right)
\end{equation}
For $R\to\infty$ the metric (\ref{ciulandario}) approaches the flat
Minkowski metric in $d=10$, while for $R\to 0$ it approaches the
following metric:
\begin{equation}\label{adsroba}
      \text{ds}^2_{10|orb}\, \stackrel{R\to 0}{\approx} \, \underbrace{R^2
    \,\,  \text{ds}^2_{\mathrm{Mink}_{1,3}} \, + \, \frac{dR^2}{R^2}}_{\mathrm{AdS}_5} \, +
    \,\underbrace{  \text{ds}^2_{\mathbb{S}^5}}_{\mathbb{S}^5}
\end{equation}
Let us now cosider the asymptotic behavior of the Ricci-flat metric
on $\mathcal{O}_{\mathbb{P}^2}(-3)$. In order to obtain a precise
comparison with the flat orbifold case the main technical point is
provided by the transcription of the $\mathbb{S}^5$-metric in terms
of coordinates well adapted to the Hopf fibration:
\begin{equation}\label{opfo}
    \mathbb{S}^5 \, \stackrel{\pi}{\longrightarrow} \, \mathbb{P}^2
    \quad ; \quad \forall\;  p \in \mathbb{P}^2 \quad \pi^{-1}(p) \sim
    \mathbb{S}^1
\end{equation}
To this effect let $Y=\left\{u,v\right\}$ be a pair of complex
coordinates for $\mathbb{P}^2$ such that the standard Fubini-Study
metric on this compact 2-fold is given by:
\begin{equation}\label{cinesino}
     \text{ds}^2_{\mathbb{P}^2}\, = \, g^{\mathbb{P}^2}_{ij^\star}\, dY^i \,  d\bar{Y}^{j^\star}\, \equiv  \,
   \frac{ {dY}\cdot d\bar{Y}}{1+Y\cdot \bar{Y}}-\frac{\left(\bar{Y}\cdot d {Y}\right)\left(Y\cdot d\bar{Y}\right)}
    {\left(1+Y\cdot \bar{Y}\right)^2}
\end{equation}
the corresponding K\"ahler 2-form being $\mathbb{K}_{\mathbb{P}^2}\,
= \, \frac{{\rm i}}{2\pi} \,g^{\mathbb{P}^2}_{ij^\star}\, dY^i \,\wedge \,
d\bar{Y}^{j^\star}$. Introducing the one form: $ \Omega \, =
\,\frac{{\rm i} \left({Y}\cdot d\bar{Y}-\bar{Y} \cdot dY\right)}{2
\left(1+Y\cdot \bar{Y}\right)}$ whose exterior derivative is the
K\"ahler 2-form,  $d\Omega \, = \, 2\pi \,\mathbb{K}_{\mathbb{P}^2}
$, the metric of the five-sphere in terms of these variables is the
following one: 
\begin{equation}\label{ds5S}
     \text{ds}^2_{\mathbb{S}^5} \, = \,  \text{ds}^2_{\mathbb{P}^2} \, + \,
    \left(d\varphi+ \Omega\right)^2
\end{equation}
where the range of the coordinate $\varphi$ spanning the
$\mathbb{S}^1$ fiber is $\varphi \in \left[0,2\pi\right]$. In this
way the flat metric on the metric cone on $\mathbb{S}^5$,
namely (\ref{fiasaporita}) can be rewritten as follows:
\begin{equation}\label{camillino}
     \text{ds}^2_{\mathbb{C}^3} \, = \, dR^2 + R^2 \,  \text{ds}^2_{\mathbb{P}^2} \,
    + R^2  \left(d\varphi + \Omega\right)^2
\end{equation}
\subsubsection{Comparison of the Ricci-flat metric with the orbifold
metric}\label{Rgrande} In order to compare the exact Ricci-flat
metric streaming from the K\"ahler potential (\ref{carriolarotta})
with the metric (\ref{camillus}) it suffices to turn to toric coordinates
\begin{equation}\label{toricume}
z_1 = u \sqrt[3]{w}\, , \quad  z_2\, = \,  v \sqrt[3]{w}\, , \quad
z_3 \, = \, \sqrt[3]{w} \quad ; \quad \Sigma \, = \,
\left(1+\varpi\right) \, {\mathfrak{f}}^{1/3} \quad ; \quad
\varpi=|u|^2+|v|^2 \quad ; \quad {\mathfrak{f}} \, = \, |w|^2
\end{equation}
The toric coordinates $\{u,v\}\equiv Y$ span the exceptional divisor
$\mathbb{P}^2$ while $w$ is the fiber coordinate in the bundle.
Setting:
\begin{equation}\label{salmoiraghi}
    w \, = \, e^{{\rm i}\psi}  |w| \, = \, e^{{\rm i}\psi}  \, \left(\frac{R^2}{1 + |u|^2 + |v|^2}\right)^{\ft
    32}
\end{equation}
we obtain:
\begin{eqnarray}\label{caniggia}
     \text{ds}^2_{Rflat} & = & h(R) dR^2 \, + \, f(R) \,  \text{ds}^2_{\mathbb{P}^2}
    \, + \, g(R) \, \left(d\psi + 3 \, \Omega \right)^2 \nonumber\\
      f(R) & = &\sqrt[3]{R^6+1} \quad ; \quad  h(R) \, = \, \frac{R^4}{\left(R^6+1\right)^{2/3}} \quad ;
    \quad g(R) \, = \, \frac{R^6}{9 \left(R^6+1\right)^{2/3}}
    \quad
\end{eqnarray}
From eq.~(\ref{caniggia}) we derive the asymptotic form of the metric
for large values of $R$, namely:
\begin{equation}\label{sibillata}
   \text{ds}^2_{Rflat}\, \stackrel{R \to \infty}{\approx}\,
   {dR}^2+R^2 \text{ds}_{\mathbb{P}_2}^2+R^2 \left(\frac{\text{d$\psi $}}{3}+\Omega \right)^2
\end{equation}
The only difference between eq.~(\ref{camillino}) and
eq.~(\ref{sibillata}) is the range of the angular value $\varphi\, =
\,\frac{\psi}{3}$. Because of the original definition of the angle
$\psi$, the new angle $\varphi\in\left[0,\frac{2\pi}{3}\right]$ takes
one third of the values. This means that the asymptotic metric cone
is quotiened by $\mathbb{Z}_3$ as it is natural since we resolved
the singularity $\mathbb{C}^3/\mathbb{Z}_3$.
\subsubsection{Reduction to the exceptional
divisor}\label{ecceziunale} The other important limit of the
Ricci-flat metric is its reduction to the exceptional divisor
$\mathcal{ED}$. In the present case the only fixed point for the
action of $\Gamma=\mathbb{Z}_3$ on $\mathbb{C}^3$ is provided by the
origin $z_{1,2,3}=0$ which, comparing with eq.~(\ref{toricume}),
means $w=0 \Rightarrow \mathfrak{f}=0$. This is the equation of the
exceptional divisor which is created by the blowup of the unique
singular point.  In the basis of the complex toric coordinates
$Y^i\, \equiv \{u,v,w\}$, the K\"ahler metric derived from the
K\"ahler potential (\ref{carriolarotta}) has the following
appearance:
\begin{equation}\label{conigliofritto}
 g_{ij^\star}^{Rflat} \, = \,    \left(
\begin{array}{ccc}
 \frac{v \bar{v}+\mathfrak{f} (\varpi +1)^4+1}{(\varpi +1)^2 \left(\mathfrak{f} (\varpi
   +1)^3+1\right)^{2/3}} & -\frac{v \bar{u}}{(\varpi +1)^2 \left(\mathfrak{f} (\varpi
   +1)^3+1\right)^{2/3}} & \frac{w (\varpi +1)^2 \bar{u}}{3 \left(\mathfrak{f} (\varpi
   +1)^3+1\right)^{2/3}} \\
 -\frac{u \bar{v}}{(\varpi +1)^2 \left(\mathfrak{f} (\varpi +1)^3+1\right)^{2/3}} &
   \frac{u \bar{u}+\mathfrak{f} (\varpi +1)^4+1}{(\varpi +1)^2 \left(\mathfrak{f}
   (\varpi +1)^3+1\right)^{2/3}} & \frac{w (\varpi +1)^2 \bar{v}}{3 \left(\mathfrak{f}
   (\varpi +1)^3+1\right)^{2/3}} \\
 \frac{u (\varpi +1)^2 \bar{w}}{3 \left(\mathfrak{f} (\varpi +1)^3+1\right)^{2/3}} &
   \frac{v (\varpi +1)^2 \bar{w}}{3 \left(\mathfrak{f} (\varpi +1)^3+1\right)^{2/3}} &
   \frac{(\varpi +1)^3}{9 \left(\mathfrak{f} (\varpi +1)^3+1\right)^{2/3}} \\
\end{array}
\right)
\end{equation}
where the invariants $\mathfrak{f},\varpi$ are defined in equation
(\ref{toricume}). Hence the reduction of the metric to the exceptional
divisor is obtained by setting $\mathrm{d}w=\mathrm{d}{\bar w}=0$ in
the line element $\text{ds}^2_{Rflat}\, \equiv \, g_{ij^\star}^{Rflat}  dY^i \,
\,d\bar{Y}^{j^\star} $ and performing the limit
$\mathfrak{f}\to 0$ on the result. We obtain:
\begin{equation}\label{tersilla}
     \text{ds}^2_{Rflat}\, \stackrel{\mathcal{ED}}{\longrightarrow} \,
    \text{ds}_{\mathbb{P}^2}^2 \, \equiv \,
    \frac{\mathrm{d}v\left( \mathrm{d}{\bar v}+u\,{\bar u} \mathrm{d}{\bar v}-u{\bar v}\,\mathrm{d}{\bar u}\right)+\mathrm{d}u
    \left( \mathrm{d}{\bar u}+ v\,\bar{v}\mathrm{d}{\bar u}-{\bar u}\,v\,\mathrm{d}{\bar v}\right)}{\left(1+u \bar{u}+v \bar{v}\right)^2}
\end{equation}
which is the standard Fubini-Study metric on $\mathbb{P}^2$ obtained
from the K\"ahler potential:
\begin{equation}\label{fubinistudy}
    \mathcal{K}_{\mathbb{P}^2}^{FS}(\varpi)\, = \,  \log\left(1+\varpi\right)
\end{equation}
As we see, the metric on the exceptional divisor obtained from the
Ricci-flat metric has no memory of the Fayet Iliopoulos (or
stability parameters) $p,q$ which characterize instead the HKLR
metric obtained from the Kronheimer construction. This is obvious
since the Ricci-flat metric does not depend on $p,q$. On the other
hand the HKLR metric, that follows from the K\"ahler potential
(\ref{HKLRformul}), strongly depends on the Fayet Iliopoulos
parameters
 \(\zeta_1=p\, , \, \zeta_2=q\) and one naturally expects that the reduction of $\text{ds}^2_{HKLR}$
 to the exceptional divisor will inherit such a dependence. Actually
 this is not the case since the entire dependence from $p,q$ of the HKLR K\"ahler
 potential, once reduced to $\mathcal{ED}$, is localized in an overall
 multiplicative constant and in an irrelevant additive constant. This
 matter of fact is conceptually very important in view of our conjecture that the Ricci-flat metric is completely
 determined, by means of the Monge-Amp\`{e}re equation, from the
 K\"ahler metric on the exceptional divisor, as it is determined by the
 Kronheimer construction. In the present case where, up to a
 multiplicative constant, \textit{i.e. a homothety} there is only one Ricci-flat metric on
 $\mathcal{O}_{\mathbb{P}^2}(-3)$ with the prescribed isometries, our conjecture
 might be true only if the reduction of the HKLR metric to the
 exceptional divisor is unique and $p,q$-independent, apart from
 overall rescalings. It is very much reassuring that this is
 precisely what actually happens.
\section{The case $Y\to \mathbb{C}^3/\mathbb{Z}_4$ and the general problem of determining a Ricci-flat metric}

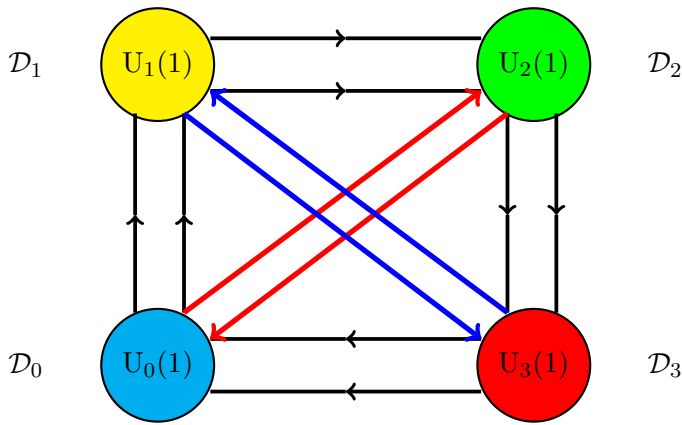
\begin{figure}
\begin{center}
\begin{tikzpicture}[scale=0.50]
\draw [thick] [fill=yellow] (-5,0) circle (1.5 cm) ; \node at (-5,0)
{$\mathrm{U_1(1)}$}; \draw [thick] [fill=green] (5,0) circle (1.5
cm); \node at (5,0) {$\mathrm{U_2(1)}$};  \draw [black,line
width=0.05cm] [->] (-3.6,0.7)-- (0,0.7); \draw [black,line
width=0.05cm] (0,0.7)-- (3.6,0.7); \draw [black,line width=0.05cm]
(3.6,-0.7)-- (0,-0.7); \draw [black,line width=0.05cm][<-]
(0,-0.7)-- (-3.6,-0.7);
 \node at
(8.5,0) {$\mathcal{D}_2$}; \node at (-8.5,0) {$ \mathcal{D}_1$};
\draw [thick] [fill=cyan] (-5,-8) circle (1.5 cm) ; \node at (-5,-8)
{$\mathrm{U_0(1)}$}; \draw [thick] [fill=red] (5,-8) circle (1.5
cm); \node at (5,-8) {$\mathrm{U_3(1)}$};  \draw [black,line
width=0.05cm] (-3.6,-8.7)-- (0,-8.7); \draw [black,line
width=0.05cm][->] (3.6,-7.3)--(0,-7.3); \draw [black,line
width=0.05cm][->] (3.6,-8.7)--(0,-8.7) ; \draw [black,line
width=0.05cm] (-3.6,-7.3)--(0,-7.3) ;
 \node at
(8.5,-8) {$\mathcal{D}_3$}; \node at (-8.5,-8) {$ \mathcal{D}_0$};
\draw [black,line width=0.05cm] [->](-4.3,-6.6)-- (-4.3,-4.0); \draw
[black,line width=0.05cm][->] (-5.6,-6.6)--(-5.6,-4.0); \draw
[black,line width=0.05cm] (-4.3,-4)-- (-4.3,-1.3); \draw [black,line
width=0.05cm](-5.6,-4)--(-5.6,-1.3);
\draw [black,line width=0.05cm] (4.3,-6.6)-- (4.3,-4.0); \draw
[black,line width=0.05cm](5.6,-6.6)--(5.6,-4.0); \draw [black,line
width=0.05cm][->](4.3,-1.3)-- (4.3,-4) ; \draw [black,line
width=0.05cm][->](5.6,-1.3)--(5.6,-4);
\draw [red, line width=0.07cm] [->] (-4.3,-6.6) --(3.6,-0.7);
\draw [red, line width=0.07cm] [->] (4.3,-1.3) --
(-3.6,-7.3); 
\draw [blue, line width=0.07cm] [->] (4.3,-6.6) --(-3.6,-0.7);
\draw [blue, line width=0.07cm] [->] (-4.3,-1.3) --(3.6,-7.3) ;
\end{tikzpicture}
\caption{\label{c3z4quiverotto} The quiver diagram describing the
$\mathbb{C}^3/\mathbb{Z}_4$ singular quotient and codifying its
resolution via K\"ahler quotient \`{a} la Kronheimer. The same
quiver diagram codifies the construction of the corresponding gauge
theory for a stack of D3-branes. Each node is associated with one of
the 4 irreducible representations of $\mathbb{Z}_4$ and in each node
we located one of the $\mathrm{U_i(1)}$ groups with respect to which
we perform the K\"ahler quotient. This is the case of one
D3-brane. For $N$ D3--branes, all gauge groups $\mathrm{U_i(1)}$ are
promoted to $\mathrm{U_i(N)}$.}
\end{center}
\end{figure}

The next case of interest to us at present is the
resolution $Y\to \mathbb{C}^3/\mathbb{Z}_4$ whose associated
Kronheimer construction was studied in detail in
\cite{noietmarcovaldo}. (A study of  $ \mathbb{C}^3/\mathbb{Z}_4$ as a non-complete intersection affine variety
in $\mathbb C^9$ is presented in the Appendix.) The corresponding MacKay quiver is displayed
in fig.\ref{c3z4quiverotto}. Differently from the case of the
resolution $Y\to \mathbb{C}^3/\mathbb{Z}_3$ studied in section
\ref{parteOP23}, here the HKLR K\"ahler metric cannot be derived
explicitly since the moment map equations form a system of
algebraic equations of higher degree. Yet as it was explained in
\cite{noietmarcovaldo} one can work out the restriction of such metric
to the  compact component of the exceptional divisor which is the
second Hirzebruch surface $\mathbb{F}_2$. Indeed it was shown that
the quotient singularity $\mathbb{C}^3/\mathbb{Z}_4$ can be
completely resolved by ${\rm tot} K_{\mathbb{F}_{2}}$ 
\cite{noietmarcovaldo}, that denotes the total space of the
canonical bundle over the second Hirzebruch surface.
\par
Hence the main goal we would like to achieve is the construction of
a Ricci-flat K\"ahler metric on ${\rm tot} K_{\mathbb{F}_{2}}$  which
restricted to the base $\mathbb{F}_2$ of the bundle hopefully
coincides with K\"ahler metric on the same surface provided by the
Kronheimer construction.
\par
 Being a non-compact Calabi-Yau variety the
existence of a Ricci-flat K\"ahler metric on 
${\rm tot} K_{\mathbb{F}_{2}}$  is not implied by the classic Yau theorem, valid for
smooth compact manifolds. To ask whether Ricci-flat metrics do
exist, one has to specify boundary conditions.
 We will be interested in metrics that, just as in the
previous example, are asymptotically conical, namely of the
form\footnote{Note that without specifying the boundary conditions
there can exist more than one Ricci-flat metric. Explicit examples
of non-asymptotically conical Ricci-flat metrics in six real
dimensions can be found in \cite{Lu:2004ya}.}
\begin{equation}\label{bequeta}
     \text{ds}^2 (Y) \approx dR^2 + R^2  \text{ds}^2 (X_5)
\end{equation}
for a suitable radial coordinate approaching $R\to \infty$.
Essentially by definition, $ \text{ds}^2 (X_5)$ is a Sasaki-Einstein metric
on a compact manifold (or orbifold) $X_5$. Then we fix the boundary
conditions for our metric by requiring that asymptotically it
approaches the cone over $\mathbb{S}^5/\mathbb{Z}_4$. With this
boundary condition\footnote{The results in \cite{Joyce-QALE} require some more
precise estimate on the fall-off of the metric at infinity.} the
 theorems in \cite{Joyce-QALE} imply the existence of a unique Ricci-flat
K\"ahler metric in every K\"ahler class of the resolved variety $Y$.
Analogous existence results for isolated quotient singularities $\mathbb{C}^m/\Gamma$ were given in  \cite{Joyceisolated}
and later extended in \cite{goto} and \cite{vanCoevering:2008fn} for crepant resolutions of general isolated conical singularities.
See also \cite{vancoe2} for of applications of the general existence results in the toric context, including the resolution of the conical singularities on
the $Y^{p,q}$ Sasaki-Einstein five-manifolds \cite{Gauntlett:2004yd}.
\par
The existence results are analogous to Yau's theorem in the compact
case. In fact, recently there has been some renewed interest and
activity in this area, with some new results concerning for example
the existence of Sasaki-Einstein manifolds, outside the toric realm.
These results are related to the idea of ``stability''. For
reference, recent work on this subject include
\cite{Xie:2019qmw,Fazzi:2019gvt}.
\par
For many purposes, knowledge of the existence of a metric, together
with some of its key properties, can be sufficient for extracting
interesting physical information. This is true also in the case of
the AdS/CFT correspondence. However, if one is interested in
constructing the metrics explicitly, namely write them down in some
coordinate systems, then the existence theorems are not helpful,
because they are not constructive (as far as we know).
\par The
classic examples of explicit Ricci-flat K\"ahler metrics in real
dimension four include Eguchi-Hanson, Gibbons-Hawking, Taub-NUT,
Atiyah-Hitchin. In real dimension six, for a long time the resolved
and deformed metrics on the conifold singularity constructed by
Candelas and de la Ossa \cite{Candelas:1989js}
were the only (non trivial) known examples of explicit
Ricci-flat K\"ahler metrics. The so-called ``resolved conifold'' metric is a metric
on the total space of the vector bundle ${\cal O}(-1)\oplus {\cal
O}(-1)\to \mathbb{P}^1$, the isometry group is $\mathrm{SU(2)\times
SU(2)\times U(1)}$ and  asymptotically it approaches the
cone over the Sasaki-Einstein manifold $\mathrm{T^{1,1}}$ (with the
same isometry).  In other cases,  different kind of resolutions
exist, where instead of a $\mathbb{P}^1$ one replaces the
singularity with a compact four-dimensional manifold (or orbifold)
${\cal M}_4$. A general ansatz that yields explicit Ricci-flat
K\"ahler metrics was constructed by Page and Pope (in any dimension)
\cite{Page:1985bq}, but  this is somewhat
limited as it \emph{assumes} that the metric induced on ${\cal M}_4$
is K\"ahler-Einstein. Explicit K\"ahler-Einstein metrics on smooth
four-dimensional manifolds are known only for ${\cal M}_4 =
\mathbb{P}^2$ and   ${\cal M}_4 = \mathbb{P}^1\times \mathbb{P}^1$.
The former leads to the construction of an explicit Ricci-flat
K\"ahler metric on the total space of ${\cal O}_{\mathbb{P}^2}(-3)\simeq \mathrm{tot}K_{\mathbb{P}^2}$, 
which is the resolution of the quotient singularity
$\mathbb{C}^3/\mathbb{Z}_3$ and was fully described in section
\ref{parteOP23} (see also\cite{Calabi-Metriques}). The latter leads
to the construction of an explicit Ricci-flat K\"ahler metric on the
total space of $\mathrm{tot}K_{\mathbb{P}^1\times \mathbb{P}^1}$,
which is the resolution of the conical singularity
(conifold)$/\mathbb{Z}_2$. The corresponding Sasaki-Einstein
manifolds at infinity are, respectively, $\mathbb{S}^5/\mathbb{Z}_3$
and $\mathrm{T^{1,1}}/\mathbb{Z}_2$.
  For the case of $\mathrm{tot}K_{\mathbb{P}^1\times \mathbb{P}^1}$, a generalisation was constructed \cite{dueparpandoro}, 
  namely  an explicit Ricci-flat K\"ahler metric that depends on the two independent
K\"ahler classes parameters: this construction however uses the
$\mathrm{SU(2)\times SU(2)\times U(1)}$ symmetry and as a result the
metric is co-homogeneity one, although it does not fit in the ansatz
of \cite{Calabi-Metriques} and \cite{Page:1985bq}. Recently, the ansatz of \cite{Calabi-Metriques,Page:1985bq} was used to produce explicit Ricci-flat K\"ahler metrics on the canonical bundle of generalised flag manifolds \cite{CorreaGrama}. Extensions that include the dependence on several   K\"ahler class parameters have appeared in  \cite{vancoeflag,bykovflag}.

\subsection{The Ricci-flat K\"ahler metric on  ${\rm tot} K_{\mathbb{F}_{1}}$}
The metric that we shall present in the sequel has some distinctive features that are shared with an
explicit Ricci-flat K\"ahler metric  on   ${\rm tot} K_{\mathbb{F}_1}$,
where $ \mathbb{F}_1$ is the first Hirzebruch
surface, \emph{i.e.}, the first del Pezzo surface
$\mathrm{dP}_1$, constructed in  \cite{Oota:2006pm}. This metric is many ways ``more complicated'' than
all the other metrics mentioned above. Let us summarise some of its
salient properties:
\begin{enumerate}
\item Asymptotically it approaches the cone over the
Sasaki-Einstein manifold\footnote{Incidentally, $Y^{2,1}$ can also
be viewed as circle bundle over $\mathbb{F}_1$. See section 5 of
\cite{Martelli:2004wu}.} $Y^{2,1}$.
\item The isometry group is $\mathrm{SU(2)\times U(1) \times U(1)}$.
\item It is cohomogeneity two. In particular, there is a homogeneous
base, given by a round $\mathbb{P}^1$, and then the metric
depends non-trivially on two coordinates.
\item It is \emph{toric}, in that there is a $\mathrm{U(1)^3 \in  SU(2)\times U(1) \times U(1)}$
subgroup of isometries that leaves invariant the K\"ahler form, and contains the torus
of the toric three-fold ${\rm tot} K_{\mathbb{F}_1}$. This group allows one
 to introduce three moment map coordinates and
three angular coordinates
(``action-angle'' coordinates system).
\item It also possesses an additional ``hidden symmetry'' corresponding to the existence of a so-called Hamiltonian two-form
\cite{Martelli:2007pv}, that implies the existence of a coordinate
system (called ``orthotoric'') in which the metric components are
all given in terms of functions of one variable.
\item Imposing this extra symmetry however, comes at the price of loosing one of the two K\"ahler class parameters.
Indeed it was later demonstrated in  \cite{Bykov:2017mgc} that the
two-parameter metric (that is known to exist thanks to the general
theorems of \cite{goto,vanCoevering:2008fn}) does not posses such Hamiltonian
two-form.
\item The metric induced on exceptional divisor ${\cal M}_4=  \mathbb{F}_1$ is obviously K\"ahler,
but it is not Einstein. Indeed,   a K\"ahler-Einstein
metric on  $\mathbb{F}_1$ does not exist.
\item In  \cite{Martelli:2007pv} (further explored in detail in  \cite{Martelli:2007mk})
it was shown that this metric is part of a family of (in general
only \emph{partial}\footnote{This means that for general $p$ and $q$
the compact divisor ${\cal M}_4$ has orbifold singularities
\cite{Martelli:2007mk,Martelli:2007pv}. This is because the metric
ansatz is ``too simple'' to account for all the necessary K\"ahler
class parameters; but completely resolved metrics are known to
exist \cite{vancoe2}. For the special case $p=2$ and $q=1$ the metric is
completely smooth. We also note that in \cite{Martelli:2007pv} were
constructed different types of partial resolutions, corresponding to
various ``chambers''. Moreover, the paper discusses the case of general
dimension, while for our purposes we shall focus on  the case of real dimension
$d=6$.})
resolutions of the conical Ricci-flat metrics on the whole family of
$Y^{p,q}$ Sasaki-Einstein manifolds.
\item  In \cite{Bykov:2017mgc} it is given a relation between the orthotoric coordinates and
a set of complex coordinates that is well adapted to the complex
structure of  ${\rm tot} K_{\mathbb{F}_{1}}$, with one complex
coordinate on the non-compact fiber $\mathbb{C}$, one coordinate on
the fiber $\mathbb{P}^1$ in $\mathbb{F}_1$ and one coordinate on the
base $\mathbb{P}^1$ in $\mathbb{F}_1$.
\item A set of local complex coordinates explicitly related to the
orthotoric coordinates was given in section  2.2 of \cite{Martelli:2007pv}.
It would be interesting to work out the relation between these and the complex coordinates defined in \cite{Bykov:2017mgc}.
\end{enumerate}
Since, similarly to $ \mathbb{F}_1$, also  $ \mathbb{F}_2$
does not admit a K\"ahler-Einstein metric, the Ricci-flat metric on
${\rm tot} K_{\mathbb{F}_{2}}$  cannot be found through the Calabi ansatz 
\cite{Calabi-Metriques,Page:1985bq}. We expect the Ricci-flat metric on
${\rm tot} K_{\mathbb{F}_{2}}$ to share many features with that on
${\rm tot} K_{\mathbb{F}_{1}}$ , summarised above.  One difference is
that at infinity it must approach the cone over the Sasaki-Einstein
orbifold $\mathbb{S}^5/\mathbb{Z}_4$, as opposed to the cone over
the Sasaki-Einstein manifold $Y^{2,1}$. The Ricci-flat metric on
${\rm tot} K_{\mathbb{F}_{2}}$ will also be toric and moreover it
should have again isometry group $\mathrm{SU(2)\times U(1)\times
U(1)}$. This immediately implies that the metric should be
co-homogeneity two and in practice it leads to PDE's in two
variables. For example, one can write the Monge-Ampere equation for
the K\"ahler potential as a PDE in two variables, or similarly the
corresponding equation for the symplectic potential. Without further
assumptions, these equations are unlikely to be solvable in closed
form.
\par
A natural assumption to make is that the metric admits a Hamiltonian
two-form, namely that it can be put in the orthotoric form. This is
natural because the partial resolution of all the $Y^{p,q}$
singularities arise in this ansatz, with $p=2$, $q=1$ giving the
complete resolution above. Strictly speaking the
$p>q>0$ should hold, however, it is known that by performing a scaling
limit of the $Y^{p,q}$ Sasaki-Einstein metrics, one can recover the
limiting cases $Y^{p,p}=\mathbb{S}^5/\mathbb{Z}_{2p}$ and
$Y^{p,0}=\mathrm{T^{1,1}}/\mathbb{Z}_{p}$, suggesting that the
partial resolution metrics may also be extended to these regimes of
parameters\footnote{In fact, in Appendix A of \cite{Martelli:2007pv}
the metric ansatz of \cite{Page:1985bq} is recovered in a limit.}.

\section{A general set up for a metric ansatz with  separation of variables}

\label{duefunziesecche} In the sequel we begin by considering a metric on a
6-dimensional manifold $\mathcal{M}_6$ which is K\"ahler and by
construction admits SU(2)$\times $U(1)$\times $U(1) as an isometry
group. This metric depends on two functions $\Upsilon (s)$ and
$P(t)$ of two real coordinates $s,t$ invariant with respect to the
isometry group. The other coordinates are four angles, with ranges and periodicities specified according with
the following summary table:
\begin{equation}\label{iniz1}
\begin{array}{cccccc}
 s\leq -3\, ,  &  -\tfrac{3}{2} \leq t\leq 0,&  0 \leq  \theta \leq \pi  
, & \phi \in [0,2\pi]\, ,  & \tau \in [0,2\pi]\, ,  & \chi \in [0,\frac{3}{2}\pi]\\
\end{array}
\end{equation}
The metric, which is defined by means of the following vielbein
\begin{eqnarray}
\label{iniz2}
E^1&=&\frac{1}{2} \sqrt{s t} {d\theta}\nonumber\\
E^2&=&\frac{1}{2} \sqrt{s t} \sin\theta {d\phi} \nonumber\\
E^3&=&\frac{1}{2} \,\sqrt{\frac{s-t}{3+s}}\,\Upsilon (s) d s\nonumber\\
E^4&=&\frac{1}{2} \,\sqrt{t-s} P(t)\, dt\nonumber\\
E^5&=&-\frac{1}{\sqrt{\frac{s-t}{3+s}}\,\Upsilon
(s)}\,\left[-\frac{1}{2} t \left({d\tau}+(1-\cos{\theta})
 d\phi -\frac{2 {d\chi}}{3}\right)+{d\chi}\right]\nonumber\\
E^6&=&-\frac{1}{\sqrt{t-s} P(t)} \left[-\frac{1}{2} s
\left({d\tau}+(1-\cos{\theta}) d{\phi}-\frac{2
{d\chi}}{3}\right)+{d\chi}\right]
\end{eqnarray}
is derived, by generalization, from the \textit{orthotoric metrics}
discussed in\footnote{In particular, see the line element (4.1) in \cite{Martelli:2007mk},  after correcting some typos in that expression. 
The relation to our coordinates is given by $t=y-1$, $s=x-1$. Moreover, we have $\theta_\mathrm{here}= \theta_\mathrm{there}$,  
$\phi_\mathrm{here}= \phi_\mathrm{there}$, as well as $\chi_\mathrm{here}= \tau_\mathrm{there}$, $\tau_\mathrm{here}=2\psi_\mathrm{there}+\tfrac{2}{3}\tau_\mathrm{there}$.}
\cite{Martelli:2007mk,Martelli:2007pv}
where the relation of latter with the metrics on Sasakian 5-manifolds
$Y^{p,q}$ is also  presented. Although in those references it was assumed that $p>q$, presently we will consider setting $p=q=2$ and show that this yields an
orthotoric metric that we shall identify as a Ricci-flat K\"ahler
metric on ${\rm tot} K_{{\mathbb WP}[112]}$. The asymptotic metric corresponds to a cone over the
limiting case $Y^{2,2}=\mathbb{S}^5/\mathbb{Z}_4$ of the Sasaki-Einstein manifolds $Y^{p,q}$ \cite{Gauntlett:2004yd}.
\par
The line-element:
\begin{eqnarray}
    \text{ds}_{\text{ort}}^2&=&\sum _{i=1}^6
    \left(E^i\right){}^2\label{iniz3}\\
 &=&\frac{1}{4} s \, t \,
\left({d\theta}^2+\sin^2\theta{d\phi}^2 \right)+\frac{(s-t)
 \Upsilon (s)^2}{4 (3+s)} \, d s^2
+\frac{1}{4} (t-s) P(t)^2 \, \text{dt}^2\nonumber\\
&&+\frac{(3+s)}{(s-t) \Upsilon
(s)^2} 
\left[-\frac{1}{2} t \left(\text{d$\tau
$}+(1-\cos{\theta}) {d\phi}-\frac{2 \text{d$\chi
$}}{3}\right)+{d\chi}\right]^2\nonumber\\
&&+\frac{1}{(t-s)
P(t)^2}
\left[-\frac{1}{2} s \left(\text{d$\tau $}+(1-\cos{\theta})
{d\phi}-\frac{2 \text{d$\chi $}}{3}\right)+{d\chi}\right]^2\label{iniz9}
\end{eqnarray}
is K\"ahlerian by construction since it admits the following closed
K\"ahler 2-form:
\begin{eqnarray}\label{iniz4}
\mathbb{K}_{\text{ort}}&=&E_1\wedge E_2+E_3\wedge E_5+E_4\wedge
E_6\nonumber\\
&=&\frac{1}{2} \left\{t \left[-\frac{1}{2} \cos\theta
\text{ds}\wedge {d\phi}+\frac{1}{2} \left(ds\wedge
{d\tau}-\frac{2}{3} ds\wedge {d\chi}\right)+\frac{1}{2} ds\wedge {d\phi}\right] -ds\wedge {d\chi} \right\}\nonumber\\
&&+\frac{1}{2} \left\{ s \left[-\frac{1}{2} \cos\theta
\text{dt}\wedge {d\phi}+\frac{1}{2} \left(\text{dt}\wedge
{d\tau}-\frac{2}{3} \text{dt}\wedge {d\chi}\right)+\frac{1}{2} \text{dt}\wedge {d\phi}\right]-\text{dt}\wedge {d\chi}\right\}\nonumber\\
&&+\frac{1}{4} s t \sin{\theta} {d\theta}\wedge
{d\phi}
\end{eqnarray}
Indeed $\mathbb{K}_{\text{ort}}$ is closed by construction and it is
a K\"ahler 2-form since we have:
\begin{equation}\label{iniz5}
    \mathbb{K}_{\text{ort}}=\frac{1}{2}\sum _{i=1}^6 \sum _{j=1}^6
\mathbb{J}_{\text{ij}}E^i\land E^j
\end{equation}
where:
\begin{equation}\label{iniz6}
    \mathbb{J}\,=\,\left(
\begin{array}{cccccc}
 0 & 1 & 0 & 0 & 0 & 0 \\
 -1 & 0 & 0 & 0 & 0 & 0 \\
 0 & 0 & 0 & 0 & 1 & 0 \\
 0 & 0 & 0 & 0 & 0 & 1 \\
 0 & 0 & -1 & 0 & 0 & 0 \\
 0 & 0 & 0 & -1 & 0 & 0 \\
\end{array}
\right)\quad;\quad \mathbb{J}^2= - \mathbf{1}
\end{equation}
is an antisymmetric tensor which squares to minus the identity,
namely it is a frame-index complex structure tensor.
\par
It should be noted that the K\"ahler form in eq.~(\ref{iniz4}) is
independent from the two functions $\Upsilon(s)$ and $P(t)$, namely
it is universal for an entire class of metrics.

\subsection{The orthotoric metric on $K_{{\mathbb WP}[112]}$}
 Within the general scope of the above described setup we have that
the metric (\ref{iniz9}) is Ricci-flat for the following choice of
the two functions parameterizing the line-element:
\begin{equation}\label{iniz10}
  \Upsilon (s)\,=\,
  \sqrt{\frac{-s}{\frac{2}{3}\,s^2-s+3}}\quad ;\quad
  P(t)\,=\, \frac{1}{\sqrt{-\frac{2}{3}t^2-t}}
\end{equation}
With the choice (\ref{iniz10}), from eq.~(\ref{iniz9}) we obtain:
\begin{eqnarray}\label{iniz11}
\text{ds}_{{\rm tot}  K_{{\mathbb WP}[112]}}^2&=&\frac{1}{4} \left\{\frac{4
\left(\frac{2 s^2}{3}+s+\frac{9}{s}\right)
\left[\left(\frac{t}{3}+1\right) {d\chi}-\frac{1}{2}
t [(1-\cos{\theta}){d\phi} +{d\tau}]\right]^2}{t-s}+\right.\nonumber\\
&&\left.+\frac{4 t (2 t+3) \left(\left[\frac{s}{3}+1\right)
{d\chi}-\frac{1}{2} s [(1-\cos\theta){d\phi} +{d\tau}]\right]^2}{3 (s-t)}\right.\nonumber\\&& \left.+s
t \left(\sin^2\theta{d\phi}^2+d\theta^2\right)+\frac{ds^2 (t-s)}{\frac{2 s^2}{3}+s+\frac{9}{s}}+\frac{3
dt^2 (s-t)}{t (2 t+3)}\right\}
\end{eqnarray}
The reason for the subscript ${\rm tot}  K_{{\mathbb WP}[112]}$  is that the Ricci
flat metric (\ref{iniz11}) turns out to be defined over the total
space of the canonical bundle of the (singular) projective space
${{\mathbb WP}[112]}$ namely on
${\rm tot}  K_{{\mathbb WP}[112]}$.
\par
It is a simple matter to verify that asymptotically, for $s\to -\infty$, the metric (\ref{iniz11}) is indeed approximatively conical, and therefore   Quasi-ALE \cite{Joyce-QALE}. To see this, one can set
$s=-\tfrac{2}{3}R^2$, so that  
 \begin{equation}
 \text{ds}_{{\rm tot}  K_{{\mathbb WP}[112]}}^2  \, \stackrel{R \to \infty}{\approx}\, dR^2 + R^2 \text{ds}^2_{X_5} 
 \end{equation}
at leading order in $R$.  Since the metric is  Ricci-flat K\"ahler, and it takes the form of a cone over a five-dimensional space, it follows that locally  the five-dimensional metric $\text{ds}^2_{X_5}$
is a Sasaki-Einstein metric. In appendix  \ref{Dariosappendix} we discuss the metric $\text{ds}^2_{X_5}$ in more detail, showing that $X_5=\mathbb{S}^5/\mathbb{Z}_4$, with a specific $\mathbb{Z}_4$ action.

As we will show in sect. \ref{genreductia},  the metric induced by (\ref{iniz11}) on the
exceptional divisor ${{\mathbb WP}[112]}$ is the same as
the one obtained on that space while resolving a
$\mathbb{C}^3/\mathbb{Z}_4$ orbifold singularity by means of the
Kronheimer construction localized on the unique type III  wall
$\mathcal{W}_2$ displayed by its chamber structure (see sect. 6.4 of
\cite{noietmarcovaldo}).

\subsection{Integration of the complex structure and the complex
coordinates}\label{leccalecca} In their algebraic geometry
description, the varieties of the type here considered are complex
threefolds $\mathcal{K}_3$ that are canonical bundles of some
compact K\"ahler two-fold $\mathcal{D}_2$ which, on its turn, is
the total space of a line-bundle over $\mathbb{P}_1$:
\begin{equation}\label{iniz12}
   \mathcal{M}_6=\mathcal{K}_3\overset{\pi }{\longrightarrow
}\mathcal{D}_2 \overset{\tilde{\pi }}{\longrightarrow } \mathbb{P}_1
\end{equation}
This hierarchical structure implies a hierarchy in the complex
coordinates that can be organized and named in the following way
according with the nomenclature of \cite{noietmarcovaldo}:
\begin{equation}\label{iniz13}
\begin{array}{rclcrcl}
  u &=& \text{coordinate on the }\mathbb{P}_1\text{ base of }\mathcal{D}_2
  &;&
  v &=& \text{coordinate on the fibers of }\mathcal{D}_2\\
  w &=& \text{coordinate on the fibers of }\mathcal{K}_3&\null&\null\\
  \end{array}
\end{equation}
This structure is reflected in the integration of the complex
structure that can be deduced from the combination of the K\"ahler
2-form with the metric.
\subsubsection{The path to the integration}
Indeed, having the metric and the K\"ahler form we can construct the
complex structure tensor. Then we try to integrate the complex
structure we have found. This is very important in order to organize
the fibred structure of the manifold. First from eq.~(\ref{iniz2})
one reads off the vielbein $E_{\mu }^i$ defined as:
\begin{eqnarray}\label{iniz14}
  \pmb{E}^i &=& E_{\mu }^i d x^{\mu }\quad ; \quad x^{\mu } \,=\,
   \{s,t,\theta ,\phi ,\tau ,\chi \}
\end{eqnarray}
The $6\times 6$ matrix $E_{\mu }^i$ depends only on the $s, t$
variables and on the angle $\theta $ (as we will see $\theta $ can
be traded for the coordinate $\rho =\tan\frac{\theta}{2}$ and in the
symplectic formalism it is a moment variable). The true angular
variables are the phases of the three complex coordinates namely
$\phi $, $\tau $, $\chi $. As a next step one introduces the inverse
vielbein which is just the matrix inverse of $E_{\mu }^i$ according
with the definition
\begin{equation}\label{iniz15}
    E_{\mu }^i E_j^{\mu } =\delta _j^i
\end{equation}
This enables us to write the differentials of the coordinates as
linear combinations of the vielbein ${dx}^{\mu }= E_j^{\mu }
\pmb{E}^i$.
\paragraph{\sc The complex structure tensor in
coordinate indices $\mathbb{JW}$.} Using the vielbein matrix and its
inverse we can convert the frame indices of the complex structure
tensor to coordinate ones and we get:
\begin{equation}\label{iniz17}
\mathbb{JW}_\mu{}^\nu=\left(
\begin{array}{cccccc}
 0 & 0 & 0 & 0 & -\frac{1}{3} \Upsilon (s)^2 & -\frac{s \Upsilon (s)^2}{6+2 s} \\
 0 & 0 & 0 & 0 & -\frac{1}{3} (3+t) P(t)^2 & -\frac{1}{2} t P(t)^2 \\
 0 & 0 & 0 & \csc{\theta} & -\tan\frac{\theta }{2} & 0 \\
 -\frac{2 (3+s) t \sin^2\frac{\theta }{2}}{(s-t) \Upsilon (s)^2} & -\frac{2 s \sin^2\frac{\theta }{2}}{(t-s)
P(t)^2} & -\sin{\theta} & 0 & 0 & 0 \\
 -\frac{(3+s) t}{(s-t) \Upsilon (s)^2} & -\frac{s}{(t-s) P(t)^2} & 0 & 0 & 0 & 0 \\
 \frac{2 (3+s) (3+t)}{3 (s-t) \Upsilon (s)^2} & \frac{2 (3+s)}{3 (t-s) P(t)^2} & 0 & 0 & 0 & 0 \\
\end{array}
\right)
\end{equation}
\subsubsection{Integration of the autodifferentials}
The matrix $\mathbb{JW}$ has three eigenvectors corresponding to the
eigenvalue ${\rm i}$ and three corresponding to the eigenvalue $-{\rm i}$ (their
complex conjugates). The three eigenvectors corresponding to ${\rm i}$ are
the rows of the following matrix
\begin{equation}\label{iniz18}
    Y_{\mu }^i=\left(
\begin{array}{cccccc}
 \frac{{\rm i} s \Upsilon (s)^2}{6+2 s} & \frac{1}{2} {\rm i} t P(t)^2 & 0 & 0 & 0 & 1 \\
 \frac{1}{3} {\rm i} \Upsilon (s)^2 & \frac{1}{3} {\rm i} (3+t) P(t)^2 & {\rm i} \tan\frac{\theta }{2} & 0 & 1 & 0 \\
 0 & 0 & - {\rm i} \csc{\theta} & 1 & 0 & 0 \\
\end{array}
\right)
\end{equation}
combining with the differentials $dY^i \, =
\,Y_{\mu}^i dx^{\mu }$ we obtain three closed one-forms $
   d\left(dY^i\right)\,=\, 0$
that can be integrated to yield the three-complex variables $u,v$
and $w$.
\paragraph{\sc The coordinate $u$} is obtained from the integration of
$dY^3$:
\begin{equation}\label{iniz21}
   u =\exp{{\rm i}\int (-{\rm i} \csc{\theta} d\theta +d\phi)}
= \tan\frac{\theta }{2} e^{{\rm i}\phi}
\end{equation}
\paragraph{\sc The coordinate $v$} is obtained from the integration of
$dY^2$:
\begin{equation}\label{iniz22}
    v =\exp {\rm i}\int \left(d\tau+\frac{1}{3} {\rm i} (3+t)
P(t)^2 \,dt + {\rm i}  \tan\frac{\theta}{2}\,
d\theta+\frac{1}{3} {\rm i} \Upsilon (s)^2 \, ds \right)= 
\cos^2\frac{\theta }{2} H(t) \, \Psi (s)e^{{\rm i} \tau }
\end{equation}
where we have introduced the following new functions of $t$ and $s$:
\begin{equation}\label{iniz23}
    H(t)=\exp\left(-\frac{1}{3}\int_{\text{const}}^t(3+x)P(x)^2dx\right)
\quad;\quad \Psi(s)=\exp\left(-\frac{1}{3}\int_{-\infty }^s \Upsilon
(x)^2 \, dx \right)
\end{equation}
\paragraph{\sc The coordinate $w$} is obtained from the integration of
$dY^1$. Here we are not assisted by $\mathrm{SU(2)}$
invariance to define the exact coefficient in front of the
differential. We choose a coefficient that appears reasonable from
the result and what we obtain is either the coordinate $w$ of other
approaches or a power $w^a$. In the sequel comparing with the
construction from the iterative procedure we will see what is the
correct identification of the power $a$. At the beginning our
educated guess suggests the use of a coefficient $4/3$. So we set
\begin{equation}\label{iniz24}
w=\exp {\rm i}\frac{4}{3} \,\int \left( d\chi +\frac{1}{2} {\rm i} \, t \,
P(t)^2 dt\, +\, \frac{{\rm i} \, s \,  \Upsilon (s)^2}{6+2 s}\, ds\right) =
\Phi(s) K(t) e^{{\rm i}\frac{4}{3}\chi}
\end{equation}
where we have introduced the new functions:
\begin{equation}\label{iniz25}
    \Phi(s)=\exp\left(-\frac{2}{3}\int _{-\infty }^s\frac{x
}{x+3}\Upsilon (x)^2 dx \right)\quad;\quad
K(t)=\exp\left(-\frac{2}{3}\int_{\text{const}}^t x\,
P(x)^2\,dx\right)
\end{equation}
One necessary property that must be possessed by the function
$\Phi(s)$ is:
\begin{equation}\label{iniz26}
    \Phi(-3) = 0
\end{equation}
which defines the exceptional divisor at $w = 0$

Notice that with the ranges of the coordinates that we specified  in (\ref{iniz1}), we see that $u$ is a complex coordinate on a $\mathbb{P}^1$, while $\nu$ and $w$ are complex coordinates on two copies of $\mathbb{C}$.

\section{AMSY symplectic formalism and transcription
of the metric in this formalism} \label{sympatia} According to the formalism introduced by Abreu 
\cite{abreu} and developed by Martelli, Sparks and Yau 
\cite{Martelli:2005tp}, in the case of toric K\"ahler varieties
of complex dimension $n$, one can find moment maps  $\mu^i$ and
angular variables $\Theta_i$ such that the K\"ahler 2-form takes the
universal form:
\begin{equation}\label{buko1}
    \mathbb{K}=\sum _{i=1}^n d\mu ^i\wedge d\Theta_i
\end{equation}
At the same time there exist a function $G(\mu ^i)$ of the $n$ real
moment variables, named the \textit{symplectic potential}, such that
the metric takes the following universal form:
\begin{equation}\label{buko2}
    \text{ds}_{\text{symp}}^2= G_{ij} d\mu^i d\mu^j
    +\left(G^{-1}\right){}^{ij }d\Theta_i
d\Theta_j
\end{equation}
where by definition: $G_{ij}\equiv \partial _{i,j}G\text{
}$  is the Hessian of the symplectic potential and  $
\left(G^{-1}\right){}^{ij}$ is the inverse of the Hessian
matrix.
\par
In our case the three angular variables are $ \pmb{\Theta } =
\{\phi ,\tau ,\chi \}$ and the K\"ahler form is given by
$\mathbb{K}$ as defined in eq. (\ref{iniz4}). Transforming the
pseudo angle $\theta $ to the variable $\rho $ by setting $ \theta
\,=\, 2 \arctan\rho $ and implementing such change of variables in
the K\"ahler form we obtain:
\begin{eqnarray}\label{buko5}
  \mathbb{K} &=& \frac{1}{12} \left[3 t ds\wedge
  d\tau +\frac{6 t \rho ^2 ds\wedge d\phi}{1+\rho ^2}
  -2 (3+t) ds\wedge d\chi+3 s \left(dt\wedge
d\tau +\frac{2 \rho ^2 dt\wedge d\phi }{1+\rho
^2}\right)\right.\nonumber\\&&\left. -2 (3+s) dt\wedge d\chi
+\frac{12 s t \rho d\rho \wedge d\phi }{\left(1+\rho
^2\right)^2}\right]
\end{eqnarray}
which is compatible with eq.~(\ref{buko1}) if the coefficient of each
of the three angular variables $\tau$, $\chi$ , $\phi$  is a closed differential
that can be integrated to a single new moment coordinate function of
the real coordinates $\rho,s,t$. Hence we introduce the vector of
moments:
\begin{equation}\label{buko6}
    \mu =\{\mathfrak{u},\mathfrak{v},\mathfrak{w}\}
\end{equation}
and the K\"ahler 2-form (\ref{buko5}) can be rewritten as:
\begin{equation}\label{buko7}
    \mathbb{K}=d\mathfrak{u}\wedge d\phi +d\mathfrak{v}\wedge
d\tau +d\mathfrak{w}\wedge d\chi
\end{equation}
provided we have defined the coordinate transformation:
\begin{equation}\label{buko8}
   \mathfrak{u}\,=\,\frac{s \, t \,\rho ^2}{2+2 \rho ^2}\quad ; \quad
   \mathfrak{v}\,=\, \frac{s\, t}{4}\quad ; \quad \mathfrak{w}=\frac{1}{6} \left(-3 t-s (3+t)\right)
\end{equation}
The unique inverse transformation of the above coordinate change is
the following one:
\begin{equation}\label{buko9}
\begin{array}{rclcrcl}
  \rho &=& \frac{\sqrt{\mathfrak{u}}}{\sqrt{-\mathfrak{u}+2 \mathfrak{v}}}
  &;&
  t &=& \frac{1}{6} \left(-4 \mathfrak{v}-6
\mathfrak{w}+\sqrt{-144 \mathfrak{v}+(4 \mathfrak{v}+6
\mathfrak{w})^2}\right)\\
\null &\null&\null&\null&  s &=& \frac{1}{3} \left(-2 \mathfrak{v}-3
\mathfrak{w}-\sqrt{4 (-9+\mathfrak{v}) \mathfrak{v}+12 \mathfrak{v}
\mathfrak{w}+9 \mathfrak{w}^2}\right)\\
\end{array}
\end{equation}
The new real coordinates are named
$\mathfrak{u}$,$\mathfrak{v}$,$\mathfrak{w}$ with gothic letters
since they are the symplectic counterparts of the complex
coordinates $u,v,w$ yet, differently from the latter, we do not need
the complex structure to find them and hence they are independent
from the metric.
\subsection{Transcription of the metric in the toric symplectic form}
At this point we try to rewrite the metric depending on the two
functions:
\begin{equation}\label{duefunzie}
    M(s,t) \equiv \sqrt{\frac{s-t}{3+s}}\, \Upsilon(s) \quad ; \quad
    \Phi(s,t) \equiv \sqrt{t-s} \, P(t)
\end{equation}
in the symplectic form (\ref{buko2}). Setting:
\begin{equation}\label{buko10}
    M(s,t)\,=\, \mathfrak{M}(\mathfrak{v},\mathfrak{w})\equiv \mathfrak{M}\quad;\quad
    \Phi (s,t)\,=\,\mathfrak{F}(\mathfrak{v},\mathfrak{w})\equiv \mathfrak{F}\quad ;\quad
\Omega \, \equiv \,4 (-9+\mathfrak{v}) \mathfrak{v}+12 \mathfrak{v}
\mathfrak{w}+9 \mathfrak{w}^2
\end{equation}
we easily derive that $ds^2_{ort}$ takes the form (\ref{buko2})
with the following matrix $\mathcal{G}_{ij}$:
\begin{eqnarray}\label{buko11}
  \mathcal{G}_{11} &=& -\frac{\mathfrak{v}}{\mathfrak{u}^2-2 \mathfrak{u} \mathfrak{v}}
  \nonumber\\
  \mathcal{G}_{12} &=& \frac{1}{\mathfrak{u}-2 \mathfrak{v}}\nonumber\\
  \mathcal{G}_{13} &=& 0 \nonumber\\
  \mathcal{G}_{22} &=& \frac{1}{9} \left[-\frac{9 \mathfrak{u}}{\mathfrak{u}
   \mathfrak{v}-2 \mathfrak{v}^2}+\frac{\mathfrak{F}^2 \left(-2
   \mathfrak{v}-3 \mathfrak{w}+\sqrt{\Omega }+9\right)^2}{\Omega
   }+\frac{\mathfrak{M}^2 \left(2 \mathfrak{v}+3
   \mathfrak{w}+\sqrt{\Omega }-9\right)^2}{\Omega }\right] \nonumber\\
  \mathcal{G}_{23} &=& \frac{1}{6 \Omega }\left\{
 \mathfrak{F}^2 \left[8 \mathfrak{v}^2+\mathfrak{v} \left(24
   \mathfrak{w}-4 \sqrt{\Omega }-54\right)-3 (2 \mathfrak{w}-3)
   \left(\sqrt{\Omega }-3
   \mathfrak{w}\right)\right] \right.\nonumber\\
   &&\left.+\mathfrak{M}^2 \left[8 \mathfrak{v}^2+\mathfrak{v}
   \left(24 \mathfrak{w}+4 \sqrt{\Omega }-54\right)
   +3 (2 \mathfrak{w}-3) \left(3 \mathfrak{w}+\sqrt{\Omega }\right)\right] \right\} \nonumber\\
  \mathcal{G}_{33} &=& \frac{1}{16 \Omega
   }\left[ 4 \mathfrak{F}^2 \left(-2 \mathfrak{v}-3
   \mathfrak{w}+\sqrt{\Omega }\right)^2+4 \mathfrak{M}^2 \left(2
   \mathfrak{v}+3 \mathfrak{w}+\sqrt{\Omega }\right)^2\right]
   \end{eqnarray}
It remains to be seen if we are able to retrieve the symplectic
potential from which the above matrix is obtained through double
derivatives.
\par
With some integrations and some educated guesses we find that the
form (\ref{buko11}) of the matrix can be reproduced if we write the
symplectic potential as follows:
\begin{equation}\label{buko12}
    G(\mathfrak{u},\mathfrak{v},\mathfrak{w}) \, = \,
    G_0(\mathfrak{u},\mathfrak{v})+\mathcal{G}(\mathfrak{v},\mathfrak{w})
\end{equation}
where
\begin{equation}\label{buko13}
    G_0(\mathfrak{u},\mathfrak{v})=\frac{1}{2} \left(-\mathfrak{u}+\mathfrak{u} \log
    \mathfrak{u}\right)-\mathfrak{v} \log \mathfrak{v}
    +\left(-\frac{\mathfrak{u}}{2}+\mathfrak{v}\right)
\log(-\mathfrak{u}+2 \mathfrak{v})
\end{equation}
and where $\mathcal{G}(\mathfrak{v},\mathfrak{w})$ is some 
function of the two fibre coordinates
$\mathfrak{v}$,$\mathfrak{w}$ only. With this choice the matrix
$G_{\text{ij}}$ becomes:
\begin{equation}\label{buko14}
    G_{ij}=\left(
\begin{array}{ccc}
 -\frac{\mathfrak{v}}{\mathfrak{u}^2-2 \mathfrak{u} \mathfrak{v}} &
 \frac{1}{\mathfrak{u}-2 \mathfrak{v}} & 0 \\
 \frac{1}{\mathfrak{u}-2 \mathfrak{v}} &
  -\frac{\mathfrak{u}}{\mathfrak{u} \mathfrak{v}-2 \mathfrak{v}^2}+\mathcal{G}^{(2,0)}(\mathfrak{v},\mathfrak{w})
& \mathcal{G}^{(1,1)}(\mathfrak{v},\mathfrak{w}) \\
 0 & \mathcal{G}^{(1,1)}(\mathfrak{v},\mathfrak{w}) &
 \mathcal{G}^{(0,2)}(\mathfrak{v},\mathfrak{w}) \\
\end{array}
\right)
\end{equation}
and the full-fledged expression of the line element can be obtained
by substitution. Comparing the obtained result with
eq.~(\ref{buko11}) we easily see that the functions
$\mathfrak{M}(\mathfrak{v}, \mathfrak{w})=M(s,t)$ and
$\mathfrak{F}(\mathfrak{v}, \mathfrak{w}) = \Phi(s,t)$ can be
expressed in terms of the derivatives
$\mathcal{G}^{(2,0)}(\mathfrak{v},\mathfrak{w})$,
$\mathcal{G}^{(0,2)}(\mathfrak{v},\mathfrak{w})$, but in order to
avoid other functions we get a second order differential constraint
on the symplectic potential $\mathcal{G}(\mathfrak{v},
\mathfrak{w})$ that relates its mixed derivatives to
$\mathcal{G}^{(2,0)}(\mathfrak{v},\mathfrak{w})$,
$\mathcal{G}^{(0,2)}(\mathfrak{v},\mathfrak{w})$. This differential
is expressed in a simpler way by means of the original coordinates
$s,t$. We shall presently derive it. We anticipate that its solution
very strongly limits the possibilities so that it has to be
discarded. In other words we have to accept a generic function
$\mathcal{G}(\mathfrak{v}, \mathfrak{w})$ and try to match it with
the boundary conditions on the exceptional divisor.
\subsection{Orthotoric separation of variables and the symplectic
potential} In order to compare the generic metric in symplectic
formalism provided by the symplectic potential displayed in eq.s
(\ref{buko12}), (\ref{buko13}) with the following two-function
metric\footnote{At this level we do not require $M(s,t)$ and
$\Phi(s,t)$ to have the specific form of eq.~(\ref{duefunzie})}:
\begin{eqnarray}
\label{scaramanzia}
   \text{ds}^2_{2fun} &=& \frac{1}{4} s t \left(d\phi ^2 \sin ^2\theta+d\theta ^2\right)+\frac{1}{M(s,t)^2}
  \left[d\chi -\frac{1}{2} t \left(d\phi  (1-\cos (\theta
   ))+d\tau -\frac{2 d\chi
   }{3}\right)\right]^2\nonumber\\
   &&+\frac{1}{\Phi (s,t)^2}{\left[d\chi -\frac{1}{2} s \left(d\phi
   (1-\cos{\theta})+d\tau -\frac{2 d\chi
   }{3}\right)\right]^2}+\frac{1}{4} dt^2 \Phi(s,t)^2+\frac{1}{4} ds^2
   M(s,t)^2\nonumber\\
\end{eqnarray}
 we make the following steps. First we regard
the function $\mathcal{G}(\mathfrak{v},\mathfrak{w})$ as a function
only of $t$ and $s$, as it is evident from the transformation rule
(\ref{buko9}), and we write: $\mathcal{G}(\mathfrak{v},\mathfrak{w})
\, \equiv \, \Gamma(t,s) $. By means of the transformation
(\ref{buko9}) we can rewrite the generic metric (\ref{buko2})
produced by the symplectic potential (\ref{buko12}-\ref{buko13}) in
terms of the variables $s,t$, instead of
$\mathfrak{v},\mathfrak{w}$. The result coincides with $\text{ds}^2_{2fun}$
as given in eq.~(\ref{scaramanzia}) if the following conditions hold
true:
\begin{eqnarray}\label{santecchio}
  &&\frac{\partial^2}{\partial s^2}\Gamma(t,s) \,=\, \frac{1}{4}  \,
  M(s,t)^2 \quad ; \quad
  \frac{\partial^2}{\partial t^2}\Gamma(t,s) \,=\, \frac{1}{4}  \, \Phi(s,t)^2 \nonumber \\
  &&\frac{\partial}{\partial s}\Gamma(t,s)\, - \,\frac{\partial}{\partial t}\Gamma(t,s)
  +(s-t)\frac{\partial^2}{\partial s\,\partial t}\Gamma(t,s) \,=\, 0
\end{eqnarray}
The first two equations in (\ref{santecchio}) just provide the
identification of the two functions $M(s,t)$ and $\Phi(s,t)$ in
terms of second order derivatives of the symplectic potential. On
the other hand the last equation of (\ref{santecchio}) is a very
strong constraint on the function $\Gamma(t,s)$ which severely
restricts the available choices of $\Gamma(t,s)$.
\subsection{The symplectic potential of the Ricci-flat orthotoric metric
on ${\rm tot}  K_{{\mathbb WP}[112]}$} In the case of the canonical
bundle ${\rm tot} K_{{\mathbb WP}[112]}$, whose Ricci--flat metric
is given by eq.~(\ref{iniz11}),  eq.s (\ref{iniz10}) imply
\begin{eqnarray}\label{brugola}
    \Gamma ^{(0,2)}(t,s) & = & \frac{(s-t) \Upsilon (s)^2}{4
    (s+3)} \, = \, -\frac{s (s-t)}{4 (s+3) \left(\frac{2
    s^2}{3}-s+3\right)}\quad \Rightarrow \quad \Upsilon (s)^2 = \frac{s }{\left(\frac{2
    s^2}{3}-s+3\right)}\nonumber\\
    \quad \Gamma ^{(2,0)}(t,s) & = &\frac{1}{4} P(t)^2 (t-s)\, =
    \, \frac{s-t}{4 \left(\frac{2 t^2}{3}+t\right)} \quad \Rightarrow
    \quad P(t)^2 \, = \, -\frac{1}{\frac{2}{3}\, t\,\left(t+3\right)}
\end{eqnarray}
By means of two double integrations and modulo linear functions in
$s,t$ (they are irrelevant for the metric) we determine the explicit
form of the potential $\Gamma (t,s)$:
\begin{eqnarray}\label{gazpacho}
  \Gamma_{{{\mathbb WP}[112]}} (t,s) &=& \frac{1}{224}
  \left\{-7\left[(3 t-st+3s)
  \log (2 s^2-3s+9)+2 (s+3) (t+3) \log{(s+3)}-8 s t \log{t}\right.\right.\nonumber\\
   &&\left.\left.+2 (2 s+3) (2
   t+3) \log (2 t+3)\right]-6 \sqrt{7} (s t+s+t-6) \arctan\frac{3-4 s}{3 \sqrt{7}}\right\}
\end{eqnarray}
The function $\Gamma_{{{\mathbb WP}[112]}}(t,s)$ satisfies by
construction the differential constraint encoded in the third of
eq.s~(\ref{brugola}). Using the transformation rule (\ref{buko9}) we
can rewrite it as a function of the symplectic variables
$\mathfrak{v},\mathfrak{w}$. In this way we arrive at the following
symplectic potential where we have used the liberty of  adding
linear functions of $\mathfrak{v}$ or  $\mathfrak{w}$ to obtain the
most convenient form of its reduction to the exceptional divisor,
located at $\mathfrak{w}=\ft 32$. The function
\begin{eqnarray}
\label{miracolina}
  \mathcal{G}_{{{\mathbb WP}[112]}}\left(\mathfrak{v},\mathfrak{w}\right)  &=&
  \frac{1}{224} \left\{ 7 \left[6 (2 \mathfrak{w}-3) \log
   \left(-\sqrt{(2 \mathfrak{v}+3 \mathfrak{w})^2-36
   \mathfrak{v}}-2 \mathfrak{v}-3 \mathfrak{w}+9\right)\right.\right.\nonumber\\
   &&\left.\left.+16
   \mathfrak{v} \log\left(\sqrt{(2 \mathfrak{v}+3
   \mathfrak{w})^2-36 \mathfrak{v}}-2 \mathfrak{v}-3
   \mathfrak{w}\right)^2\right.\right.\nonumber\\
   &&\left.\left.-2 (8 \mathfrak{v}-12
   \mathfrak{w}+9) \log \left(\sqrt{(2 \mathfrak{v}+3
   \mathfrak{w})^2-36 \mathfrak{v}}-2 \mathfrak{v}-3
   \mathfrak{w}+\frac{9}{2}\right)\right.\right.\nonumber\\
   &&\left.\left.+2 (4 \mathfrak{v}+3
   \mathfrak{w}) \log \left(\frac{1}{567} \left[4 \sqrt{(2
   \mathfrak{v}+3 \mathfrak{w})^2-36 \mathfrak{v}}+8
   \mathfrak{v}+12 \mathfrak{w}+9\right]^2+1\right)\right]\right.\nonumber\\
   &&\left.-4
   \sqrt{7} (4 \mathfrak{v}-3 (\mathfrak{w}+3))\, \text{arctan}\left(\frac{4 \sqrt{(2 \mathfrak{v}+3 \mathfrak{w})^2-36
   \mathfrak{v}}+8 \mathfrak{v}+12 \mathfrak{w}+9}{9
   \sqrt{7}}\right)\right.\nonumber\\
   &&\left.-(8 \mathfrak{v}+9) \log
   \frac{34359738368}{823543}+2 \sqrt{7} (8
   \mathfrak{v}-27)\, \arctan\frac{5}{\sqrt{7}}\right\}
\end{eqnarray}
is expressed in terms of elementary transcendental functions, yet it has the
remarkable property of satisfying the 
Monge-Amp\`{e}re equation for Ricci-flatness, so that it may be called ``the miraculous function".
On the exceptional divisor
it  reduces  to
\begin{equation}
\label{bondorino}
  \mathcal{D}_{{{\mathbb WP}[112]}}(\mathfrak{v})\,  \equiv\,
  \mathcal{G}_{{{\mathbb WP}[112]}}\left(\mathfrak{v},\ft 32\right)\, = \,
  \frac{1}{16} \left[ 8 \mathfrak{v} \log \left(16
   \mathfrak{v}^2\right)+(9-8 \mathfrak{v}) \log
   \left(\frac{9}{2}-4 \mathfrak{v}\right)\right].
\end{equation}

\section{K\"ahler metrics on Hirzebruch surfaces and their canonical bundles}
\label{ciurlebba} For the case of the canonical bundle on
$\mathbb{F}_2$, which is  the complete resolution of the
$\mathbb{C}^3/\mathbb{Z}_4$ singularity, we have additional
information that is relevant and inspiring for the general case.
\par
Let us summarize the main points. According to the results of
\cite{noietmarcovaldo} there is a well adapted system of complex
coordinates that arise from the toric analysis of
$\mathbb{C}^3/\mathbb{Z}_4$ and of its resolution. These coordinates
are named as follows: $
  z_i \, = \, \left\{  u, \, v, \, w \right\}$
and are defined on a dense open chart reaching all components of the
exceptional divisor. Their interpretation was already anticipated in
eq. (\ref{iniz13}) and it is the following. The coordinate $w$ spans
the fibers in the canonical bundle
$Y\,\stackrel{\pi}{\longrightarrow} \, \mathbb{F}_2 $ while $u,v$
span a dense open chart for the base manifold (\textit{i.e.} the
compact component $\mathbb{F}_2$ of the exceptional divisor
$\mathcal{ED}$). In particular since $\mathbb{F}_2$ is a
$\mathbb{P}^1$ bundle over $\mathbb{P}^1$, namely $ \mathbb{F}_2
\stackrel{\pi}{\longrightarrow} \, \mathbb{P}^1$, the coordinate $u$
is a standard Fubini-Study coordinate for the base $\mathbb{P}^1$
while $v$ spans a dense open chart of the fibre $\mathbb{P}^1$.
This set of coordinates can be used for any $\mathbb{F}_n$
Hirzebruch surface with $n \geq 1$. The action of the isometry group
(\ref{gisometry}) on these coordinates was described in
\cite{noietmarcovaldo} and it is as follows:
\begin{eqnarray}
   \forall\, \mathbf{g} &=& \left(
          \begin{array}{cc}
            a & b \\
            c & d \\
          \end{array}
        \right) \, \in \, \mathrm{SU(2)}\quad : \quad
    \mathbf{g}\left(u,v,w\right) \, = \,\left(\frac{a \, u + b}{c\, u +
    d}, \quad v \, \left(c \,u+d\right)^{n},\quad w\right)
    \nonumber\\
\forall\, \mathbf{g} &=& \exp[{\rm i}\,\theta_1] \, \in \, \mathrm{U(1)}_v
\quad\quad : \quad \mathbf{g}\left(u,v,w\right) \, = \,\left(u,
\quad \exp[{\rm  i}\,\theta_1]\, v \,\quad w\right)
\nonumber\\
\forall\, \mathbf{g} &=& \exp[{\rm i}\,\theta_2] \, \in \, \mathrm{U(1)}_w
\quad \quad: \quad \mathbf{g}\left(u,v,w\right) \, = \,\left(u,
\quad \, v,\quad \exp[{\rm i}\,\theta_2]\, w\right) \label{ciabattabuona}
\end{eqnarray}
The above explicit action of the isometry group on the $u,v,w$
coordinates suggests the use of an invariant real combination \begin{equation}\label{invarpi}
    \varpi_n \, \equiv \, \left(1+\mid u\mid^2\right)^n \, \mid
    v\mid^2
\end{equation}
and the assumption that the K\"ahler potential
$\mathcal{K}_{\mathbb{F}_n}$ of the K\"ahler metric
$\mathbf{g}_{\mathbb{F}_n}$ should be a function (up to trivial
terms $\mathrm{Re} f(z)$) only of $\varpi_n$:
\begin{equation}\label{lamasco}
    \mathcal{K}_{\mathbb{F}_n}\, = \, G_n(\varpi_n)
\end{equation}
The function $G_n(\varpi_n)$ should also depend on two parameters
(we name them $\ell$,$\alpha$) which are associated to the volumes
of the two homology cycles of $\mathbb{F}_n$, respectively named
$C_1$ and $C_2$ that also form a basis for the homology group of the
total space $Y$, namely the canonical bundle on $\mathbb{F}_n$.
Indeed the homology of $Y$ coincides with the homology of the base
manifold $\mathbb{F}_n$.
\par
Introducing the K\"ahler two form:
\begin{equation}\label{pantufla}
    \mathbf{K}_{\mathbb{F}_n} \, \equiv \, \frac{{\rm i} }{2\pi} \partial
    \bar{\partial} \mathcal{K}_{\mathbb{F}_n}
\end{equation}
we need to find:
\begin{eqnarray}
  &&\int_{C_1}\mathbf{K}_{\mathbb{F}_n} \,=\,  \frac{9}{16} \, \alpha
  \,\ell \quad ; \quad
  \int_{C_2}\mathbf{K}_{\mathbb{F}_n} \,=\, \frac{9}{16} \, (2 +\alpha ) \, \ell
  \label{palliativo}
\end{eqnarray}
where $\ell$ is a dimensionful parameter providing the scale and
$\alpha$ is some dimensionless parameter parameterizing the ratio
between the two volumes. The two toric cycles $C_{1,2}$ are
respectively defined by the following two equations:
\begin{equation}\label{frescone}
    C_1 \, \Leftrightarrow \, v=0 \quad ; \quad C_2 \, \Leftrightarrow \, u=0
\end{equation}
As pointed out in \cite{noietmarcovaldo}, in addition to the above
two properties of the K\"ahler form, if we consider the Ricci
two-form of the K\"ahler metric on $\mathbb{F}_n$
\begin{equation}\label{riccionetto}
    \mathbf{Ric}_{\mathbb{F}_n} \, = \, \frac{{\rm i} }{2\pi} \partial
    \bar{\partial} \log\left[\mbox{det}\left(\mathbf{g}^{\mathbb{F}_n}\right)
    \right]\quad ; \quad \mathbf{g}^{\mathbb{F}_n}_{ij^\star} \, =
    \, \partial_i\partial_{j^\star}\,\mathbf{K}_{\mathbb{F}_n}
    \quad \quad i=1,2 \quad j^\star = 1^\star,2^\star
\end{equation}
we must find:
\begin{eqnarray}
&&\int_{C_1}\mathbf{Ric}_{\mathbb{F}_n} \, =\,  2-n \quad ; \quad
  \int_{C_2}\mathbf{Ric}_{\mathbb{F}_n} \, =\,  2
  \label{curativo}
\end{eqnarray}
It appears that eq.s (\ref{palliativo}-\ref{curativo}) are strong
constraints on the function $G_n(\varpi_n)$. It is interesting to
see how they are realized in the metric on $\mathbb{F}_2$ obtained
from the Kronheimer construction. We will show this below.
\subsection{The metric on $\mathbb{F}_2$ induced by the Kronheimer
construction}\label{KronoinducoF2} In \cite{noietmarcovaldo}, relying on the
Kronheimer construction, we have
constructed an analytically defined K\"ahler metric on the total
space of the canonical bundle of $\mathbb{F}_2$. The K\"ahler potential has only an implicit
definition as the largest real root of a sextic equation. Yet its
reduction to the compact exceptional divisor, which is indeed  the
2nd Hirzebruch surface, is explicit and the K\"ahler potential of
this metric can be exhibited in closed analytic form. We think that
this information is very important  for the comparison between the
parameters of the Ricci-flat metric appearing in supergravity with
those emerging in the Kronheimer construction that are the Fayet
Iliopoulos parameters of the dual gauge theory.
\par
Following the chamber structure discussed in \cite{noietmarcovaldo}
we choose the chamber VI defined by the following inequalities on the
three Fayet Iliopoulos parameters $\zeta_{1,2,3}$:
\begin{equation}\label{cameraperdue}
   \zeta _1-\zeta _2-\zeta _3<0 \quad;\quad -\zeta _1+\zeta _2-\zeta _3>0 \quad;\quad -
   \zeta _1-\zeta _2+\zeta _3<0
\end{equation}
and chamber VIII, defined instead by the following ones:
\begin{equation}\label{camerapertre}
   \zeta _1-\zeta _2-\zeta _3<0 \quad;\quad -\zeta _1+\zeta _2-\zeta _3<0 \quad;\quad -
   \zeta _1-\zeta _2+\zeta _3<0
\end{equation}
Inside those two chambers we make the choice:
\begin{equation}\label{cardellino}
    \zeta_1\, = \, \zeta_3 \, = \,  r \quad ; \quad \zeta_2
    \, = \, (2+\alpha) r \quad ; \quad \quad r>0
\end{equation}
For $\alpha >0$ we are in chamber VI, while for $\alpha <0$ we are
in chamber VIII. For $\alpha=0$ we are instead on the wall where the
non singular variety:
\begin{equation}
Y \, \equiv \,
{\rm tot}  K_{{\mathbb F}_2}
\label{Ylisciadefi}
\end{equation}
degenerates in
\begin{equation} Y_3 \, \equiv \,  {\rm tot}  K_{{\mathbb WP}[112]}
\label{Y3defi}
\end{equation}
denoting by ${\rm tot} K_{\mathcal{M}}$ the total space of the
canonical bundle of a K\"ahler manifold (or orbifold) $\mathcal{M}$.
\par
The solution of the moment map equations for the two independent
moment maps reduced to the exceptional divisor by performing the
limit $w\to 0$ is the following one:
\begin{equation}
 T_1=T_3 \, =\, {\sqrt{\frac{\sqrt{\alpha ^2+6 \alpha  \varpi
+\varpi
   (\varpi +8)}+\alpha +\varpi }{2(\alpha +2) \sqrt{{\varpi
   }/{v \bar{v}}}}}} \, ; \,
  T_2 \, =\, \frac{\sqrt{\alpha ^2+6 \alpha  \varpi +\varpi ^2+8 \varpi }+3
   \alpha +\varpi +4}{2 \alpha ^2+6 \alpha +4}
\end{equation}
The complete K\"ahler potential of the quotient is made of two
addends, the pull-back on the constrained surface of the K\"ahler
potential of the flat ambient metric plus the logarithmic term:
\begin{equation}\label{carolingio}
    \mathcal{K}_{quotient} \, = \, \mathcal{K}_{0} \, + \, \underbrace{\zeta_I
    \mathfrak{C}^{IJ} \, \log T_J}_{\mathcal{K}_{log}} \quad ; \quad  \mathfrak{C}^{IJ}
    \, = \, \left(
              \begin{array}{ccc}
                2 & -1 & 0 \\
                -1 & 2 & -1 \\
                0 & -1 & 2 \\
              \end{array}
            \right)
\end{equation}
In the present case we explicitly find:
\begin{equation}\label{carondimonio}
  \mathcal{K}_{0} \, = \,
 2 \frac{\alpha  \left(\sqrt{\alpha ^2+6 \alpha  \varpi
   +\varpi  (\varpi +8)}+2 \varpi +1\right)+\sqrt{\alpha ^2+6
   \alpha  \varpi +\varpi  (\varpi +8)}+\alpha ^2+3 \varpi
   }{\sqrt{\alpha ^2+6 \alpha  \varpi +\varpi  (\varpi
   +8)}+\alpha +\varpi }
\end{equation}
and
\begin{eqnarray}\label{bartolomeus}
    \mathcal{K}_{log}& = & 2 (\alpha +1) \log\frac{\sqrt{\alpha ^2
    +6 \alpha  \varpi
   +\varpi ^2+8 \varpi }+3 \alpha +\varpi +4}{2 \alpha ^2+6 \alpha
   +4}
   -2 \alpha  \log {\sqrt{\frac{\sqrt{\alpha
   ^2+6 \alpha  \varpi +\varpi  (\varpi +8)}+\alpha +\varpi
   }{2(\alpha +2) \sqrt{{\varpi }/{v\bar{v}}}}}}\nonumber\\
\end{eqnarray}
By explicit calculation we were able to verify that the K\"ahler
potential of the quotient $\mathcal{K}_{quotient}$ yields a metric
satisfying all the constraints (\ref{palliativo}-\ref{curativo}). We
show this in section \ref{F2Kalrenormo}.
\section{Reduction to the exceptional divisor}
\label{genreductia} In this section we consider the reduction to the
exceptional divisor for a generic metric of the class described in
section \ref{duefunziesecche}, emphasizing that the K\"ahler metric
induced on the divisor is completely determined by the real function
$P(t)$ of the real variable $t$. We carefully consider what are the
differential constraints on such a function required by the topology
and complex structure of the second Hirzebruch surface
$\mathbb{F}_2$ showing that they are all met by the $P(t)$ function
that one obtains by localizing the generalized Kronheimer
construction of the $ \mathbb{C}^3/\mathbb{Z}_4$ singularity
resolution on the exceptional divisor.
\subsection{The reduction}
The reduction to the exceptional divisor is obtained in the K\"ahler
form and in the metric by setting $s = -3$. The K\"ahler form on the
 divisor is the following one
\begin{equation}\label{divisor1}
    \mathbb{K}_{\mathcal{ED}}=\frac{1}{12} \left(-9 \,t \, \sin{\theta}
{d\theta}\wedge {d\phi}-9 \,dt\wedge
{d\tau}+9 \, (\cos{\theta} -1)\, dt\wedge d\phi\right)
\end{equation}
while the metric is the following one:
\begin{equation}\label{divisor2}
    \text{ds}_{\mathcal{ED}}^2=\,-\frac{3\,t}{4}({d\theta}^2 + \sin^2\theta {d\phi}^2)+
    \frac{1}{4} (t+3) P(t)^2 {dt}^2 +\frac{9 [{d\tau}+(1-\cos\theta){d\phi} ]^2}{4 (t+3) P(t)^2}
\end{equation}
and it is completely determined by the function $P(t)$. For the
choice:
\begin{equation}\label{divisor3}
    P(t)=\left(-\frac{2}{3}t^2-t\right)^{-\frac{1}{2}}
\end{equation}
it is the metric on the orbifold ${\mathbb WP}[112]$
while for other choices of $P(t)$, obtainable from the Kronheimer
construction, $\text{ds}_{\mathcal{ED}}^2$ can indeed be a good
K\"ahler metric on the second Hirzebruch surface $\mathbb{F}_2$.
\par
From eq.~(\ref{divisor1}) specifying the K\"ahler 2-form of the
exceptional divisor and eq.~(\ref{divisor2}) providing its K\"ahler
metric, we immediately work out also the complex structure tensor
that has the following appearence:
\begin{equation}\label{lomellina}
    \mathbb{J}_{\mathcal{ED}}\, = \,\left(
\begin{array}{cccc}
 0 & 0 & 0 & -\frac{1}{3} (t+3) P(t)^2 \\
 0 & 0 & \csc\theta & -\tan{\theta\over 2} \\
 {6 \sin^2{\theta\over 2} \over (t+3) P(t)^2} &
   -\sin\theta & 0 & 0 \\
 \frac{3}{(t+3) P(t)^2} & 0 & 0 & 0 \\
\end{array}
\right)
\end{equation}
\subsection{Topology and the functions of the $t$ coordinate}
We have two important informations on the topology of
$\mathbb{F}_2$, which provide an extremely selective test in order
to know whether a certain metric is indeed defined on
$\mathbb{F}_2$ or on some different twofold, may be degenerate.
The tests are related with the integrals of the K\"ahler 2-form
$\mathbb{K}$ and of the Ricci 2-form $\mathbb{R}$ic on the two toric
curves $C_{1,2}$ respectively defined by the vanishing of either
coordinate $(u,v)$
\begin{equation}\label{divisor4}
    C_1=\{v=0\}\quad ; \quad C_2=\{u=0\}
\end{equation}
Indeed, as we illustrated in section \ref{ciurlebba} we must find
\begin{equation}\label{divisor5}
 \int _{C_1}\mathbb{K}\neq \infty \quad ;\quad\int
_{C_2}\mathbb{K}\neq \infty \quad ; \quad \int _{C_1}\mathbb{K}\neq
0\quad ;\quad \int _{C_2}\mathbb{K}\neq 0 \quad ; \quad
\int_{C_1}\text{$\mathbb{R}$ic}=0 \quad ; \quad \int
_{C_2}\text{$\mathbb{R}$ic}=2
\end{equation}
The explicit reduction of the K\"ahler form
$\mathbb{K}_{\mathbb{F}_2}$ to the two cycles $C_1$ and $C_2$ is
very simple when $\mathbb{K}_{\mathbb{F}_2}$ is written in the basis
of the real coordinates ($t$,$\theta $,$\tau $,$\phi $). Indeed in
order to set $v=0$ we have just to look for the zeros of the above
defined function $H(t)$ that depends by integration from $P(t)$. Let
us suppose that $H(-|t_{max}| )=0$. We obtain the reduction of the
K\"ahler form to the cycle $C_1$ by setting $t = - |t_{max}| $ =
const $<$ 0, while we get the reduction to the cycle $C_2$ by
setting $\theta $ = 0.
\begin{equation}\label{divisor6}
  \mathbb{K}\mid _{C_1} \,=\, \frac{3}{4}\,|t_{max}|\, \sin{\theta}\,
  {d\theta}
  \wedge
{d\phi} \quad ; \quad
 \mathbb{K}\mid _{C_2} \,=\, -\frac{3}{4}{dt}\wedge {d\tau}
\end{equation}
Hence we see that in order to get $\mathbb{F}_2$ as exceptional
divisor we need two conditions, that are necessary, although not
sufficient.
\begin{enumerate}
  \item $|t_{max}| $ $\neq $ 0
  \item the range of
  the coordinate $t$ must be finite
  [$-|t_{max}| $,$-|t_{min}|$] in order to get a finite size for the cycle $C_2$
\end{enumerate}
If the zero of the function $H(t)$ is at t=0 we immediately know
that there is a degeneration and this is indeed the case of
${\mathbb WP}[112]$.
\par
If we integrate the complex structure of the exceptional divisor
displayed in eq.~(\ref{lomellina}) with the same method we used
for the whole 6-dimensional space,  we find that the coordinate $u$
is exactly the same as in eq.~(\ref{iniz21}), while for $v$ we find:
\begin{equation}\label{vigevano}
   v \, = \,  H(t) \cos^2\frac{\theta}{2}e^{{\rm i} \tau }
\end{equation}
Comparison with the result for $v$ in the entire space
(eq.s~(\ref{iniz22}-\ref{iniz23})) tells us that the function
$\Psi(s)$ must be finite and non vanishing at $s = -3$ in order to
have a consistent reduction to the divisor:
\begin{equation}\label{cornicione}
    \Psi(-3) = 1  \quad ; \quad   \Phi(-3) = 0
\end{equation}
The normalization $\Psi(-3)$ = 1  can always be obtained by an
irrelevant rescaling in the definition of $v$ if $-3$ is not a zero
of $\Psi(s)$  while it must be a zero of $\Phi(s)$.
\subsubsection{Interpretation of the function $H(t)$}
From the explicit integration of the complex structure we obtain a
very important interpretation of the function $H(t)$ in relation
with the complex K\"ahler geometry of the exceptional divisor.
 Since the K\"ahler metric on this two-fold has
isometry $\mathrm{SU(2)\times U(1)}$, $\mathrm{SU(2)}$ acting on the
$u$ variable by linear fractional transformation and on $v$ by
multiplication with the $u$-compensator $(c u+d)^2$, as described in
eq.s (\ref{ciabattabuona}),  the K\"ahler potential $\mathcal{K}$
can be a function only of the  invariant combination $\varpi\equiv
\varpi_2$ defined in eq.~(\ref{invarpi}). Relying on the
representation of $u$ and $v$ derived from the integration of the
complex structure we easily obtain:
\begin{equation}\label{cannolosiculo}
   \varpi \, =\,  \cos ^4\frac{\theta }{2}\left(\tan
   ^2\frac{\theta}{2}+1\right)^2 H(t)^2 \, = \, H(t)^2
\end{equation}
It follows that:
\begin{equation}\label{gnoccodipatata}
    t\, = \, H^{-1}(\sqrt{\varpi})
\end{equation}
where $H^{-1}$ denotes the inverse function. Since the range of
$\sqrt{\varpi}$ is $[0,\infty]$, it is necessary that the inverse
function $H^{-1}$ maps the semi-infinite interval $[0,\infty]$
in a finite one $[-|t_{max}|,-|t_{min}|]$ defined by:
\begin{equation}\label{carriolino}
    -|t_{max}| \, = \, \lim_{\varpi \to 0}\, H^{-1}(t) \quad ; \quad -|t_{min}| \, = \,
    \lim_{\varpi \to \infty}\, H^{-1}(t)
\end{equation}
\subsubsection{Topological constraints on the function $P(t)$}
Given
the above  topology results characterizing the second Hirzebruch
surface and considering the metric of the divisor as given in
eq.~(\ref{divisor2}) and its K\"ahler form (\ref{divisor1}) we
immediately obtain the conditions on the function $P(t)$. Indeed, while
calculating the Ricci form we can specify integral differential
conditions on $P(t)$ from the values of its periods mentioned above.
We know the explicit form of the complex structure on the
exceptional divisor that is obtained by reduction to $s = -3$ of the
complex structure pertaining the full 6-dimensional manifold
$\mathcal{M}_6$. The complex structure of the exceptional divisor
was displayed in eq.~(\ref{lomellina}). The Ricci form can be
calculated by setting its antisymmetric components equal to $\text{$\mathbb{R}$ic}_{ij}$ = $J_i^{\text{
  }k}$$R_{\text{kj}}$ where $R_{kj}$ is the standard Ricci tensor.
  In this way we obtain the following general result that exclusively
depends on the function $P(t)$:
\begin{equation}\label{divisor8}
\text{$\mathbb{R}$ic}_{\mathcal{ED}}\,=\,\mathfrak{A}(t)\sin{\theta}\,
{d\theta }\wedge {d\phi}+\mathfrak{B}(t)\, \sin ^2\frac{\theta
}{2}\, {dt}\wedge {d\phi}\,+\, \mathfrak{C}(t){dt}\wedge {d\tau}
\end{equation}
where $\mathfrak{A}(t),\mathfrak{B}(t),\mathfrak{C}(t)$ are
functions of the $t$-variable expressed as rational functions of
$P(t)$ and its first and second derivative with simple $t$-dependent
coefficient. We do not write them explicitly for shortness. Then the
Ricci 2-form can be easily localized on the two cycles $C_1$ and
$C_2$, yielding:
\begin{equation}\label{divisor9}
\text{$\mathbb{R}$ic}|_{C_1} \, =\,
\mathfrak{A}(-|t_{max}|)\,\sin{\theta}\, d\theta\wedge d\phi\quad
;\quad
  \text{$\mathbb{R}$ic}|_{C_2}\, =\,\mathfrak{C}(t)\, dt \wedge d\tau
\end{equation}
Hence, in order to realize the second Hirzebruch surface not only
the range of $t$ must have finite extrema [$-|t_{max}| $,$-|t_{min}|
$] but we should also have:
\begin{equation}\label{divisor11}
    \mathfrak{A}(-|t_{max}|)=0 \quad ; \quad \int_{-|t_{min}|}^{-|t_{max}| } \mathfrak{C}(t) \, dt = 2
\end{equation}
\subsubsection{The relation between the function $P(t)$ and the
K\"ahler potential $\mathcal{K}(\varpi)$ of the exceptional
divisor}

Our goal is that of determining a Ricci-flat metric on the canonical
bundle ${\rm tot} K_{\mathbb{F}_2}$, starting from a given
\textit{bona fide} K\"ahler  metric on the second Hirzebruch
surface, described in terms of the real variables $t$, $\theta$, $\tau$, $\phi$. In the complex description, any K\"ahler metric is
determined by a suitable K\"ahler potential; given the isometries
and their realization on the chosen complex coordinates $u,v$, the
K\"ahler potential for the $\mathbb{F}_2$ surface  is a real
function  of the invariant combination $\varpi $ defined in
eq.~(\ref{invarpi}) which we generically denote $\mathcal{K}$($\varpi
$). Therefore it is important to determine the relation between the
real variables and the standard complex ones at the same time with
the relation between the K\"ahler potential $\mathcal{K}$($\varpi $)
and the function $P(t)$ which determines the metric in the real
variables. In this respect the essential point to be stressed is
that the relation between the real variables and the complex ones is
not universal and fixed once for all, rather it depends on the
choice of the K\"ahler potential or viceversa of the function
$P(t)$. Hence it is convenient to introduce a name for the inverse
function:
\begin{equation}\label{divisor13}
    H^{-1} \left(\sqrt{\varpi }\right)  = G_T(\varpi )
\end{equation}
and find its differential relation with the K\"ahler potential which
follows from a comparison between the metric as determined in
complex K\"ahler geometry from $\mathcal{K}$($\varpi $) and as
written in real variables. For convenience we rewrite the general
real form of the metric  on the exceptional divisor in the following
more compact way
\begin{equation}\label{divisor15}
  \text{ds}_{\mathcal{ED}}^2\, = \,  -\frac{3}{4} \left(t\underbrace{
 [{d\theta}^2
 +{d\phi}^2
   \sin ^2\theta]}_{\text{metric on $\mathbb{P}^1$}}+ F(t)\,{dt}^2\,
   +\frac{1}{F(t)}
   {\underbrace{[{d\tau}
   +(1-\cos{\theta}){d\phi}]}_{\text{connection
   on the $\mathrm{U(1)}$ bundle}}}^2\right) \quad ; \quad F(t)\,\equiv \, -\frac{1}{3}(t+3) P(t)^2
\end{equation}
which clearly displays the fibred structure of the exceptional
divisor.
\par
Next we convert the metric in eq.~(\ref{divisor15}) using the
substitution rule
\begin{equation}\label{divisor16}
    t\,=\, G_T(\varpi ),\quad \theta \,
=\, 2\arctan\sqrt{u\bar{u}}, \quad \tau \, =\, -\frac{1}{2} {\rm i}\,
\log\left(\frac{v}{\bar{v}}\right),\quad \phi \,=\,
-\frac{1}{2} {\rm i}\log\left(\frac{u}{\bar{u}}\right)
\end{equation}
In this way we transform the metric (\ref{divisor15}) to the complex
coordinates $u,v$ and we compare it with the generic metric obtained
from a generic K\"ahler potential $\mathcal{K}$($\varpi $). We find
that the two metrics coincide provided the following two conditions
are satisfied:
\begin{equation}\label{divisor17}
  t\, \equiv \, G_T(\varpi )\,=\, -\frac{2}{3} \varpi  \partial _{\varpi
}\mathcal{K}(\varpi )\quad ; \quad P(t)\,= \,\pm \frac{3}{2
\sqrt{-\varpi \left(3-\frac{2}{3} \varpi  \partial _{\varpi
}\mathcal{K}(\varpi )\right) \left(\partial _{\varpi
}\mathcal{K}(\varpi )+\varpi
\partial _{\varpi ,\varpi }\mathcal{K}(\varpi )\right)}}
\end{equation}
Given the K\"ahler potential $\mathcal{K}$($\varpi $), which is
supposed to depend also on a deformation parameter, the above
equation (\ref{divisor17}) allows  to rewrite the same metric in
real coordinates, provided one is able to invert the first formula,
namely, to find $\varpi $ as a function of $t$ and of the
deformation parameter $\alpha $.
\subsection{The Kronheimer K\"ahler potential for the
$\mathbb{F}_2$ surface and its associated $P(t)$ function}
\label{F2Kalrenormo} From the Kronheimer construction of the
$\mathbb{C}^3/\mathbb{Z}_4$ resolution reduced to the exceptional
divisor we have the K\"ahler potential derived in section
\ref{KronoinducoF2}. The result obtained in eq.s
(\ref{carolingio},\ref{carondimonio},\ref{bartolomeus}) can be
summarized writing the following  general form of the K\"ahler
potential:
\begin{eqnarray}
\label{divisor18}
  \mathcal{K}_{\mathbb{F}_2}(\varpi ,\alpha )&=&
   -\frac{9}{16}\left\{-4 (\alpha +1) \log
   \left(\sqrt{\alpha ^2+6 \alpha  \varpi +\varpi ^2+8
\varpi }+3 \alpha +\varpi +4\right)\right. \nonumber\\
&&-\frac{4 \left[\alpha  \left(\sqrt{\alpha ^2+6 \alpha  \varpi
+\varpi (\varpi +8)}+2 \varpi +1\right)+\sqrt{\alpha ^2+6 \alpha
\varpi +\varpi  (\varpi +8)}+\alpha ^2+3 \varpi
\right]}{\sqrt{\alpha ^2+6 \alpha  \varpi
+\varpi  (\varpi +8)}+\alpha +\varpi }\nonumber\\
&&\left. +4 \alpha  \log \sqrt{\frac{\sqrt{\alpha ^2+6 \alpha
\varpi +\varpi  (\varpi +8)}+\alpha +\varpi }{\sqrt{\varpi
}}}+8+16 \log{2}\right\}
\end{eqnarray}
where for additional  convenience we have changed  the overall
normalization of the metric multiplying by a $9/8$   factor,
  have disregarded the irrelevant addends proportional to $\log[v
\,\bar{v}]$ and   have added a convenient constant addend. For
$\alpha $=0 the surface described by the  K\"ahler metric
corresponding to the potential (\ref{divisor18}) degenerates into
the singular ${\mathbb WP}[112]$ while for other values of $\alpha $
such that $0 < \mid \alpha \mid < 1$ we have a metric on a smooth
$\mathbb{F}_2$ surface.
\subsubsection{The degenerate case   ${\mathbb WP}[112]$}
It is interesting to see how we recover the degenerate case
${\mathbb WP}[112]$ from the general case. Setting
$\alpha $=0 we obtain:
\begin{eqnarray}
\label{divisor19}
    \mathcal{K}_{ {\mathbb WP}[112]}(\varpi )& =&\frac{9}{4}
    \left[\frac{3 \varpi +\sqrt{\varpi  (\varpi
   +8)}}{\varpi +\sqrt{\varpi  (\varpi +8)}}+\log\left(\varpi
   +\sqrt{\varpi  (\varpi +8)}\right)-2-4 \log{2}
   \right]\nonumber\\
t&=&-\frac{3 \varpi ^2 \left(\varpi +\sqrt{\varpi  (\varpi
+8)}+8\right)}{\sqrt{\varpi  (\varpi +8)} \left(\varpi +\sqrt{\varpi
(\varpi +8)}\right)^2} \quad \Rightarrow \quad \varpi =\frac{8
t^2}{3 (2 t+3)}
\end{eqnarray}
This implies that the interval [0,$\infty $] of $\varpi$ is mapped
into the interval [0,-3/2] and this suffices to guarantee that the
cycle $C_1$ is contracted to zero as we have already explained.
Finally for the function $P(t)$, using the above general
formulae we get:
\begin{equation}
\label{divisor20}
   P(t)=\sqrt{\frac{-3}{2 t^2+3t}}
\end{equation}
\subsubsection{The smooth $\mathbb{F}_2$ case}
First we can verify that when $\alpha $ is either $-1$ or $-2$, the
surface degenerates, as the metric depends only on the variable
$u$ and no longer on $v$. Using the formula (\ref{divisor17}) we can
calculate $t$ and $P(t)$. We find the following relatively complicated
answer:
\begin{eqnarray}\label{divisor21}
t&=& G_T(\varpi )\,=\, \frac{N_T}{D_T} \nonumber\\
N_T &=& -3 \left\{\alpha ^4+\alpha  \varpi  \left[3 \varpi
   \left(\sqrt{\alpha ^2+6 \alpha  \varpi +\varpi  (\varpi
   +8)}+16\right)+8 \sqrt{\alpha ^2+6 \alpha  \varpi +\varpi
   (\varpi +8)}+3 \varpi ^2\right]\right.\nonumber\\
   &&\left.+4 \varpi ^2 \left(\sqrt{\alpha
   ^2+6 \alpha  \varpi +\varpi  (\varpi +8)}+\varpi
   +8\right)+\alpha ^2 \varpi  \left[6 \left(\sqrt{\alpha ^2+6
   \alpha  \varpi +\varpi  (\varpi +8)}+2\right)+19 \varpi
   \right]\right.\nonumber\\
   &&\left.+\alpha ^3 \left(\sqrt{\alpha ^2+6 \alpha  \varpi
   +\varpi  (\varpi +8)}+9 \varpi \right)\right\}\nonumber\\
D_T &=& 4 \sqrt{\alpha ^2+6 \alpha  \varpi +\varpi  (\varpi +8)}
   \left(\sqrt{\alpha ^2+6 \alpha  \varpi +\varpi  (\varpi
   +8)}+\alpha +\varpi \right)^2
\end{eqnarray}
The new function $G_T(\varpi )$ maps the interval [0,$\infty $] of
$\varpi $ into the interval $\left[-\frac{3 \alpha
}{8},-\frac{3}{8} (4+3 \alpha )\right]$ so that the range of the
negative variable $t$ is
\begin{equation}\label{divisor22}
    t\in \left[-\frac{3}{8}(4+3\alpha ),-\frac{3\alpha }{8}\right]
\end{equation}
and, as expected, the cycle $C_1$ does not shrink to zero unless
$\alpha $ = 0. Quite surprisingly the function $G_T(\varpi )$ can
be easily inverted and we find:
\begin{equation}\label{divisor23}
    \varpi \equiv H(t,\alpha )^2=\frac{64 t^2-9 \alpha ^2}{54 \alpha
+48 t+72}
\end{equation}
while for $P(t)$ we get:
\begin{equation}\label{divisor24}
P(t,\alpha )\,=\,2 \sqrt{\frac{27 \alpha ^2+432 \alpha  t+192 t
(t+3)}{(t+3) (9 \alpha +8 t+12) \left(9 \alpha ^2-64 t^2\right)}}
\end{equation}
and we verify that
\begin{equation}\label{divisor25}
   P(t,0 )=\sqrt{\frac{-3}{2 t^2+3t}}
\end{equation}
which is the correct result for the singular case ${\mathbb WP}[112]$.
In terms of the function $F(t)$
parameterizing the metric (\ref{divisor15}) we have:
\begin{equation}\label{barilla}
    F(t,\alpha) \, = \, \frac{4 \left(27 \alpha ^2+432
   \alpha  t+192 t (t+3)\right)}{3 (9 \alpha +8 t+12) \left(64 t^2-9 \alpha
   ^2\right)} \, = \, \frac{1}{2}\,\left(\frac{1}{\frac{3 \alpha }{8}+t}-\frac{1}{\frac{3}{8} (3 \alpha +4)+t}+\frac{1}{t-\frac{3
   \alpha }{8}}\right)
\end{equation}
The above structure of the function $F(t,\alpha)$ is very much
inspiring. As we see, it is just the sum of three simple poles that
are alternatively simple poles of the $dt^2$-coefficient and zeros
of the coefficient of the $(d\tau+(1-\cos(\theta) d\phi)^2$-term.
The range of the variable $t$ turns out to be the interval between
two such poles where the sign of the function $F(t)$ is the correct
one for in order for the metric  (\ref{divisor15}) to have Euclidian
signature. The three poles are:
\begin{equation}\label{polari}
    t_1 \, =\, -\frac{3 \alpha }{8} \quad ; \quad t_2 \, = \, -\frac{3}{8} (3 \alpha
    +4) \quad ; \quad t_3 \, =\, \frac{3 \alpha }{8}
\end{equation}
We also see what is the mechanism of the degeneration producing the
singular ${\mathbb WP}[112]$ case: the two poles $t_1$ and $t_3$
come to coincide and the coincidence point is zero. This produces
the vanishing of the $C_1$-cycle as we explained above.
\par
Substituting the function $F(t,\alpha)$ as given in
eq.~(\ref{barilla}) into the metric   we get a final form of a specific K\"ahler
metric on the second Hirzebruch surface which follows from the
Kronheimer construction. This metric provides the boundary condition
for the Ricci-flat metric on the canonical bundle $ {\rm tot} K_{\mathbb{F}_2}$  which must reduce to it when setting $ds=0$, $d\chi
=0$ and $s = - 3$.
\paragraph{\sc Verification of the topological conditions for the K\"ahler
metric of $\mathbb{F}_2$.} As a matter of check we calculate the
periods of the K\"ahler and Ricci 2-forms also in the real
formalism, obtaining the following expected result which holds true
for $0 < \mid \alpha \mid < 1$:
\begin{equation}\label{divisor27}
\int _{C_1}\mathbb{K} = \frac{9 \alpha }{16}\quad;\quad \int
_{C_2}\mathbb{K} = \frac{9 (2+\alpha )}{16}\quad;\quad\int
_{C_1}\text{$\mathbb{R}$ic} = 0\quad;\quad\int
_{C_2}\text{$\mathbb{R}$ic} = 2
\end{equation}
The above result for the K\"ahler form is immediate once the
function $P(t)=P(t,\alpha)$ is specified. It is instead interesting
to see the subtle way in which the result for the Ricci form is
obtained independently from the value of $\alpha$.
\par Calculating the Ricci tensor of the metric in eq.~(\ref{divisor15}) with the
function $F(t,\alpha)$ of eq.~(\ref{barilla}) we find the symmetric
matrix $\mathcal{R}ic$ which, multiplied  by the transpose of the
complex structure (\ref{lomellina}) with $P(t)$ as in
eq.~(\ref{divisor24}) produces the Ricci form
$\text{$\mathbb{R}$ic}_{{\mathcal{ED}}}$ with the structure displayed in
eq.~(\ref{divisor8}) and the following explicit expressions for the
functions $\mathfrak{A}(t)$ and $\mathfrak{C}(t)$.
\begin{eqnarray}
 \mathfrak{A}(t) &=& \frac{(8 t-3 \alpha ) (3 \alpha +8 t) \left(27 \alpha ^2 (3 \alpha +4)+512 t^3+576 (3
   \alpha +4) t^2+216 \alpha ^2 t\right)}{8 t \left(9 \alpha ^2+144 \alpha  t+64 t
   (t+3)\right)^2} \label{lappus1}\\
  \mathfrak{C}(t)&=& \frac{d}{dt} \,U(t) \quad ; \quad
  U(t) \, = \, \frac{864 (\alpha +1) (\alpha +2) \left(3 \alpha ^2+8 (3 \alpha +4) t\right)}{\left(9
   \alpha ^2+144 \alpha  t+64 t (t+3)\right)^2}-\frac{3 (3 \alpha +4)}{8
   t} \label{lappus2}
\end{eqnarray}
We immediately see that $-|t_{max}|=-\frac{3 \alpha }{8}$ is a zero
of $\mathfrak{A}(t)$ so that $\int_{C1}\mathbb{R}
\text{ic}_{\mathcal{ED}}\, = \, 0$, while setting as we must
$-|t_{min}|=-\frac{3}{8} (3 \alpha +4)$ we obtain:
\begin{equation}\label{ugolino}
        U(-|t_{max})-U(-|t_{min})\, = \, 2 \quad \Rightarrow \quad \int_{C2}\mathbb{R}
\text{ic}_{\mathcal{ED}}\, = \, 0
\end{equation}
\subsection{The exceptional divisor in symplectic coordinates.}
Considering next the description of the 6-dimensional manifold
$\mathcal{M}_6$ in terms of symplectic coordinates
$\left\{\mathfrak{u},\mathfrak{v},\mathfrak{w},\phi,\tau,\chi\right\}$
(see sect.\ref{sympatia}) we easily find that the localization $s =
-3$ of the exceptional divisor  corresponds to
$\mathfrak{w}=\frac{3}{2}$, $d\mathfrak{w}={d\chi} =\, 0$.
Hence defining
\begin{equation}\label{sympcord3}
    \mathcal{D}(\mathfrak{v})\equiv \mathcal{G}\left(\mathfrak{v},\frac{3}{2}\right),
\end{equation}
where $\mathcal{G}\left(\mathfrak{v},\mathfrak{w}\right)$ is the
variable part of the overall symplectic prepotential, we obtain that
the K\"ahler metric on the exceptional divisor has also a
description in terms of a symplectic potential given by
\begin{equation}\label{sympcord4}
    \mathfrak{D}(\mathfrak{u},\mathfrak{v})=G_0(\mathfrak{u},\mathfrak{v})
    +\mathcal{D}(\mathfrak{v})=\mathcal{D}(\mathfrak{v})+
    \left(\mathfrak{v}-\frac{\mathfrak{u}}{2}\right)
\log (2 \mathfrak{v}-\mathfrak{u})+\frac{1}{2} (\mathfrak{u} \log
\mathfrak{u}-\mathfrak{u})-\mathfrak{v} \log \mathfrak{v}
\end{equation}
with moment and angular variables $ \mu
^i=\{\mathfrak{u},\mathfrak{v}\}$, $\Theta_j=\{\phi ,\tau \}$ and
line element as follows:
\begin{equation}\label{sympcord6}
    \text{ds}_{\mathcal{ED}}^2 \, = \,D_{ij}{d\mu}^i{d\mu}^j
    + (D^{-1})^{ij}{d\Theta}_i{d\Theta}_j
\end{equation}
where the two matrices are:
\begin{eqnarray}\label{sympcord7}
&&  D_{ij} \,=\, \left(
\begin{array}{cc}
 -\frac{\mathfrak{v}}{\mathfrak{u}^2-2 \mathfrak{u} \mathfrak{v}} &
 \frac{1}{\mathfrak{u}-2 \mathfrak{v}} \\
 \frac{1}{\mathfrak{u}-2 \mathfrak{v}} &
 {\mathcal{D}''(\mathfrak{v})}-\frac{\mathfrak{u}}{\mathfrak{u}
\mathfrak{v}-2 \mathfrak{v}^2} \\
\end{array}
\right)\quad ; \quad
   (D^{-1})^{ij}\,=\, \left(
\begin{array}{cc}
 \frac{\mathfrak{u} \left(\mathfrak{v} (2 \mathfrak{v}-\mathfrak{u})
 \mathcal{D}''(\mathfrak{v})+\mathfrak{u}\right)}{\mathfrak{v}^2
 \mathcal{D}''(\mathfrak{v})}
& \frac{\mathfrak{u}}{\mathfrak{v} \mathcal{D}''(\mathfrak{v})} \\
 \frac{\mathfrak{u}}{\mathfrak{v} \mathcal{D}''(\mathfrak{v})} &
 \frac{1}{\mathcal{D}''(\mathfrak{v})} \\
\end{array}
\right)
\end{eqnarray}
Reduced to the exceptional divisor, the coordinate transformation
(\ref{buko8}) is very simple. We have:
$\mathfrak{u}=\frac{3}{4}t(-1+\cos\theta)$,
$\mathfrak{v}=-\frac{3t}{4}$. So if we declare that the function
$\mathcal{D}(\mathfrak{v}) = \Pi (t)$, is a function of $t$ we
obtain $\mathcal{D}''(\mathfrak{v})=\frac{16}{9}\Pi'' (t)$ and
replacing  these transformation in (\ref{sympcord3}-\ref{sympcord4})
we obtain that the line element in symplectic coordinates coincides
with the line element of eq.~(\ref{divisor2}) provided that:
\begin{equation}\label{sympcord11}
\mathcal{D}(\mathfrak{v}) \equiv \Pi (t)  \quad ; \quad
P(t)^2=-\frac{4 \Pi ''(t)}{t+3} \quad \Rightarrow \quad
    \Pi ''(t)\, = \, -\frac{3}{4} F(t)
\end{equation}
So the function $F(t)$ determining the K\"ahler geometry of the
exceptional divisor,  linked to its K\"ahler potential by
eq.~(\ref{divisor17}),  is just $4/3$$\times$ the second derivative
of the non-fixed part of the symplectic potential.
\paragraph{\sc The case of the  K\"ahler metric on $\mathbb{F}_2$ with generic $\alpha $.}
Applying the above scheme to the K\"ahler metric on $\mathbb{F}_2$
induced by the Kronheimer construction, namely utilizing in
eq.~(\ref{sympcord11}) $F(t)=F(t,\alpha)$ as given in
eq.~(\ref{barilla}) we obtain the following differential equation:
\begin{equation}\label{sympcord13}
   \mathcal{D}''(\mathfrak{v};\alpha)=16 \left(\frac{1}{27 \alpha -32
\mathfrak{v}+36}+\frac{1}{32 \mathfrak{v}-9 \alpha }+\frac{1}{9
\alpha +32 \mathfrak{v}}\right)
\end{equation}
which, modulo linear functions implies $
\mathcal{D}(\mathfrak{v},\alpha )\,=\,
\mathcal{D}_{\mathbb{F}_2}(\mathfrak{v},\alpha )$ where
\begin{eqnarray}\label{sympcord14}
    \mathcal{D}_{\mathbb{F}_2}(\mathfrak{v},\alpha )
    & \equiv & \frac{1}{2} \mathfrak{v} \log
\left(1024 \mathfrak{v}^2-81 \alpha ^2\right)+\frac{1}{64} (27
\alpha -32 \mathfrak{v}+36) \log (27 \alpha -32
\mathfrak{v}+36) \nonumber \\ && +  \frac{9}{32} \alpha  \,
\text{arctanh}\left(\frac{32 \mathfrak{v}}{9 \alpha
}\right)-\frac{\mathfrak{v}}{2}
\end{eqnarray}
We also find:
\begin{eqnarray}\label{sympcord15}
    \mathcal{D}_{\mathbb{F}_2}(\mathfrak{v},0)& =&
-\frac{1}{16}(-9+8\mathfrak{v})\log\left(3-\frac{8}{3}\mathfrak{v}\right)
+\mathfrak{v} \log\mathfrak{v}\quad \text{modulo a linear
function of}\,\, \mathfrak{v}
\end{eqnarray}
For $\alpha $ = 0 this is the correct result for ${\mathbb WP}[112]$.

\section{The Monge-Amp\`{e}re equation and its series expansion}
In this section we arrive at the core of the issue, {\it i.e.}
the construction of Ricci-flat metrics on the spaces we are
concerned with. The common general feature of these   is that
they are the total space of the canonical bundle of a complex
two-dimensional compact K\"ahler manifold $\mathcal{M}_4$, the \textit{exceptional divisor} when the total space is the full or partial
resolution of a quotient singularity. In this interpretation the
base of the canonical bundle is indeed the exceptional divisor
produced by the \textit{blow up} of an isolated singular point.
\par
The additional common structural feature of the Ricci-flat metrics
we want to consider is, as we already stressed several times, the
group of continuous isometries that they should possess, mentioned
in equation (\ref{gisometry}). The action of $\mathrm{G}_{iso}$ on
the three complex coordinates $u,v,w$ that originate from the
integration of the complex structure was displayed in
eq.s~(\ref{ciabattabuona}).
The presence of these isometries imposes very   stringent
constraints on the K\"ahler metric that are most efficiently handled
at the level of the potential $\mathfrak{P}$ from which the metric
can be obtained by means of derivatives. The condition of
Ricci-flatness of the metric is translated into a nonlinear
differential equation to be satisfied by the potential
$\mathfrak{P}$ that we name the Monge-Amp\`{e}re equation.
\par
As we have seen in the previous pages, there are three equivalent
formulations of the K\"ahler geometry of the toric six-dimensional
manifolds $\mathcal{M}_6$ we are concerned with:
\begin{description}
  \item[A)] The complex setup where the geometry is encoded in the
  K\"ahler potential
  $\mathfrak{P}\, = \,\mathcal{K}\left(u,v,w,\bar{u},\bar{v},\bar{w}\right)$
  \item[B)] The symplectic setup where the geometry is encoded in
  the symplectic potential $\mathfrak{P}\, =
  \,\mathcal{G}\left(\mathfrak{u},\mathfrak{v},\mathfrak{w}\right)$
  \item[C)] The hybrid setup where the geometry is encoded in the
  symplectic potential, but instead of the coordinates $v,w$ we use
  the coordinates $s$,$t$ related to them by the coordinate
  transformation (\ref{buko8}-\ref{buko9}).
\end{description}
Correspondingly there are, to begin with, two formulations of the
Monge-Amp\`{e}re equation, one for the K\"ahler potential, one for
the symplectic potential. In both cases the constraints imposed by
the chosen isometries reduce the effective potential to be a
function of only two real variables so that the Monge-Amp\`{e}re
equation is a non linear partial differential equation in two
variables. At this point the symplectic case still splits into two
versions depending on whether we employ the pure symplectic variables
or the hybrid ones $s, t$.
\par
In all formulations, as we show below, the equation has the property
that we can fix as boundary condition an arbitrarily chosen K\"ahler
metric on the exceptional divisor.
\subsection{The Monge-Amp\`{e}re equation for the K\"ahler
potential}\label{kelleropotto} We begin with the Monge-Amp\`{e}re
equation written in terms of the K\"ahler potential.  It follows
from the chosen isometries that the K\"ahler potential $\mathcal{K}$
must be a function only of the two invariants:
\begin{equation}\label{invariabilmente}
    \mathfrak{f} \, \equiv \, \mid w\mid^2 \quad ; \quad \varpi \, \equiv \,
    \left( 1 +\mid u \mid^2\right)^2 \, \mid v \mid^2 \quad \text{or}
    \quad \mathcal{T} \, \equiv \,4+\varpi - \sqrt{\varpi\,(\varpi +8)}
\end{equation}
so that we can set:
\begin{equation}\label{parruccadaorso}
    \mathcal{K} \, = \, G(\varpi,\mathfrak{f}) \quad \text{or} \quad \mathcal{K} \, = \, G(\mathcal{T},\mathfrak{f})
\end{equation}
The use of the alternative combination $\mathcal{T}$ simplifies the
K\"ahler potential in certain cases.
\par
The Monge-Amp\`{e}re equation in this setup is simply the
statement that the determinant of the K\"ahler metric is constant.
Indeed in the complex coordinate setup the hermitian Ricci tensor is
obtained from the logarithm of the metric determinant in the same
way as the K\"ahler metric is obtained from the K\"ahler potential:
\begin{equation}\label{Riccitensor}
     \text{Ric}_{ij^\star}\, = \, \frac{\partial}{\partial z^i} \,
   \frac{\partial}{\partial \bar{z}^{j^\star}} \log\left[\text{det}
   \mathbf{g}\right] \quad ; \quad \mathbf{g} \, = \, g_{ij^\star}\, = \, \frac{\partial}{\partial z^i} \,
   \frac{\partial}{\partial \bar{z}^{j^\star}}G(\mathcal{T},\mathfrak{f})
   \quad ; \quad z^i \equiv \{u,v,w\}
\end{equation}
Hence if
\begin{equation}\label{barnabino}
\text{det} \, \mathbf{g} \, = \, \kappa
\end{equation}
where $\kappa$ is a constant parameter, the Ricci tensor is
necessarily zero and we have a Ricci-flat metric. The
Monge-Amp\`{e}re equation is obtained by replacing in
eq.~(\ref{barnabino}) the expression of $\text{det} \, \mathbf{g} $
in terms of derivatives of the K\"ahler potential
$G(\mathcal{T},\mathfrak{f})$. Relying on the definition of the
invariants provided in eq.~(\ref{invariabilmente})  we obtain:
\begin{eqnarray}\label{mongolocomplesso}
   4 \mathcal{T}^3 G^{(1,0)}(\mathcal{T},\mathfrak{f} )
   \left\{G^{(0,1)}(\mathcal{T},\mathfrak{f} )
   \left[\left(\mathcal{T}^2+8 \mathcal{T}-16\right) G^{(1,0)}(\mathcal{T},\mathfrak{f} )
   +\mathcal{T}
   \left(\mathcal{T}^2-16\right) G^{(2,0)}(\mathcal{T},\mathfrak{f} )\right]\right.
   \nonumber\\
    \left.+\mathfrak{f}
   \left[G^{(0,2)}(\mathcal{T},\mathfrak{f} ) \left\{\left(\mathcal{T}^2+8
   \mathcal{T}-16\right)
   G^{(1,0)}(\mathcal{T},\mathfrak{f} )+\mathcal{T}
   \left(\mathcal{T}^2-16\right) G^{(2,0)}(\mathcal{T},\mathfrak{f}
   )\right\}\right.\right.\nonumber\\
    \left.\left.-\mathcal{T}
   \left(\mathcal{T}^2-16\right) G^{(1,1)}(\mathcal{T},\mathfrak{f}
   )^2\right]\right\}\, = \,-\kappa  (\mathcal{T}+4)^4
\end{eqnarray}
One can solve the Monge-Amp{\`e}re equation in the above form by
developing the K\"ahler potential in power series of $\mathfrak{f}
$:
\begin{equation}
G(\mathcal{T},\mathfrak{f} )= G_0(\mathcal{T})+\sum _{n=1}^{\infty }
G_n(\mathcal{T}) \,\mathfrak{f} ^n \label{eqqaB4}
\end{equation}
where $G_0(\mathcal{T})$ is the K\"ahler potential of a convenient
K\"ahler metric defined over the exceptional divisor.
\par
Indeed it is  a property of the considered system that inserting
(\ref{eqqaB4}) into the Monge-Amp{\` e}re equation
(\ref{mongolocomplesso}), the function $G_0(\mathcal{T})$
corresponding to the K\"ahler potential of the K\"ahler metric on
the exceptional divisor  is undetermined, while all the other
$G_n(\mathcal{T})$ functions can be iteratively determined in
terms of the previous $G_{k<n}(\mathcal{T})$.
\par
As we discussed before, it is quite remarkable that on the
exceptional divisor located at $s=-3$ the Ricci-flat orthotoric
metric (\ref{iniz11}) reduces  precisely to the K\"ahler metric on
${\mathbb WP}[112]$, which was obtained in \cite{noietmarcovaldo}
from the Kronheimer construction while performing the partial
resolution of the $\mathbb{C}^3/\mathbb{Z}_4$ singularity on a type
III wall.
\subsubsection{Recursive solution for the K\"ahler potential in the case
 ${\rm tot} K_{{\mathbb WP}[112]}$} 
\label{granmogul}
In this
section we present the recursive solution of the Monge-Amp\`{e}re
equation which was obtained  by means of a  dedicated MATHEMATICA
code using as 0-th order K\"ahler potential  the following one
\begin{equation}
G_0(\mathcal{T})= 4\log \mathcal{T}+\mathcal{T}\label{zeroth}
\end{equation}
which yields the Kronheimer K\"ahler metric on ${\mathbb WP}[112]$.
\par
The K\"ahler potential of the full metric on
${\rm tot} K_{{\mathbb WP}[112]}$ can be expressed as follows
\begin{equation}
G(\mathcal{T},\mathfrak{f} )=4\log \mathcal{T} +\mathcal{T}
+\frac{1}{2}\left( \sum _{k=1}^{\infty } \frac{1}{(k!)}
\frac{\mathcal{P}_{2k-2}(\mathcal{T})}{(2\mathcal{T})^{2k}(\mathcal{T}+4)^{2k-3}}(\kappa
 \mathfrak{f} )^k\right)
\end{equation}
where the symbol $\mathcal{P}_{2k-2}(\mathcal{T})$ denotes a
polynomyal of degree $2k-2$ in the variable $\mathcal{T}$. The
remarkable feature is that the coefficients of the polynomials
$\mathcal{P}_{2k-2}(\mathcal{T})$ are all integer numbers whose
decomposition into prime factors involves prime numbers of
increasing values. We show the first 6 of these intriguing
polynomials

\vskip 1.5mm

\noindent$k = 1 ~  | ~ \mathcal{P}_{0}(\mathcal{T})=   2$\\
\\
\noindent$k = 2 ~  | ~ \mathcal{P}_{2}(\mathcal{T})=   112+16 \mathcal{T}+\mathcal{T}^2$\\
\\
\noindent$k = 3 ~  | ~ \mathcal{P}_{4}(\mathcal{T})=   2 \left(10112+4000 \mathcal{T}+408 \mathcal{T}^2+30 \mathcal{T}^3+\mathcal{T}^4\right)$\\
\\
\noindent$k = 4 ~  | ~ \mathcal{P}_{6}(\mathcal{T})=6563840+4347392 \mathcal{T}+925952 \mathcal{T}^2+82624 \mathcal{T}^3+7112 \mathcal{T}^4+350 \mathcal{T}^5+8 \mathcal{T}^6   $\\
\\
\noindent$k = 5 ~  | ~ \mathcal{P}_{8}(\mathcal{T})=  3128950784+2919825408
\mathcal{T}+987267072 \mathcal{T}^2+150301696 \mathcal{T}^3+13354240\mathcal{T}^4+1313920 \mathcal{T}^5\\
{}\qquad\qquad\qquad\; +76064\mathcal{T}^6+2812 \mathcal{T}^7+49 \mathcal{T}^8$\\
\\
\noindent$k = 6 ~  | ~ \mathcal{P}_{10}(\mathcal{T})= 1980772122624+2387983728640
\mathcal{T}+1118459035648 \mathcal{T}^2+256754671616
\mathcal{T}^3\\
{}\qquad\qquad\qquad~~\;\;+32204621824 \mathcal{T}^4+3128804864
\mathcal{T}^5+331169920 \mathcal{T}^6+20666912 \mathcal{T}^7+975904
\mathcal{T}^8+28886 \mathcal{T}^9\\
{}\qquad\qquad\qquad~~\;+407 \mathcal{T}^{10}  $\\
\subsection{The Monge-Amp\`{e}re equation of Ricci-flatness for the
symplectic potential}\label{sympotto} According to \cite{abreu,Martelli:2005tp}
the condition for Ricci-flatness can be written as a differential
condition on the symplectic potential which is the following
\begin{equation}\label{mongolon1}
    \text{Det}\left[G_{ij}\right]=\text{const} \times
\text{Exp}\sum _{h=1}^n c^h\partial _hG
\end{equation}
where $c^h$ are some constants. In the case of our general metric
with isometry $\mathrm{SU(2)\times U(1)\times U(1)}$, the symplectic  form
of the Monge-Amp\`{e}re equation simplifies since we have the
particular form (\ref{buko14}) of the matrix $G_{ij}$.
Indeed we find:
\begin{equation}\label{mongolon2}
    \text{detHes} \equiv
\text{Det}\left[G_{ij}\right]=\frac{\mathfrak{v}}{\mathfrak{u}
(\mathfrak{u}-2 \mathfrak{v})}
\left[\mathcal{G}^{(1,1)}(\mathfrak{v},\mathfrak{w})^2-\mathcal{G}^{(0,2)}(\mathfrak{v}, \mathfrak{w})
\mathcal{G}^{(2,0)}(\mathfrak{v},\mathfrak{w})\right]
\end{equation}
This facilitates the study of the Ricci-flatness condition because
the coefficients $c^{\mathfrak{u}}$ and $c^{\mathfrak{v}}$ are
already fixed by the need to reproduce the $\mathfrak{u}$-dependence
of $\mathrm{detHes}$. We easily find:
\begin{equation}\label{mongolon3}
\text{Exp}\left[-2\partial_{\mathfrak{u}}G(\mathfrak{u},\mathfrak{v},\mathfrak{w})
\, -
2\partial_{\mathfrak{v}}G(\mathfrak{u},\mathfrak{v},\mathfrak{w})\,
+\,
k\,\partial_{\mathfrak{w}}G(\mathfrak{u},\mathfrak{v},\mathfrak{w})\right]
\, = \, -\frac{e^{k \mathcal{G}^{(0,1)}(\mathfrak{v},\mathfrak{w})-2
\mathcal{G}^{(1,0)}(\mathfrak{v},\mathfrak{w})}
\mathfrak{v}^2}{\mathfrak{u}^2-2 \mathfrak{u} \mathfrak{v}}
\end{equation}
Hence in the symplectic formalism the Monge-Amp\`{e}re equation for Ricci
flatness reduces  to the following relation:
\begin{equation}\label{sympmongolone}
    c e^{k \mathcal{G}^{(0,1)}(\mathfrak{v},\mathfrak{w})-2
\mathcal{G}^{(1,0)}(\mathfrak{v},\mathfrak{w})}
\mathfrak{v}+\mathcal{G}^{(1,1)}(\mathfrak{v},\mathfrak{w})^2
-\mathcal{G}^{(0,2)}(\mathfrak{v},\mathfrak{w})
\mathcal{G}^{(2,0)}(\mathfrak{v},\mathfrak{w}) = 0
\end{equation}
imposed solely on the function of two variables
$\mathcal{G}$[$\mathfrak{v}$,$\mathfrak{w}$].
\par
We have explicitly verified that the function
$\mathcal{G}_{{{\mathbb WP}[112]}}(\mathfrak{v},\mathfrak{w})$
defined in equation (\ref{miracolina}), which corresponds to the
orthotoric Ricci-flat metric on
${\rm tot} K_{{\mathbb WP}[112]}$ satisfies
eq.~(\ref{sympmongolone}) with:
\begin{equation}\label{mongolon4}
    k =- \frac{8}{3}\quad; \quad c=\frac{72 e^3}{7}
\end{equation}
\subsubsection{Discussion of the boundary condition}
As we show below, differently from the case  of the Monge-Amp\`{e}re
equation for the K\"ahler potential in the symplectic case, there is
a subtle issue concerning the choice of boundary condition to be
imposed on the function while restricting it to the exceptional
divisor. The important point is that at the level of the metric the
limit $\mathfrak{w}\to \ft 32 $ should reproduce the metric on the
divisor derived from the potential $\mathcal{D}(\mathfrak{v})\, = \,
\mathcal{G}\left(\mathfrak{v},\ft 32\right)$. There are only two
ways to obtain this. If the symplectic potential
$\mathcal{G}(\mathfrak{v},\mathfrak{w})$ is holomorphic at
$\mathfrak{w}=\ft 32$ and admits a Taylor series expansion in
$\mathfrak{w}-\ft 32$ we are obliged to impose that
$\partial_w\mathcal{G}(\mathfrak{v},\mathfrak{w})$ be a
constant at $\mathfrak{w}=\ft 32$ and this results in a recursive
solution with coefficients that are rational functions of increasing
order and can hardly define a convergent series. Furthermore the
only known solution of the Monge Amp\`{e}re equation, provided by
the function (\ref{miracolina}) corresponding to the orthotoric
metric on ${\rm tot} K_{{\mathbb WP}[112]}$ has not this
holomorphic behavior. Indeed
$\mathcal{G}_{{{\mathbb WP}[112]}}(\mathfrak{v},\mathfrak{w})$
provides a paradigmatic example of the other possible boundary
condition which foresees a logarithmic singularity of the symplectic
potential while approaching the exceptional divisor:
\begin{equation}\label{miavvicino}
    \mathcal{G}(\mathfrak{v},\mathfrak{w})\,
    \stackrel{\mathfrak{w}\to \ft 32}{\approx} \,
    \left(\mathfrak{w}-\ft 32\right) \, \log\left(\mathfrak{w}-\ft
    32\right) \, + \,\mathcal{G}_0(\mathfrak{v})
    + \mathcal{O}\left(\mathfrak{w}-\ft 32\right)
\end{equation}
In the sequel we show that with the second type of boundary
condition we can reconstruct the known solution
$\mathcal{G}_{{{\mathbb WP}[112]}}(\mathfrak{v},\mathfrak{w})$ of
equation (\ref{miracolina}) and also derive a series solution
pertaining to the smooth ${\mathbb F}_2$ case which displays the same
general features as
$\mathcal{G}_{{{\mathbb WP}[112]}}(\mathfrak{v},\mathfrak{w})$.
Unfortunately, up to the present moment we can only give numerical
evidences of the last statement.
\par
In view of what we explained above we skip the details concerning
the first type of boundary condition (holomorphicity at
$\mathfrak{w}=\ft 32$ and jump directly to the case of a
logarithmic singularity at $\mathfrak{w}=\ft 32$. Indeed, a logarithmic singularity is known to be the correct behaviour to ensures smoothness of the toric K\"ahler metrics near 
to divisors \cite{guillemin,donaldson,abreu}.

\subsection{The boundary condition with a logarithmic singularity
at $\mathfrak{w}=\ft 32$}
We implement the second type of boundary condition requiring that
following two properties should be preserved:
\begin{description}
  \item[a)] The symplectic potential
  $\mathcal{G}(\mathfrak{v},\mathfrak{w})$ has a finite limit for
  $\mathfrak{w}\to \ft 32$
  \item[b)] The limit for $\mathfrak{w}\to \ft 32$ of the bundle
  metric  should be exactly the exceptional
  divisor metric (\ref{sympcord6}-\ref{sympcord7})
\end{description}
Namely we must have:
\begin{equation}\label{limitoni}
    \lim_{\mathfrak{w}\to \ft 32} \mathcal{G}(\mathfrak{v},\mathfrak{w}) \, = \,
    \mathcal{D}(\mathfrak{v}) \quad ; \quad \lim_{\mathfrak{w}\to \ft 32} \mathrm{ds}^2_{symp} \, = \,
    \mathrm{ds}^2_{{\mathcal{ED}}}
\end{equation}
To discuss this alternative boundary condition it is convenient to
use rescaled variables defined as follows
\begin{equation}\label{xydefinio}
    x=2\mathfrak{v} \quad ; \quad y = 3\mathfrak{w} \quad ; \quad y=\frac{9}{2}+\omega
\Rightarrow  \omega = 3\, \left(\mathfrak{w} - \ft 32\right)
\end{equation}
In terms of such variables the Monge-Amp\`{e}re equation
(\ref{sympmongolone}) becomes
\begin{equation}\label{mongolfiera}
    c\, x\, \exp \left[-8 \,\mathcal{G}^{(0,1)}(x,\omega )-4
    \,\mathcal{G}^{(1,0)}(x,\omega
)\right]-\mathcal{G}^{(1,1)}(x,\omega
)^2+\mathcal{G}^{(0,2)}(x,\omega ) \mathcal{G}^{(2,0)}(x,\omega )\,
= \, 0
\end{equation}
Instead of assuming that $\mathcal{G}(x,\omega )$ is holomorphic at
$\omega=0$, we impose that it has a logarithmic singularity of the
form $\omega \log\omega$. Indeed this is the unique alternative
way in which the metric on the total space can reduce to the metric
exceptional divisor in the limit $\omega \to 0$. Furthermore this
behavior for $\mathfrak{w} \to \ft 32$ is precisely that displayed
by the symplectic potential
$\mathcal{G}_{{{\mathbb WP}[112]}}(\mathfrak{v},\mathfrak{w})$
explicitly written down in eq.~(\ref{miracolina}). Hence we assume
the following different series expansion which isolates a
logarithmic singularity at $\omega =0$:
\begin{equation}\label{ciulifusco}
    \mathcal{G}(x,\omega )\, =\, \ft 18 \,\omega \, \log[\omega]
     + \mathcal{G}_0(x) \, + \,
\sum _{k=1}^{\infty } \omega ^k \mathcal{G}_k(x)
\end{equation}
The function $\mathcal{G}_0(x)$ is free. All the functions
$\mathcal{G}_k(x)$ ($k\geq $1) are determined in terms of
$\mathcal{G}_0(x)$. For instance we have:
\begin{eqnarray}\label{carrettodeldiavolo}
  \mathcal{G}_0(x) &=& \mathcal{G}_0(x) \nonumber\\
  \mathcal{G}_1(x)&=& \frac{1}{8} \log \left(-\frac{e^{-4 \mathcal{G}_0'(x)} x}{\mathcal{G}_0''(x)}\right)\nonumber \\
  \mathcal{G}_2(x) &=& \frac{32 x^2 \mathcal{G}_0''(x){}^3-12 x \mathcal{G}_0''(x){}^2+16 x^2 \mathcal{G}_0{}^{(3)}(x) \mathcal{G}_0''(x)+2 \mathcal{G}_0''(x)-2 x \mathcal{G}_0{}^{(3)}(x)+x^2 \mathcal{G}_0{}^{(4)}(x)}{256
x^2 \mathcal{G}_0''(x){}^2} \nonumber \\
  \mathcal{G}_3(x) &=& \frac{1}{36864 x^4 \mathcal{G}_0''(x){}^4}\left(-48 x^2 \left(16 \mathcal{G}_0{}^{(3)}(x) x^2
  +7\right) \mathcal{G}_0''(x){}^4-48 x \left(12 x^3 \mathcal{G}_0{}^{(4)}(x)-5\right)
\mathcal{G}_0''(x){}^3\right.\nonumber\\
&&\left.+4 \left(144 \mathcal{G}_0{}^{(3)}(x){}^2 x^4-18
\mathcal{G}_0{}^{(5)}(x) x^4+38 \mathcal{G}_0{}^{(4)}(x) x^3+12
\mathcal{G}_0{}^{(3)}(x) x^2-11\right)
\mathcal{G}_0''(x){}^2\right.\nonumber\\
&&\left.+2 x \left(-152 x^2 \mathcal{G}_0{}^{(3)}(x){}^2+\left(72
\mathcal{G}_0{}^{(4)}(x) x^3+4\right) \mathcal{G}_0{}^{(3)}(x)\right.\right.\nonumber\\
&&\left.\left.-x \left(\mathcal{G}_0{}^{(6)}(x) x^2-6
\mathcal{G}_0{}^{(5)}(x) x+2 \mathcal{G}_0{}^{(4)}(x)\right)\right)
\mathcal{G}_0''(x)+9 x^2 \left(x \mathcal{G}_0{}^{(4)}(x)-2
\mathcal{G}_0{}^{(3)}(x)\right){}^2\right)\nonumber\\
\end{eqnarray}
\subsubsection{Recursive solution of the symplectic Monge-Amp\`{e}re equation in the case
where the smooth ${\mathbb F}_2$ surface is at the
boundary}\label{sezionedeltone} Relying on the results of the
previous subsection we consider the case where the symplectic
potential at the boundary (\textit{i.e.}~on the exceptional divisor)
is the one yielding the Kronheimer K\"ahler metric on
${\mathbb F}_2$. In terms of the $x$ variable and setting, $ \Delta
 =  \ft 98  \alpha$ the function in eq.~(\ref{sympcord14})
can be rewritten as follows:
\begin{equation}\label{G0x}
   \mathcal{G}_0(x,\Delta) \equiv \frac{1}{4} x \log \left(4 x^2-\Delta ^2\right)
   +\frac{1}{8} \,\Delta
    \log \left(\frac{\Delta+2 x}{\Delta -2 x}\right)
   + \frac{1}{8}\left(3 \Delta -2
   x+\frac{9}{2}\right) \log \left(3 \Delta -2
   x+\frac{9}{2}\right)
\end{equation}
from the formal solution discussed in the previous section we
obtain:
\begin{equation}\label{solutdelton}
  \mathcal{G}(x,\omega) =  \mathcal{G}_0(x,\Delta) +\frac{1}{8} \omega  \log \left(\frac{x (6 \Delta -4 x+9)^2}{2 e
   \left[\Delta ^2+4 x^2-6 (2 \Delta +3) x\right]}\right)
    + \frac{1}{8} \, \omega \, \log \omega  +\sum_{k=1}^\infty \frac{N_{k+1}(x,\Delta)}{D_{k+1}(x,\Delta)} \,
   \omega^{k+1}
\end{equation}
where $N_{k+1}(x,\Delta)$ and $D_{k+1}(x,\Delta)$ are polynomials
whose degrees are as follows:
\begin{equation}\label{degradi}
    \text{degree} \left[N_{k+1}(x,\Delta)\right]\, = \, 6 \, k \quad
    ; \quad \text{degree} \left[D_{k+1}(x,\Delta)\right]\, = \, 7 \, k
\end{equation}
Hence the degree of the coefficient of $\omega^{k+1}$ in the series
expansion is a rational function of $x$ of degree $-k$, a feature
that looks promising for convergence.
\par
By means of a dedicated MATHEMATICA code we can calculate the
polynomials $N_{k+1}(x,\Delta)$, $D_{k+1}(x,\Delta)$ to any desired
order. For reason of typographical space we display here  only the
first terms up to order $k=2$.  \vskip 0.2cm

\noindent$N_2\text{(x) = }-3888
x^4-864 x^5+128 x^6-5184 x^4 \Delta -576 x^5 \Delta +720 x^3 \Delta
^2-1536 x^4 \Delta ^2+480 x^3 \Delta ^3-81 \Delta ^4+162 x \Delta
^4-120 x^2 \Delta ^4-108 \Delta ^5+108 x \Delta ^5-36 \Delta ^6$
\vskip 0.2cm

\noindent$N_3\text{(x) = }49152 x^{12}-811008 x^{11} (3+2 \Delta
)-891 \Delta ^8 (3+2 \Delta )^4+108 x \Delta ^6 (3+2 \Delta )^3
\left(810+1080 \Delta +421 \Delta ^2\right)-36 x^2 \Delta ^6 (3+2
\Delta )^2 \left(23652+31536 \Delta +10961 \Delta ^2\right)+384 x^6
\Delta ^4 \left(31347+41796 \Delta +12692 \Delta ^2\right)$\\
\noindent$+1024 x^{10} \left(42039+56052 \Delta +18412 \Delta
^2\right)-1536 x^7 \Delta ^2 \left(-9963-19926 \Delta -8238 \Delta
^2+412 \Delta ^3\right)$\\
\noindent$-1536 x^5 \Delta ^4 \left(15795+31590 \Delta +19092
\Delta
^2+3368 \Delta ^3\right)$\\
\noindent$ -3072 x^9 \left(28431+56862 \Delta +41124 \Delta
^2+10568 \Delta ^3\right)+576 x^4 \Delta ^4 \left(35721+95256 \Delta
\right.$ \\ \noindent $\left.+83151 \Delta
^2+26196 \Delta ^3+1655 \Delta ^4\right)\\
+256 x^8 \left(-137781-367416 \Delta -270540 \Delta ^2-34128 \Delta
^3+23272 \Delta ^4\right)$\\
\noindent$+864 x^3 \Delta ^4 \left(-6561-21870 \Delta -17739 \Delta
^2+3402 \Delta ^3+8805 \Delta ^4+2558 \Delta ^5\right)$ \vskip 0.2cm
\vskip 0.2cm \noindent$D_2\text{(x) = }128 x^2 (-9+4 x-6 \Delta )
\left(-18 x+4 x^2-12 x \Delta +\Delta ^2\right)^2$ \vskip 0.2cm
\noindent$D_3\text{(x) = }9216 x^4 (9-4 x+6 \Delta )^2 \left(4
x^2+\Delta ^2-6 x (3+2 \Delta )\right)^4$
\vskip 0.2cm

\subsubsection{Numerical study in the case $\Delta = \ft 34$}
Since  so far we have not been able  to guess the sum of the series
in terms of elementary or higher transcendental functions,
to get some understanding of the solution we have resorted to a
numerical study of the approximants to the solution obtained by
truncating the series in eq.~(\ref{solutdelton}) to various orders
performing the plots.
\par
The relevant thing is that for the special value $\Delta = 0$ of the
parameter we know the exact sum of the series. It is provided by the
symplectic potential (\ref{miracolina}) which pertains to the case
of ${\rm tot} K_{{\mathbb WP}[112]}$.
This fortunate occurrence enables us to compare the plot of the exact
function with those of its approximants. This comparison, as we are
going to see, turns out to be quite inspiring since it elucidates
the meaning of certain oscillatory behaviors of the approximants
that are completely analogous in the case $\Delta =0$, where we know
the sum of the series and in the case $\Delta > 0$ where the sum is
unknown.
\par
In terms of the variables $x$ and $\omega$ the symplectic potential
of the orthotoric metric takes the following explicit expression:
\begin{eqnarray}\label{supercalifragilisti}
  \mathcal{G}_{\Delta=0}(x,\omega) &=&
  \frac{1}{224} \left\{7 \left[(4 \omega+16 x)  \log
   \left({9\over 2} -\sqrt{\left(x+\omega +\frac{9}{2}\right)^2-18 x}- x-
   \omega\right)\right.\right.\nonumber\\
   &&\left.\left.-2 (4 x-4 \omega -9) \log
   \left(\sqrt{\left(x+\omega +\frac{9}{2}\right)^2-18 x}-x-
   \omega \right)\right.\right.\nonumber\\
   &&\left.\left.+(4 x+2 \omega +9) \log
   \left(\frac{1}{567} \left[4 \sqrt{\left(x+\omega
   +\frac{9}{2}\right)^2-18 x}+4 x+4 \omega
   +27\right]^2+1\right)\right]\right.\nonumber\\
   &&\left. -2 \sqrt{7} (4 x-2 \omega -27)
  \arctan\frac{4 \sqrt{\left(x+\omega
   +\frac{9}{2}\right)^2-18 x}+4 x+4 \omega +27}{9
   \sqrt{7}}\right.\nonumber\\
   &&\left. -(4 x+9) \log\frac{34359738368}{823543}
   +2 \sqrt{7} (4 x-27)
   \arctan\frac{5}{\sqrt{7}}\right\}
\end{eqnarray}
For comparison we choose the series solution in the case $\Delta =
\ft 34$. This value, corresponding to $\alpha = \ft 23$, introduces
various simplifications in the solution and, for no other good
reason, provides a good reference point.
\par
\begin{figure}
\centering
\includegraphics[width=8cm]{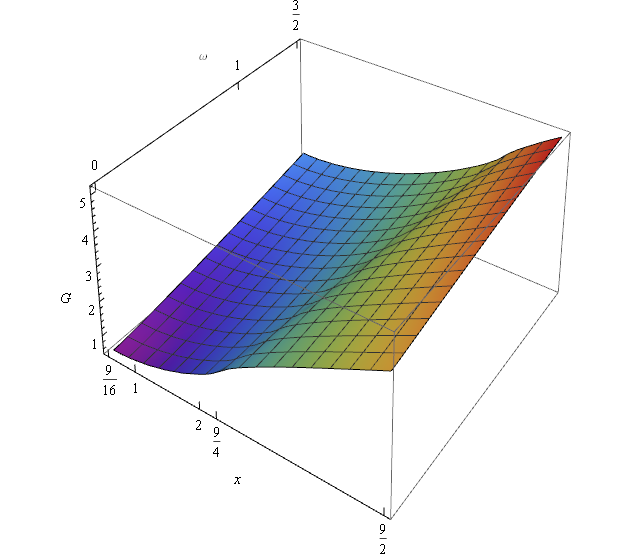}
\includegraphics[width=9cm]{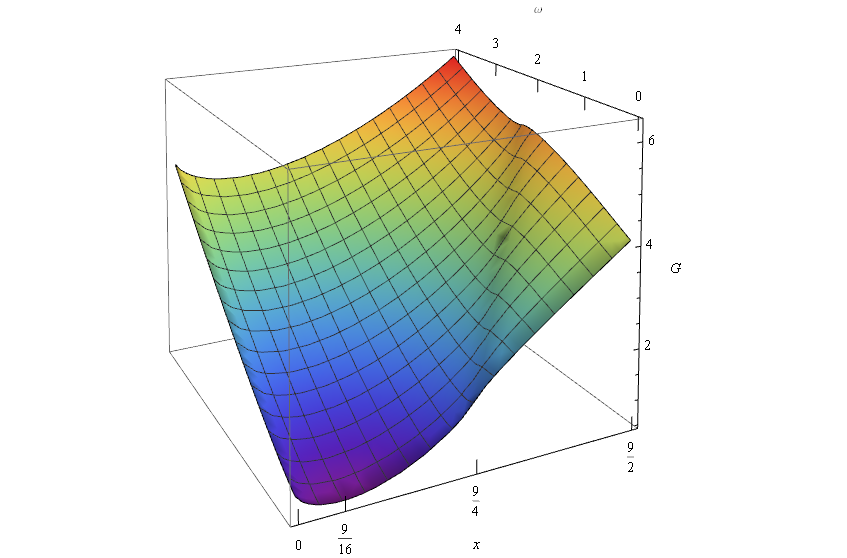}
\caption{\label{figurella1} Plots of the exact symplectic potential
$\mathcal{G}_{\Delta=0}(x,\omega)$ for small values of the distance
$\omega$ from the exceptional divisor (plot on the left) and
extending to large values (plot on the right.)}
\end{figure}

In this case the symplectic potential takes the following appearance
\begin{eqnarray}
  G_{\Delta = \ft 34}(x,\omega) &=& \frac{1}{32} \left[8 x \log
   \left(4 x^2-\frac{9}{16}\right)+(27-8 x)
   \log \left(\frac{27}{4}-2 x\right)+6 \, \text{arctanh}\frac{8
   x}{3}\right]\nonumber\\
   &&\frac{1}{8} \omega  \log
   (\omega ) +\frac{1}{8} \omega  \log
   \left(\frac{2 (27-8 x)^2 x}{e [16 x (4
   x-27)+9]}\right)+ \sum _{k=1}^{\infty } \frac{
   \hat{N}_{k+1}(x)}{\hat{D}_{k+1}(x)} \,\omega ^{k+1}
\end{eqnarray}
We omit the explicit presentation of the rational functions
$\frac{\hat{N}_{k+1}(x)}{\hat{D}_{k+1}(x)}$ that we have calculated
by means of a computer programme up to order $k=10$ and higher. We
rather present the plots of such approximants. Let us first consider
the plot of the function $\mathcal{G}_{\Delta=0}(x,\omega)$
displayed in fig.\ref{figurella1}.

As we distinctly see from the picture, the exact function, namely
the sum of the infinite series in $\omega$ defines parametrically a
perfectly smooth surface in three dimensions that however features
a nontrivial structure provided by a sort of \textit{smooth
bending} along a line that starts approximately at $x=\ft 94,\omega
= 0$ and goes up towards $x=\ft 92, \omega =\infty$.
\par
The geometrically meaning of this bending is not entirely clear,
yet one can guess that it corresponds to a transition region from a
\textit{near divisor geometry} to an asymptotic geometry that is
that of a metric cone over  the Sasakian orbifold $\mathbb{S}^5/\mathbb{Z}_4$. 
\par
It is now interesting to compare the behavior of the exact function
with its approximants obtained truncating the series to various
orders.
\begin{figure}
\centering
\includegraphics[width=8cm]{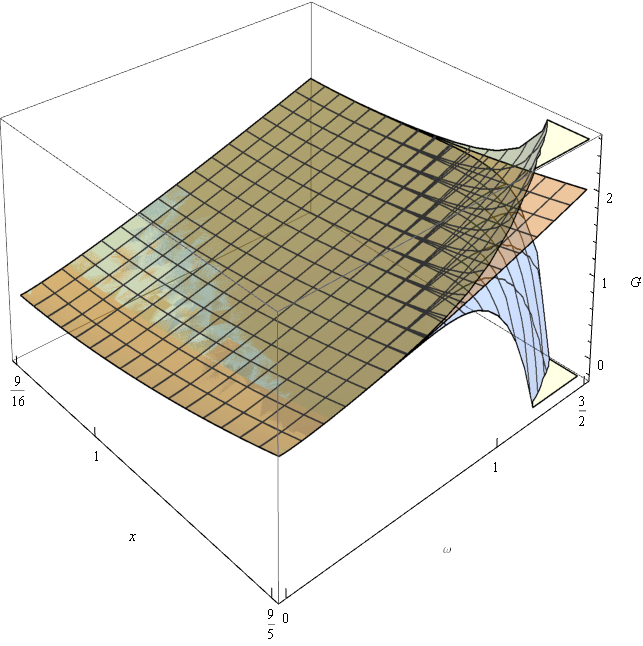}
\includegraphics[width=9cm]{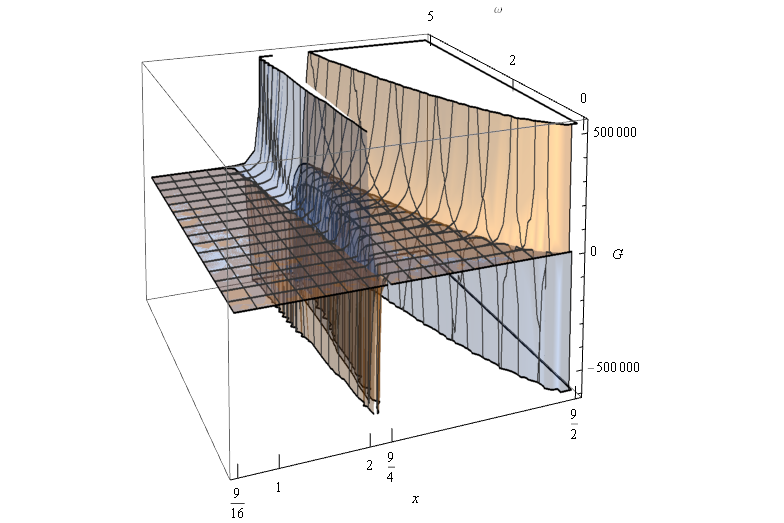}
\caption{\label{figurella2} Plots of the exact symplectic potential
$\mathcal{G}_{\Delta=0}(x,\omega)$ compared to its approximants of
order 6 and 7 respectively: on the right for small values of
$\omega$, on the left extending to large values of $\omega$.}
\end{figure}
Let us now consider the plots displayed in fig.\ref{figurella2}.
\par
In the plot on the right, the surface plotted in the middle is the
sum of the series (\textit{i.e.} the exact function), while the
other two surfaces,  respectively bending, one up, the other down,
are two consecutive approximants (the first of even order, the
second of odd order). As we clearly see, the series converges to the
exact function and does it rapidly, in the region before the bending
structure  illustrated above. As we come close to such a line of
bending the series no longer converges and its various truncations
oscillate violently creating a peculiar canyon.
\par
\begin{figure}
\centering
\includegraphics[width=8cm]{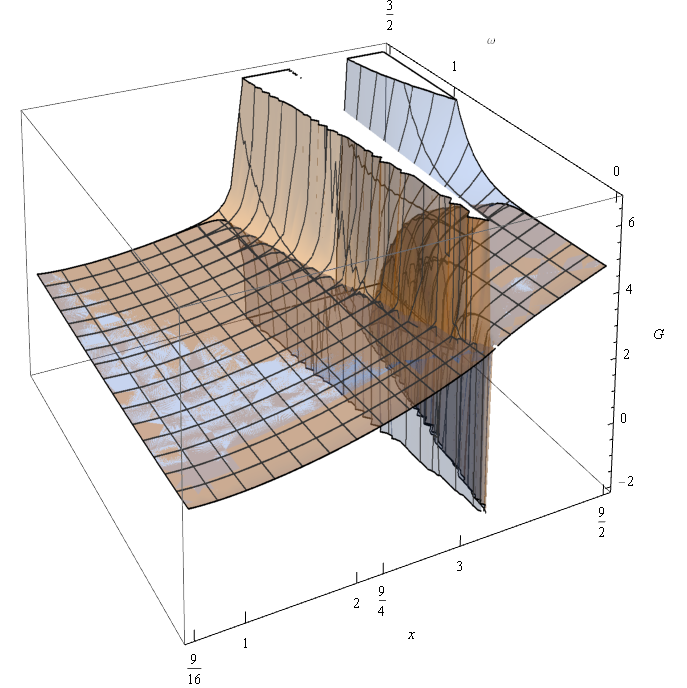}
\includegraphics[width=9.5cm]{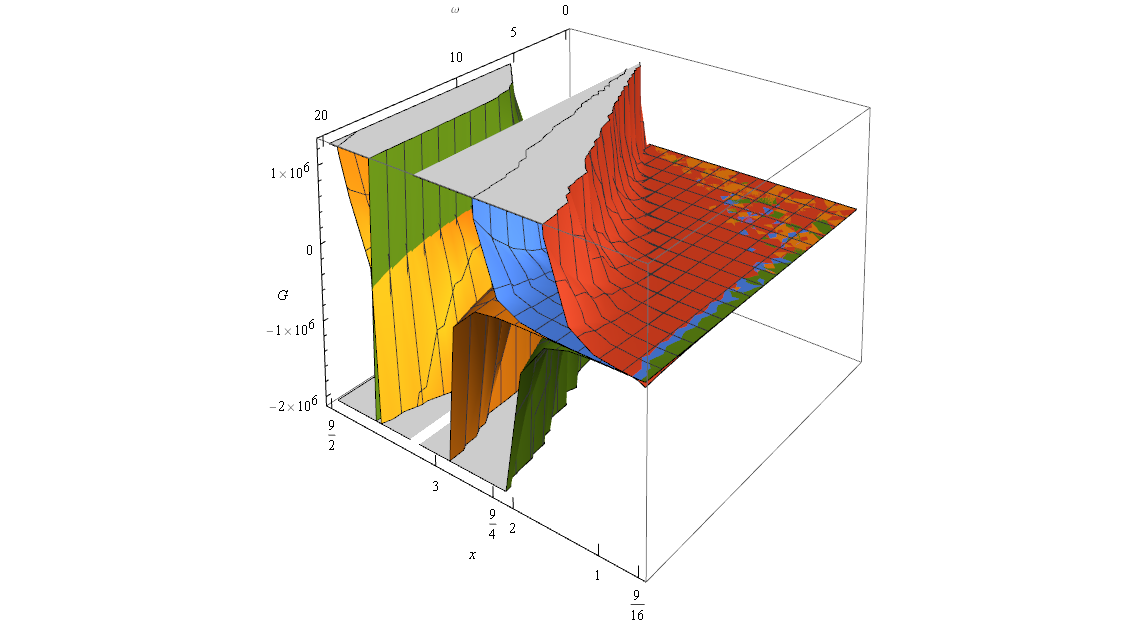}
\caption{Plots of  approximants of even and odd
order of the function $\mathcal{G}_{\Delta=\ft 34}(x,\omega)$ the
sum of whose series representation is unknown. As in the other cases
the plot on the right is for small values of $\omega$ and displays
two consecutive approximants of order 7 and 8, respectively, while
the plot on the left extends to large values of $\omega$ and
displays several approximants. \label{figurella3} }
\end{figure}

Let us now compare this behavior of the case $\Delta = 0$ with that
of the series solution for $\Delta =\ft 34$. To this effect let us
consider the figure \ref{figurella3}.
The structure of the plots of the truncated series are
qualitatively the same in the case $\Delta = \ft 34$, as they are in
the case $\Delta = 0$. Furthermore, in a completely analogous way to
the case $\Delta =0$, for small values of $\omega$ and $x$ also the
series representation of $G_{\Delta = \ft 34}(x,\omega)$ converges
rapidly to some well defined function while approaching the region
of the bending it starts oscillating.  Hence we are led to conclude
that we should be able to retrieve an analytically defined solution
of the Monge-Amp\`{e}re equation for the symplectic potential which
reduces to the Kronheimer metric on ${\mathbb F}_2$ at $\omega=0$. It
is a matter of finding some alternative way of summing the series by
a smart change of variables or by means of some smart integral
transform.
\subsection{The Hybrid version of the Monge-Amp\`{e}re equation}
The most promising setup to study the MA equation for the symplectic
potential is the hybrid one. Working in the $s,t$ coordinates defined
in eq.s~(\ref{buko8}-\ref{buko9})  and setting
$\mathcal{G}(v,w)\,=\,\Gamma(t,s)$, the equation
(\ref{sympmongolone}) is transformed into the following one:
\begin{equation}\label{carduccio}
    \frac{1}{4} c e^{\mathcal{C}} s t = \frac{64
    \mathcal{B}}{(s-t)^2}-\frac{64 \mathcal{A}}{(s-t)^4}
\end{equation}
where:
\begin{eqnarray}
\label{gamberoni}
  \mathcal{C} &\equiv& \frac{8 (s+1) \Gamma ^{(0,1)}(t,s)-8 (t+1) \Gamma ^{(1,0)}(t,s)}{s-t}\nonumber \\
  \mathcal{B} &\equiv& \Gamma ^{(0,2)}(t,s) \Gamma ^{(2,0)}(t,s) \nonumber\\
  \mathcal{A} &\equiv& \left(\Gamma ^{(0,1)}(t,s)-\Gamma ^{(1,0)}(t,s)+(s-t) \Gamma ^{(1,1)}(t,s)\right)^2
\end{eqnarray}
It is an important observation that the term $\mathcal{A}$ is the
square of the constraint whose vanishing implies the orthotoric
separation of the $s,t$ variables (see the last of
eq.s~(\ref{santecchio}). It is interesting to see how with this
separation of variables, namely when $\mathcal{A}=0$, the
differential equation (\ref{carduccio}) does indeed split in two
equations, one for the $t$ variable, the other for the $s$ variable.
On the other hand the equation for the $t$ variable implies the
$\mathbb{W}P[112]$ symplectic potential.
\par
The argument goes as follows. Generalizing the structure of the
known solution for the case $\mathbb{W}P[112]$ we introduce the following
ansatz:
\begin{equation}\label{suzzone}
  \Gamma(s,t)\, =\,  -\frac{1}{3} (2 s+3) {\Pi}(t)+\mathcal{P}(s)+(s+3) \left(\mathcal{Q}(t)-\frac{1}{16} (t+3) \log
   (s+3)\right)+s t \text{Y}_1(s)+(s+t) \text{Y}_2(s)
\end{equation}
where $\Pi(t)$ is an unknown function of $t$ that we would like to
identify with the symplectic potential of the exceptional divisor
metric and $Y_{1,2}(s)$ are also two unknown functions of $s$. On
the other hand the other two functions entering the ansatz are
integral differential functionals of $Y_{1,2}(s)$ and $\Pi(t)$,
respectively :
\begin{eqnarray}
  \mathcal{P}(s) &=& \frac{1}{16} \int \left(-16 \kappa _1-16 s^2 \text{Y}_1'(s)-32 s \text{Y}_2'(s)+s+3\right)
   \, ds \label{sintegralo} \\
 \mathcal{Q}(t) &=& \frac{1}{3} (t+3) \int \frac{(2 t+3) \Pi'(t)-2 \Pi(t)}{(t+3)^2} \,
 dt+\kappa_1 \label{tintegralo}
\end{eqnarray}
With these choices the term $\mathcal{C}$ in eq.~(\ref{gamberoni})
splits into separate functions of different variables:
\begin{eqnarray}
\label{seppioline}
  \mathcal{C} &=& \mathcal{U}(s)+T(t)\nonumber \\
  T(t) &=& \frac{8 \left(3 (t+1) \Pi'(t)+2 (t+3) \int \frac{(2 t+3) \Pi'(t)-2
   \Pi(t)}{(t+3)^2} \, dt-4 \Pi(t)\right)}{3 (t+3)} \nonumber\\
  \mathcal{U}(s) &=& \frac{1}{2} \left(-16 s^2 \text{Y}_1'(s)-16 s \text{Y}_1'(s)-16 \text{Y}_1(s)-16 s
   \text{Y}_2'(s)-16 \text{Y}_2'(s)+16 \text{Y}_2(s)+s-2 \log
   (s+3)+1\right)\nonumber\\
\end{eqnarray}
On the other hand we find that $\mathcal{A}\, = \,0$ while the
$\mathcal{B}$-term factorises as follows:
\begin{eqnarray}
\label{granchio}
  \mathcal{B} &=& \mathcal{H}(t) J(s) \\
  \mathcal{H}(t) &=& \frac{4 \Pi''(t)}{t+3} \\
  J(s) &=& \frac{16 s (s+3) \text{Y}_1''(s)+32 (s+3) \text{Y}_1'(s)+16 s \text{Y}_2''(s)+48
   \text{Y}_2''(s)-1}{s+3}
\end{eqnarray}
In this way the solution of the MA equation reduces to the solution of
two separate integral differential equations one in the $s$
variable, one in the $t$-variable:
\begin{equation}\label{duequazie}
    \frac{\lambda}{2}\, t \, \exp[T(t)] \, = \,\mathcal{H}(t)  \quad
    ; \quad \frac{\mu}{2}\, s \, \exp[\mathcal{U}(s)] \, = \,J(s)
\end{equation}
We focus on the first in the variable $t$. With rather simple
manipulations it can be reduced to an ordinary differential equation
of higher order, namely:
\begin{equation}\label{gongolame}
    \frac{8 (t+1) \Pi''(t)}{t+3}-\frac{\Pi^{(3)}(t)}{\Pi''(t)}+\frac{2
   t+3}{t^2+3 t} \, = \, 0
\end{equation}
which is a differential equation of the first order for the function
$F(t)$. Apart from an integration constant which is fixed by the
topological constraints on the periods of Ricci form, the unique
solution of eq.~(\ref{gongolame}) is $F(t,0)$ corresponding to the
geometry of ${\mathbb WP}[112]$. This shows that in order to impose 
a boundary function consistent 
with $\alpha\neq 0$ we need to
modify the ansatz (\ref{suzzone}) in such a way as to introduce a
certain $s,t$-mixing.
\section{Conclusions}
As we advocated in the introduction, the present paper is an
illustration of the conjecture \ref{cognisco} for which we have
strong support from the fact that it is verified for the value
$\Delta=0$ of the parameter in the paradigmatic case of the
$\mathbb{ C}^3/\mathbb{Z}_4$ singularity resolution. Further
numerical evidence emerges from the study of the power series
solution of the Monge-Amp\`{e}re equation in the symplectic
potential formulation. This latter in its hybrid version seems to
provide the most promising approach since  different series
expansions might be glued together to prolong the solution beyond
the valleys of oscillations.
\par
Assuming that in due time our conjecture can be transformed into a
proof,  we would like to stress its relevance. According to our view
point, Conjecture \ref{cognisco} provides a precise mathematical relationship to realise the
gauge/gravity correspondence in a proper way. The generalized
Kronheimer construction fixes all the items of the gauge theory on
the brane world-volume: field content, gauge group, flavor symmetries and
interactions. As maintained by Conjecture \ref{cognisco}, the same Kronheimer
construction determines, via the Monge-Amp\`{e}re equation, also the
Ricci-flat K\"ahler metric to be used in the construction of the
dual D3-brane solution of supergravity. If \ref{cognisco} is
proved we can say that, 
for the class of theories realised on D3 branes at $\mathbb{C}^3/\Gamma$ Calabi-Yau singularities,
 the McKay quiver determines uniquely both sides of the correspondence.
\section*{Acknowledgments} 
We thank M. Graffeo and V. Peragine for help with the computations
involved in   Appendix A. M.~B. would like to acknowledge illuminating discussions on quivers and dimers with D.~Bufalini, S.~Mancani, A.~Antinucci and F.~Riccioni.

\appendix
\section{The affine variety  $\mathbb C^3/\Z_4$}
We study the quotient $\mathbb C^3/\Z_4$ as an affine variety, i.e.,
as a closed (in the Zariski topology) subset of an affine space
$\C^D$ cut by an ideal $I$ of the polynomial ring $
\C[x_1,\dots,x_D]$. In particular, we show that $\mathbb C^3/\Z_4$
is not a (schematic) complete intersection.
\paragraph{\sc Reminder on complete intersections.} We recall the notions of {\em set-theoretic} and {\em schematic} intersection.
%
\begin{defin} Let $X \subset \C^D$ be an affine variety, and denote by $d$ its codimension in $\C^D$.
\begin{enumerate} \itemsep=-2pt \item
$X$ is a set-theoretic complete intersection if it is cut by $d$
equations as a subset of $\C^D$.
\item $X$ is a schematic complete intersection if it is cut by $d$ equations as an affine variety.
In other terms, if $A(X)$ is the coordinate ring of $X$ (the ring of
regular functions on $X$), then $A(X) = \C[x_1,\dots,x_D]/I$, where
the minimal number of generators of the ideal $I$ is $d$.
\end{enumerate}
\end{defin}
It turns out that all  quotients $\C^3/\Gamma$, where $\Gamma$ is a
finite  abelian subgroup of $\operatorname{SL}_3(\C)$, are
set-theoretic complete intersections, and therefore  so is the case
for $\mathbb C^3/\Z_4$. However, we shall not prove this fact here,
and rather concentrate on proving that $\mathbb C^3/\Z_4$ is not a
schematic complete intersection.
\paragraph{\sc Finding the equations.}
An affine toric variety $X$ can be expressed as
\begin{equation}\label{semig} X_\sigma = \operatorname{Specm}\, \C[S_\sigma]  \end{equation}
where $\sigma\subset N\otimes \R$ is a strongly convex polyhedral
cone, $N$ is  a lattice, and $S_\sigma$ is the semigroup $S_\sigma =
\sigma^\vee\cap M$, with $ \sigma^\vee$ the dual cone to $\sigma$;
$M$ is the dual of the lattice $N$. Specm denotes the maximal
spectrum, i.e., the set of maximal ideals of $\C[S_\sigma]$ with the
Zariski topology. Basically following \cite{CoxLS}, we delineate a
procedure to find the equations for the affine toric variety $X$. We
remind that a {\em Hilbert basis} $\mathcal H_\sigma$ for the
semigroup $S_\sigma$ is a minimal set of generators for $S_\sigma$
which contains the rational generators of the rays of $\sigma^\vee$.
 Define $D =  \# \mathcal H_\sigma$. Then the elements of  $\mathcal H_\sigma$ are related by $D-n$ relations,
which generate an ideal $I_{\sigma,0}$ of $\C[x_1,\dots,x_D]$. Given
two ideals $I$, $J$ in a ring $R$, the saturation of $I$ with
respect to $J$ is defined as
$$I:J^\infty = \{a\in R \, \vert \, a^NJ\subset I \ \mbox{for} \ N\gg 0.\}$$
 Then one proves that the  ideal $I_\sigma$ of $X_\sigma$ in $\C^D$ is  the saturation of $I_{\sigma,0}$ wih the respect to the
ideal $$K_\sigma=(x_1\cdots x_D) \subset \C[x_1,\dots,x_D].$$
\begin{remark} If $I$, $J$ are ideals in $\C[x_1,\dots,x_D]$,
the affine variety corresponding to the ideal $I:J^\infty$  is
$$ V(I:J^\infty) = \overline{V(I)\setminus V(J)}$$
where the closure is taken in $\C^D$ (in the Zariski topology), and
for every ideal $L$ in   $\C[x_1,\dots,x_D]$, $V(L)$ denotes the
closed set in $\C^D$ corresponding to $L$. \label{C2Z2ALE}
\end{remark}
\paragraph{\sc Equations for $X=\C^3/\Z_4$.} Now we check that $\C^3/\Z_4$ is {\em not} a schematic complete intersection, as noted in
\cite{Hanany-Seong}. Realizing $X$ as in equation \eqref{semig} we
can take for $\sigma$ the cone with generators $(1,0,0)$,
$(-1,2,0)$, $(0,-1,2)$ in the latttice $N=\Z^3$. The dual cone
$\sigma^\vee$ has rational generators $(4,2,1),$ $ (0,2,1),$ $
(0,0,1)$ in $M\simeq \Z^3$. A Hilbert basis of $S_\sigma$  is
obtained by adding the lattice points
$$ (1,1,1), \ (0,1,1), \ (1,2,1), \ (2,1,1), \ (2,2,1), \ (3,2,1) .$$
 Assigning variables $x_1,\dots,x_9$ to these lattice points we
obtain that $I_{\sigma,0}$ is generated by the 6 equations
\begin{gather*}
x_1x_8-x_9^2=0,\quad x_2x_9^2-x_8^3=0, \quad x_3x_9^2-x_7^2x_8 = 0, \\
x_4x_9-x_7x_8=0, \quad x_5x_9^2-x_7x_8^2=0, \quad x_6x_9-x_8^2=0
\end{gather*}
Saturating this ideal with respect to $K=(x_1\cdots x_9)$ one sees
that $I_\sigma$ is generated by the 20 quadratic equations (the
equation needed to cut $X$ from $\C^9$ with the correct schematic
structure):
\begin{equation}\label{eqs}
\begin{array}{ccccc}
x_8^2-x_6x_9; & x_7x_8-x_4x_9; & x_6x_8-x_2x_9; & x_4x_8-x_4x_5;  & x_1x_8-x_9^2; \\
x_6x_7-x_5x_9; & x_5x_7-x_3x_8; & x_4x_7-x_3x_9; & x_2x_7-x_5x_8 ; & x_6^2-x_2x_8 ; \\
x_4x_6-x_5x_8; & x_1x_6-x_8x_9; & x_4x_5-x_3x_6; & x_1x_5-x_4x_9; & x_4^2-x_3x_8; \\
x_2x_4-x_5x_6; & x_1-x_7x_9; & x_2x_3-x_5^2; & x_1x_3-x_7^2 ; &
x_1x_2 - x_6x_9.
\end{array}
\end{equation}
 These are a {\em minimal set} of generators. So $X$ is the
intersection of 20 quadrics in $\C^9$. All these quadrics are
singular along their intersection with a plane of codimension 3
(when their equation contains a square) or 4 (when their equation
does not contain a square). The dimension of the singular locus is 6
and 5 respectively (not 5 and 4!)
It may be interesting to see what variety does the ideal
$I_{\sigma,0}$ describe. To this end one computes the primary
decomposition of the ideal \cite{eisen}. This yields 5   ideals; one
is radical, and  coincides with $I_\sigma$,  so that one component
of the variety is $\C^3/\Z_4$. The other ideals are generated by
monomials, and correspond to (intersections of) coordinate planes of
different dimensions, counted with multiplicities.


\section{The orbifold $\mathbb S^5/\mathbb{Z}_4$}
\label{Dariosappendix}

 Setting $s=-\tfrac{2}{3}R^2$  with $r\to \infty$ in the metric  (\ref{iniz11}), it is straightforward to verify that this  takes the approximate form 
 \begin{equation}
 \text{ds}_{{\rm tot}  K_{{\mathbb WP}[112]}}^2   \stackrel{R\to \infty}{\approx}  dR^2 + R^2 \text{ds}^2_{X_5}
\end{equation}
at leading order in $R$.  Since the metric is  Ricci-flat K\"ahler, and it takes the form of a cone over a five-dimensional space, it follows that locally the five-dimensional metric $\text{ds}^2_{X_5}$
is Sasaki-Einstein. 
Below we shall show that globally, this is precisely a Sasaki-Einstein metric on the orbifold $\mathbb{S}^5/\mathbb{Z}_4$.

In the coordinates used in the paper, the five-dimensional metric reads
\begin{eqnarray}\label{first5dmetric}
 \text{ds}^2_{X_5}&=& 
 - \frac{t}{6} \left(\sin^2\theta{d\phi}^2+d\theta^2\right)
-\frac{dt^2}{2t (2 t+3)}\left.-\frac{2 t (2 t+3)}{9} \left[ \frac{d\chi}{3}-\frac{1}{2}  [(1-\cos\theta){d\phi} +{d\tau}]\right]^2\right.\nonumber\\
&& +\frac{4}{9}
\left[\left(\frac{t}{3}+1\right) {d\chi}-\frac{1}{2}t [(1-\cos{\theta}){d\phi} +{d\tau}]\right]^2
\end{eqnarray}
After introducing the  new coordinate $\sigma\in [0,\frac{\pi}{2}]$ as 
 \begin{equation}
 \label{ttosigma}
 t = - \frac{3}{2}\sin^2\sigma
\end{equation}
it becomes
\begin{eqnarray}\label{second5dmetric}
 \text{ds}^2_{X_5}&=&  d\sigma^2 + \frac{\sin^2\sigma}{4} \left(\sin^2\theta{d\phi}^2+d\theta^2\right)
+\frac{\sin^2\sigma\cos^2\sigma}{4} \left[ \frac{2}{3}d\chi - d\tau-d\phi + \cos\theta{d\phi}\right]^2\nonumber\\
&& +\frac{1}{9}\left[2 d\chi-\frac{3}{2}\sin^2\sigma\left( \frac{2}{3}d\chi - d\tau-d\phi + \cos\theta{d\phi}\right)\right]^2
\end{eqnarray}
and one can check that this is indeed locally a Sasaki-Einstein metric, where the first line is a K\"ahler-Einstein metric. In order to uncover the relation 
with the metric on the five-sphere $\mathbb{S}^5$, it is convenient to redefine the angular coordinates as 
\begin{equation}
\label{mycoc}
\tilde \phi  =  \phi\, ,   \qquad  \beta  =  \tfrac{2}{3}\chi - \phi - \tau \, , \qquad \psi  =  2\chi \, ,
\end{equation}
with inverse 
\begin{equation}
\label{myinversecoc}
\phi  =  \tilde \phi\, ,   \qquad  \tau  =  \tfrac{1}{3}\psi- \phi - \beta \, , \qquad \chi  =  \tfrac{1}{2}\psi \, ,
\end{equation}
where, after performing the change of coordinates, we can drop the tilde on $\tilde \phi$ and simply continue to denote this as $\phi$. The metric then reads
\begin{eqnarray}\label{third5dmetric}
 \text{ds}^2_{X_5}&=&  d\sigma^2 + \frac{\sin^2\sigma}{4} \left(\sin^2\theta{d\phi}^2+d\theta^2\right)
+\frac{\sin^2\sigma\cos^2\sigma}{4} \left( d\beta + \cos\theta d\phi\right)^2\nonumber\\
&& +\frac{1}{9}\left[d\psi-\frac{3}{2}\sin^2\sigma\left( d\beta+ \cos\theta{d\phi}\right)\right]^2
\end{eqnarray}

It is well-known (and simple to verify) that taking $\phi \sim \phi + 2\pi$ and $\beta \sim \beta +4\pi$, with $\theta \in [0,\pi]$ and $\sigma\in [0,\tfrac{\pi}{2}]$, the first line is  the standard Einstein metric on 
$\mathbb{P}^2$. Moreover, with $\psi \sim  \psi + 6\pi$, the five-dimensional metric is the round metric on $\mathbb{S}^5$, viewed as the total space of a circle bundle 
$\mathbb{S}^5 \overset{\pi }{\longrightarrow} \mathbb{P}^2$, normalised 
so to obey the equation
\begin{align}
R_{ij}^{X_5} & =  4 g_{ij}^{X_5}
 \end{align} 

On the other hand, we are not free to chose the ranges of the coordinates, but these are inherited from the ranges of the original coordinates $(\phi,\tau,\chi)$, fixed in (\ref{iniz1}). 
From the change of coordinates (\ref{mycoc}), it follows\footnote{The periodicities of $\phi$ and $\psi$ are obvious. The simplest way to determine the periodicity of $\beta$ is by demanding that the total volume of the three-torus with coordinates $(\phi,\chi,\tau)$ is preserved by the coordinate transformation (\ref{mycoc}).} that we must enforce the following periodicities:
 \begin{equation}\label{myperiods}
\phi \sim  \phi +2\pi\, , \qquad 
\beta  \sim \beta + 2\pi \, , \qquad  \psi \sim  \psi  +3\pi\, \, , 
\end{equation}
thus suggesting that globally the space is an orbifold $\mathbb{S}^5/\mathbb{Z}_4$. However, the precise form of the $\mathbb{Z}_4$ action is not transparent from these considerations.

Next, we will show that the $\mathbb{Z}_4$ action is precisely the correct action inherited from the $\mathbb{C}^3/\mathbb{Z}_4$ orbifold singularity. 
We start  with three standard complex coordinates $(z_1,z_2,z_3)$ on $\mathbb{C}^3$ and consider the following change of coordinates
 \begin{align}
z_1  = R \sin\sigma \cos\tfrac{\theta}{2} e^{{\rm i}\left(-\frac{\beta+\phi}{2} + \frac{\psi}{3}\right)}\, , \qquad 
z_2  = R  \sin\sigma \sin\tfrac{\theta}{2} e^{{\rm i}\left(\frac{\phi-\beta}{2} +\frac{\psi}{3}\right)}\, , \qquad 
z_3  = R \cos\sigma e^{{\rm i}\frac{\psi}{3}}\, , 
\label{pasquale}
 \end{align}
where 
\begin{align}
|z_1|^2+|z_2|^2+|z_3|^2& =R^2
 \end{align} 
It can be checked that the metric induced at $R=1$,
\begin{align}
{\rm d s}^2_5 &= (|dz_1|^2+|dz_2|^2+|dz_3|^2) |_{R=1}
 \end{align} 
 coincides with (\ref{third5dmetric}), and more generally, the six-dimensional metric is the cone ${\rm  d s}^2_{\rm cone} =dR^2+R^2 \text{ds}^2_{X_5}$.
 
  To see that these are good coordinates on $\mathbb{C}^3$ we can also view it as $\mathbb{C}^3\simeq \mathbb{R}^6 = \mathbb{R}^2\oplus \mathbb{R}^2\oplus \mathbb{R}^2$, by defining 
 \begin{align}
z_1  =  \rho_1 e^{{\rm i}\varphi_1}\, , \qquad z_2  =  \rho_2 e^{{\rm i}\varphi_2}\, , \qquad z_3  =  \rho_3 e^{{\rm i}\varphi_3}
 \end{align} 
 so that the induced metric reads
 \begin{align}
{\rm  d s}^2_{\rm cone} &= |dz_1|^2+|dz_2|^2+|dz_3|^2 =d\rho_1^2 + \rho_1^2 d\varphi_1^2+d\rho_2^2 + \rho_2^2 d\varphi_2^2+d\rho_3^2 + \rho_3^2 d\varphi_3^2
  \end{align} 
  where for $\mathbb{C}^3$ the ranges of the coordinates are now   $\rho_i\in [0,+\infty)$ and $\varphi_i \sim \varphi_i + 2\pi$,  for $i=1,2,3$. 
  Defining $y_i=\tfrac{1}{2}\rho_i^2$, this gives the standard metric in symplectic-toric  coordinates, with  K\"ahler form
  \begin{align} 
  \mathbb{K}_{\mathbb{C}^3} = d y_1\wedge d\varphi_1+ d y_2\wedge d\varphi_2+d y_3\wedge d\varphi_3
    \end{align} 
  From this it is clear that  the $U(1)^3$ torus action on $\mathbb{C}^3$ 
     \begin{align} 
   (z_1,z_2,z_3) \to  (\lambda_1 z_1,\lambda_2z_2,\lambda_3 z_3)
    \end{align}
 with $|\lambda_i|=1$, $\lambda_i=e^{{\rm i}c_i}$ descends on the $\varphi_i$ coordinates to 
     \begin{align} 
     \label{pippo}
  (\varphi_1,\varphi_2,\varphi_3) \to  (\varphi_1+c_1,\varphi_2+c_2,\varphi_3+c_3)
    \end{align} 
 Notice that on $\mathbb{C}^3$  the periodicities of the two sets of angular coordinates are consistent with  the change of  coordinates
   \begin{align} 
  \varphi_1  = -\frac{\beta+\phi}{2} + \frac{\psi}{3}\, , \qquad 
\varphi_2  = \frac{-\beta+\phi}{2} + \frac{\psi}{3}\,  \,  , \qquad 
\varphi_3  = \frac{\psi}{3}\, , 
 \end{align} 
 with inverse
  \begin{align} 
\phi   = -\varphi_1+\varphi_2 , \qquad 
\beta = -\varphi_1-\varphi_2 +2\varphi_3\,  \,  , \qquad 
\psi  = 3\varphi_3\, , 
 \end{align} 
  as we have
    \begin{align} 
  (2\pi)^3 = \int d\varphi_1 d\varphi_2 d\varphi_3 = \frac{1}{6} \int d\phi d\beta d\psi = \frac{1}{6} (2\pi)(4\pi)(6\pi)
   \end{align}

Let us now reformulate the standard orbifold action of a discrete group $\Gamma\in SU(3)$ on $\mathbb{C}^3$ with the corresponding 
action on $\mathbb{S}^5$ in the above $(\phi, \beta,\psi)$ coordinates. We will restrict to $\Gamma=\mathbb{Z}_n$ for simplicity. 
  In the $(z_1,z_2,z_3)$ coordinates on $\mathbb{C}^3$, a $\mathbb{Z}_n$ orbifold action is defined by the identification
   \begin{align} 
   (z_1,z_2,z_3)\sim  (\omega_n^{a_1} z_1,\omega_n^{a_2} z_2,\omega_n^{a_3} z_3)
    \end{align} 
  where $\omega_n$ is a $n$-th root of unity. The requirement that $\mathbb{Z}_n\in SU(3)$ implies that 
\begin{equation}
 a_1+a_2+a_3=0\quad \mathrm{mod}\quad  n
 \end{equation}
  Using (\ref{pippo}),  the above  orbifold action  implies the following  identification in the $\varphi_i$ coordinates
    \begin{align} 
   (\varphi_1,\varphi_2,\varphi_3) \sim (\varphi_1+a_1\tfrac{2\pi}{n},\varphi_2+a_2\tfrac{2\pi}{n},\varphi_3+a_3\tfrac{2\pi}{n})
    \end{align} 
   and, equivalently, the following identification in the $(\phi, \beta,\psi)$ coordinates
      \begin{align} 
      \label{pluto}
(\phi, \beta,\psi) \sim (\phi+(-a_1+a_2)\tfrac{2\pi}{n},\beta+(-a_1-a_2+2a_3)\tfrac{2\pi}{n},\psi+3a_3\tfrac{2\pi}{n})
    \end{align}

   The simplest example is the $\mathbb{C}^3/\mathbb{Z}_3$ orbifold,  with  $\mathbb{Z}_3$ action on $\mathbb{C}^3$ given by  
   $$ (z_1,z_2,z_3) \sim ( e^ {{\rm i}\tfrac{2\pi}{3}} z_1,e^ {{\rm i}\tfrac{2\pi}{3}}z_2,e^ {{\rm i}\tfrac{2\pi}{3}} z_3),$$ which using (\ref{pluto}) corresponds simply to $   \psi \sim  \psi + 2\pi$.
   In this case, the metric (\ref{third5dmetric}), taking  $\phi \in [0,2\pi]$, $\beta \in [0,4\pi]$, is the metric  on $\mathbb{S}^5/\mathbb{Z}_3$. 
   This space can also be viewed as the unit circle bundle inside $\mathcal{O}_{\mathbb{P}^2}(-3)$, namely
   the total space of the canonical  line bundle over $\mathbb{P}^2$.

Let us now discuss our main example, the orbifold $\mathbb{C}^3/\mathbb{Z}_4$. In the table below we summarise the action of the three non-trivial elements of $\mathfrak{g} \in \mathbb{Z}_4$, including the identifications both in the $(\varphi_1,\varphi_2,\varphi_3)$  and the $(\phi, \beta,\psi)$ coordinates.
\begin{equation}\label{Z4actiononHopfcoordinates}
\begin{array}{|c|c|c|c|}
\hline
\mathfrak{g} : (z_1,z_2,z_3) & \{a_1,a_2,a_3 \} & (\varphi_1,\varphi_2,\varphi_3)\sim  &  (\phi, \beta,\psi) \sim  \\
\hline
({\rm i},{\rm i}, -1) &  \{1,1,2 \} & (\varphi_1+\tfrac{\pi}{2},\varphi_2+\tfrac{\pi}{2},\varphi_3+\pi)  & (\phi, \beta+\pi,\psi+3\pi)\\
(-1,-1, 1) &  \{2,2,0 \} & (\varphi_1+\pi,\varphi_2+\pi,\varphi_3) & (\phi, \beta+2\pi,\psi)\\
(-{\rm i},-{\rm i}, -1) &  \{3,3,2 \} & (\varphi_1+\tfrac{3\pi}{4},\varphi_2+\tfrac{3\pi}{4},\varphi_3+\pi)& (\phi, \beta+3\pi,\psi+3\pi)\\
\hline
  \end{array}
\end{equation}
As we see, in either of these two sets of angular coordinates the identifications are not diagonal. In the  coordinates $(\phi, \beta,\psi)$ the clearest identification is the action of (junior) element  $\{2,2,0 \}$, 
which implies  that the base space, with metric  in the first line of (\ref{third5dmetric}), is $\mathbb{P}^2/\mathbb{Z}_2$. 
 The action of the (junior) element  $\{1,1,2 \}$
means that as $\beta$ goes half way around its circle, the coordinate $\psi$ goes once around the $\psi$-circle, with period $3\pi$.
The action of the (senior) element  $\{3,3,2 \}$   is simply a consequence of the previous two. 

In order to  clarify the orbifold action on $S^5$, it is useful to adopt a set of angular coordinates in which the $\mathbb{Z}_4$  action 
is diagonal. It is then simple to verify that this is achieved precisely by the original coordinates $(\phi,\tau,\chi)$ defined in (\ref{iniz1}).
We summarise this diagonal action in the table below, where for convenience we defined $\gamma \equiv \tfrac{4}{3}\chi$.
\begin{equation}\label{Z4actiononpapercoordinates}
\begin{array}{|c|c|c|}
\hline
\mathfrak{g} : (z_1,z_2,z_3) & \{a_1,a_2,a_3 \} &  (\phi, \tau,\gamma) \sim  \\
\hline
({\rm i},{\rm i}, -1) &  \{1,1,2 \} & (\phi, \tau,\gamma+2\pi)\\
(-1,-1, 1) &  \{2,2,0 \} & (\phi, \tau+2\pi,\gamma)\\
(-{\rm i},-{\rm i}, -1) &  \{3,3,2 \} & (\phi, \tau+2\pi,\gamma+2\pi)\\
\hline
  \end{array}
\end{equation}
This shows that the indeed, the $\mathbb{Z}_4$ action on $\mathbb{C}^3$ induces the correct $\mathbb{Z}_4$ action on the
asymptotic metric on $\mathbb{S}^5$.

In order to further clarify the orbifold action on $\mathbb{S}^5$, it is convenient
to rewrite the metric  (\ref{third5dmetric}) in the form of a   circle fibration over a base space, 
that turns out to be precisely  $\mathbb{W} P [112]$. In particular, rearranging the terms in  (\ref{third5dmetric})  we find
 \begin{align}
 \text{ds}^2_{X_5}& =  \widetilde{\text{ds}}^2_{\mathbb{W} P [112]}+ \frac{1}{16}  (1+3\cos^2\sigma)    \left[  d\gamma + \frac{2\sin^2\sigma}{1+3\cos^2\sigma} (d\tau - \cos\theta d\phi)\right]^2 
\label{bingo}
 \end{align}
where\footnote{Interestingly, precisely this metric was found in \cite{Gauntlett:2004zh} as a limiting case of a more general  one-parameter family of 
smooth metrics on $\mathbb{F}_2$, in the context of AdS$_5$ solutions of eleven-dimensional supergravity. See eq.  (5.7) of this reference.} 
  \begin{align}
 \widetilde{\text{ds}}^2_{\mathbb{W} P [112]}&=   d\sigma^2 + \frac{1}{4}\sin^2\sigma  \left(d\theta^2 + \sin^2\theta d\phi^2 \right)   +  \frac{\sin^2 \sigma\cos^2\sigma}{1+3\cos^2\sigma}\left( d\tau- \cos\theta d\phi\right)^2\, ,
\label{oldwcp2metric}
\end{align}
which clearly displays the fact that $S^5/\mathbb{Z}_4$ arises as the total space of a circle fibration over $\mathbb{W} P [112]$, equipped with the metric (\ref{oldwcp2metric}).
 We decorated this metric with a tilde to distinguish it 
 from the different metric on $\mathbb{W} P [112]$, that we discuss in the main body of the paper, namely the metric  
 (\ref{divisor2}) induced on the exceptional divisor by the Ricci-flat metric  ({\ref{iniz11}). 
 Below we will rewrite the latter metric in different coordinates, to facilitate 
 the comparison with the metric in (\ref{oldwcp2metric}). 
 Let us discuss briefly how to see that the underlying (singular) variety to the metric defined in (\ref{oldwcp2metric}) is indeed  $\mathbb {W}P[112]$.
 With the ranges of coordinates and periodicities  $\sigma \in [0,\frac{\pi}{2}]$, $\theta \in [0,\pi]$, $\phi \in [0,2\pi]$, $\tau\in [0,2\pi]$ we see that near to $\sigma\approx 0$ the metric develops an
 $\mathbb{R}^4/\mathbb{Z}_2$ singularity (it is a cone over the Lens space $\mathbb{S}^3/\mathbb{Z}_2$), while near to $\sigma\approx \tfrac{\pi}{2}$, the space shrinks smoothly to $\mathbb{S}^2 \times \mathbb{R}^2$.
 Following a reasoning analogous to that in the main body of the paper, one can see that there exists only one non-trivial two-cycle
 \begin{equation}
 C_2\quad  \Leftrightarrow \quad \{ \theta = {\rm constant}, ~\phi = {\rm constant}\}
 \end{equation}
 while the other two-cycle of $\mathbb{F}_2$, that would be defined  by $C_1  \Leftrightarrow \{t=t_{max}=-\tfrac{3}{2}\sin^2\sigma_{max}\neq 0\}$ is shrunk to zero size in the above metric\footnote{Since we have not established whether the above metric is K\"ahler, 
 the simplest way to see this is probably to consider the explicit one-parameter familiy of metrics on $\mathbb{F}_2$ desingularising (\ref{oldwcp2metric}), presented in \cite{Gauntlett:2004zh}.}.

From the metric (\ref{bingo}) we now read off the connection one-form 
 \begin{align}
\widetilde {\cal A} & \equiv    \frac{2\sin^2\sigma}{1+3\cos^2\sigma} (d\tau - \cos\theta d\phi)
\label{bingobongo}
\end{align} 
whose associated first Chern class can be integrated on $C_2$ to give
 \begin{align}
\frac{1}{2\pi}\int_{C_2} d \widetilde  {\cal A} & = 2
\label{pimpa}
\end{align} 
showing that indeed this is a connection on the unit circle bundle inside the canonical bundle of $\mathbb{W} P [112]$.

We conclude this appendix by writing the metric on the exceptional divisor 
 (\ref{divisor2}) induced 
 by the Ricci-flat metric  ({\ref{iniz11})  in a form that makes more transparent the comparison with  the above discussion. Using (\ref{ttosigma}) we have that 
 \begin{align}
 \text{ds}^2_{\mathcal{ED}} &= \frac{9}{4} \left[  (1+\cos^2\sigma)  d\sigma^2 +\frac{1}{2}\sin^2\sigma  \left(d\theta^2 + \sin^2\theta d\phi^2 \right)  +  \frac{\sin^2 \sigma\cos^2\sigma}
 {1+\cos^2\sigma}\left( d\tau - \cos\theta d\phi\right)^2\right]
 \label{exceptwcp2}
\end{align}
from which the similarity with the metric (\ref{oldwcp2metric}) is apparent.
 For completeness, let us also display the behaviour of the  orthotoric metric (\ref{iniz11}) near to the exceptional divisor. Setting $s=-3-\rho^2$ in (\ref{iniz11}), for $\rho\to 0$ we have 
    \begin{align}
    \label{nomoremetrics}
 \text{ds}_{{\rm tot}  K_{{\mathbb WP}[112]}}^2& \stackrel{\rho \to 0}{\approx}   \text{ds}^2_{\mathcal{ED}} + \frac{3(1+\cos^2\sigma ) }{8}  \left[ d\rho^2    +  \rho^2  \left[ d \gamma + \frac{2\sin^2\sigma}{1+\cos^2\sigma}
 (d\tau - \cos\theta d\phi)\right]^2 \right]
 \end{align}
 where the angular variables $(\phi,\tau,\gamma\equiv \tfrac{4}{3}\chi)$ are precisely those defined in (\ref{iniz1}), which all have canonical $2\pi$-periodicities.
This shows that the metric  (\ref{iniz11})  is smooth in the neighborhood of the exceptional divisor 
 $\mathcal{ED}={\mathbb WP}[112]$, in particular locally it has the topology of ${\mathbb WP}[112]\times \mathbb{C}^2$.
The  connection one-form 
 \begin{align}
 {\cal A} & \equiv    \frac{2\sin^2\sigma}{1+\cos^2\sigma} (d\tau- \cos\theta d\phi)
\label{bingobongo}
\end{align} 
read off from (\ref{nomoremetrics}) has first Chern class again given by
 \begin{align}
\frac{1}{2\pi}\int_{C_2} d  {\cal A} & = 2
\label{pimpa}
\end{align} 
as it should be.  

To summarise,  in this appendix we have shown that the orbifold action of $\mathbb{Z}_4$ on $\mathbb{S}^5$, induced by the $\mathbb{C}^3/\mathbb{Z}_4$ quotient,  is not 
diagonal  in the canonical coordinates where the Sasaki-Einstein metric on  $\mathbb{S}^5$ can be viewed as a $U(1)$ fibration over $\mathbb{P}^2$ with its K\"ahler-Einstein metric.  
This action is  \emph{diagonalised}  precisely by   the coordinates
$(\phi,\tau,\tfrac{4}{3}\chi)$ used in the main part of the paper, and adapting the metric to these coordinates, it takes the form of a $U(1)$ fibration 
 over  ${\mathbb WP}[112]$, with the non-Einstein metric (\ref{oldwcp2metric}). 
This is precisely the unit circle bundle in the canonical line bundle over ${\mathbb WP}[112]$. The metric on the exceptional divisor of the partial resolution, 
induced from the orthotoric Ricci-flat metric, is a similar, but manifestly \emph{different}  non-Einstein metric (\ref{exceptwcp2}).


\bigskip
\noindent{\bf Statement about conflict of interest}: On behalf of all authors, the corresponding author states that there is no conflict of interest.

 \end{document}